\documentclass[11pt]{article}

\usepackage{amssymb}
\usepackage{amsmath}
\usepackage[dvips]{epsfig}
\usepackage[small]{caption}
\usepackage{graphicx}
\usepackage{mathrsfs}
\usepackage{wrapfig}
\usepackage{caption}
\usepackage{subcaption}
\usepackage{float}
\usepackage[colorlinks=true,linkcolor=blue,citecolor=red]{hyperref}

\newtheorem{Lemma}{Lemma}[section]
\newtheorem{Theorem}{Theorem}
\newtheorem{Proposition}[Lemma]{Proposition}

\newtheorem{Remark}[Lemma]{Remark}

\newtheorem{Hypothesis}[Lemma]{Hypothesis}

\newenvironment{Proof}[1][Proof]%
 {\begin{trivlist} \item[]{\textbf{#1.} }}%
{\hspace*{\fill}$\rule{.4\baselineskip}{.4\baselineskip}$\end{trivlist}}

 {\begin{trivlist}\item[]\textbf{Acknowledgments }}{\end{trivlist}}

\makeatletter\@addtoreset{figure}{section}\makeatother

\makeatletter \@addtoreset{equation}{section} \makeatother

\newcommand{\R}{\mathbb{R}}
\newcommand{\C}{\mathbb{C}}

\newcommand{\Z}{\mathbb{Z}}

\def\Re{\mathop{\mathrm{Re}}}
\def\Rg{\mathop{\mathrm{Rg}}}

\def\cok{\mathop{\mathrm{CoKer}}}
\def\dim{\mathop{\mathrm{dim}}}

\def\un{\mathrm{n_-}}

\newcommand{\bbF}{\mathbb{F}}
\newcommand{\bbL}{\mathbb{L}}
\newcommand{\qprime}{{\prime\prime\prime\prime}}
\newcommand{\tprime}{{\prime\prime\prime}}
\newcommand{\pprime}{{\prime\prime}}
\newcommand{\caL}{\mathcal{L}}
\newcommand{\caO}{\mathcal{O}}
\def\caE{\mathcal{E}}

\newcommand{\rmd}{\mathrm{d}}
\newcommand{\rme}{\mathrm{e}}
\newcommand{\rmi}{\mathrm{i}}
\newcommand{\id}{\mathrm{\,Id}\,}
\newcommand{\ad}{\mathrm{ad}}

\newcommand{\sgn}{\mathrm{sgn}}
\renewcommand{\span}{\mathrm{\,span}}
\renewcommand{\leq}{\leqslant}
\renewcommand{\geq}{\geqslant}
\def\beq{\begin{equation}}
\def\eeq{\end{equation}}
\def\cF{{\cal F}}
\def\eps{\epsilon}
\newcommand{\Frac}[2]{\frac{\textstyle #1}
                           {\textstyle #2}}

\newcommand{\rmnum}[1]{\romannumeral #1}
\newcommand{\Rmnum}[1]{\expandafter\@slowromancap\romannumeral #1@}

\newfam\bifam
\font\tenbi=cmmib10 scaled \magstep1 \font\sevenbi=cmmib10 at 11pt
\font\fivebi=cmmib10 at 6pt \textfont\bifam = \tenbi
\scriptfont\bifam = \sevenbi \scriptscriptfont\bifam= \fivebi

\begin{document}

\thispagestyle{empty}

\title{\bf Existence of pearled patterns in the planar Functionalized Cahn-Hilliard equation}

\author{Keith Promislow and Qiliang Wu\\
\textit{\small Department of Mathematics, Michigan State University}\\
\textit{\small 619 Red Cedar Road East Lansing, MI 48824}}

\date{\small \today}

\maketitle

\begin{abstract}
\noindent The functionalized Cahn-Hilliard (FCH) equation supports planar and circular bilayer interfaces as equilibria
which may lose their stability through the pearling bifurcation: a periodic, high-frequency, in-plane modulation of the bilayer thickness. 
In two spatial dimensions we employ spatial dynamics and a center manifold reduction to reduce the FCH equation to an 8th order ODE system.
A normal form analysis and a fixed-point-theorem argument show that the reduced system admits a degenerate 1:1 resonant normal form, from
which we deduce that the onset of the pearling bifurcation coincides with the creation of a two-parameter family of pearled equilibria which
are periodic in the in-plane direction and exponentially localized in the transverse direction. \end{abstract}


 \hrule
 {\small
%
%
%
  }


\section{ The Functionalized Cahn-Hilliard equation}
Amphiphilic materials are typically small molecules which contain both hydrophilic and hydrophobic components. This class of materials includes
surfactants, lipids, and block copolymers. Their propensity to spontaneously assemble network morphologies has drawn scientific attention for more than a 
century, \cite{amphiphilic2000}.  While amphiphilic materials are ubiquitous in organic settings, where lipid bilayers form cell membranes and many organelles,
their widespread use as charge separators in energy conversion devices is more recent.  Network morphologies must be distinguished from single layer interfaces that are
typical of binary metals and other purely hydrophobic blends. While single layer interfaces separate a phase $A$ from a phase $B$, network morphologies are
comprised of thin regions of a phase $B$ which interpenetrate, and typically percolate through, a domain dominated by phase $A$. The Cahn-Hilliard
free energy, proposed in 1958, \cite{CH}, has been very successfully employed as a model of single layer morphology in hydrophobic blends, and 
its gradient flows accurately describe their  evolution. Models of amphiphilic mixtures, such as  \cite{TS-87} and \cite{GS-90}, have been proposed.
The functionalized Cahn-Hilliard free energy; see \cite{KeithBrian-09, KeithNiGr-11, keithdai_2013}, 
is a special case of these earlier models that supports stable network morphologies including co-dimension one bilayers and co-dimension 
two pores as well as pearled morphologies and defects such as end-caps and junctions. Rigorous results for the FCH free energy include the existence of bilayer structures, \cite{doelmanpromislow_2013}, and an analysis of their bifurcation structure, \cite{keithgreg_2014}, in particular the pearling bifurcation  
which initiates  changes in the co-dimension of the underlying morphology, and is commonly observed in  amphiphilic polymer blends; see \cite{szostak,hayward}. The goal of this paper is to rigorously establish the existence of pearled bilayers, as modulations to stationary bilayers, 
in the planar FCH equation.

Amphiphillic mixtures, such as emulsions formed by adding a minority fraction of an oil and soap mixture
to water, form network morphologies due to the tendency of the surfactant phase, e.g. soap, to enhance the formation of interfaces.
To model the network formation, the authors of \cite{TS-87} and \cite{GS-90} were motivated by small-angle X-ray scattering (SAXS) data to include 
a higher-order term in the usual Cahn-Hilliard expansion for the free energy. Viewing the mixture as a binary phase, where $u\in H^2(\Omega)$
denotes the volume fraction of surfactant contained within the bounded material domain $\Omega\subset\R^3$, they proposed a free energy of the form
 \vskip -0.2in
\beq
 \cF(u) := \int_\Omega f(u) + \eps^2A(u)|\nabla u|^2 + \eps^2B(u) \Delta u+  C(u) (\eps^2\Delta u)^2\,dx,
\eeq
where for well-posedness $C>0$ and the dimensionless parameter $\varepsilon \ll 1$ dictates the ratio of the interfacial width to a
characteristic size of $\Omega$. Assuming zero-flux boundary conditions, integration by parts on the $A(u)$ term permits a re-writing of 
the energy in the completed-square form
\beq
 \cF(u) = \int_\Omega  C(u)\left(\eps^2\Delta u- \Frac{\overline{A}-B}{2C} \right)^2 +
 f(u) -\Frac{(\overline{A}-B)^2}{4C(u)}\, dx,
\eeq
where $\overline{A}$ is a primitive of $A$.  To simplify the form we replace $C(u)$ with $\frac12$, relabel the potential within the
squared term by $W'(u)$, and scale the potential outside the squared term as $\delta P(u)$ with $\delta\ll1$, yielding
\beq
\cF(u) = \int_\Omega \frac{1}{2}\left( \eps^2 \Delta u-W'(u)\right)^2 + \delta P(u)\, dx.
\eeq
The first term is the square of the variational derivative of a Cahn-Hilliard type free energy, and 
the strongly degenerate case $\delta=0$, has the special property that its global minimizers are precisely the 
{\sl critical points} of the corresponding Cahn-Hilliard energy. A variant of this case was proposed
as a target for $\Gamma-$convergence analysis by De Giorgi; see \cite{DeGiorgi}.

The strong functionalized Cahn-Hilliard free energy corresponds to the distinguished limit $\delta=\varepsilon$, 
a choice of potential $P$ which incorporates the functionalization parameters $\eta_1>0$, $\eta_2\in \R$ in the form 
\begin{equation}
\label{FCHE}
 \mathcal{F}(u)=\int_{\Omega}\frac{1}{2}\left(\varepsilon^2 \Delta u- W^\prime(u)\right)^2 - \varepsilon \left(\eta_1 \varepsilon^2|\nabla u|^2 + \eta_2 W(u) \right)\rmd x,
\end{equation}
and require  the $C^\infty$-smooth potential $W: \R\rightarrow\R$   to be a double well potential with two minima at $u=-1$ and $u=m>0$ and one local maximum
at $u=0$. The minima have unequal depths, normalized so that $W(-1)=0>W(m)$ and the well is non-degenerate in the sense that  
$\mu_-:=W^{\prime\prime}(-1)>0$,  $\mu_+:=W^{\prime\prime}(m)>0$, and $\mu_0:=W^{\prime\prime}(0)<0$.
With these assumptions $u=-1$ is associated to a bulk solvent phase, while the value of $u+1>0$ is proportional to the density of the amphiphilic phase.

 The strong FCH equation is the $H^{-1}$ gradient flow of the FCH energy \eqref{FCHE}, which takes the form
\begin{equation}
 \label{e:FCH}
 u_t=\Delta\frac{\delta \mathcal{F}}{\delta u}=\Delta\left((\varepsilon^2\Delta - W^{\prime\prime}(u) + \varepsilon \eta_1)
 (\varepsilon^2\Delta u - W^{\prime}(u))+\varepsilon\eta_d W^\prime(u)\right),
\end{equation}
where $\eta_d:=\eta_1-\eta_2$. The gradient flow is mass-preserving when subject to zero-flux boundary conditions; see \cite{doelmanpromislow_2013} for details. 
We focus on the stationary strong-FCH equation which takes the form
 \begin{equation}
 \label{e:sFCH}
 (\varepsilon^2\Delta - W^{\prime\prime}(u) + \varepsilon \eta_1)
 (\varepsilon^2\Delta u - W^{\prime}(u))+\varepsilon\eta_d W^\prime(u)=\varepsilon \gamma,
 \end{equation}
 subject to zero-flux boundary conditions. The constant $\gamma$ can be thought of as a Lagrange multiplier arising from mass conservation.

The FCH equation is known to support families of bilayer solutions, \cite{doelmanpromislow_2013}, which can be unstable to either pearling or meandering bifurcations.
Pearling refers to periodic modulations of the thickness of the bilayer, while the meander modes are associated with the curvature driven motion of the 
underlying bilayer interface. In this work, we provide a fully rigorous proof of the existence of spatially periodic patterns which arise after the onset of the pearling 
bifurcation. We  restrict our attention to planar domains $\Omega\subseteq\R^2$, proving the major existence results
in the spatially extended case $\Omega=\R^2$. The construction of a bilayer morphology requires a choice of a smooth,
closed, co-dimension one interface $\Gamma\subset\Omega$ that is far from self intersection. We address
two simple choices of interface: the extended flat bilayer, corresponding to $\Gamma_f=\{(s,0)\,\bigl |\, s\in\R\}$, and the circular bilayer of radius $R_0>0$,
corresponding to $\Gamma_{R_0}:=\{(R_0\cos{\theta}, R_0\sin{\theta})\, \big|\, \theta\in[0,2\pi)\}.$ Our construction applies spatial dynamics techniques, 
 a center-manifold-reduction argument, and a normal form transformation to the stationary, strong-FCH equation, yielding an 8th order ODE system, 
 which weakly couples the four dimensional pearling  subspace and the four dimensional meander subspace. To prove the existence, we  
 restrict to the pearling subspace, yielding a four-dimensional reduced system, called the pearling normal form (PNF), \eqref{e:PNF},
 \begin{equation*}
\begin{cases}
\dot{C_1}&=\rmi (1+\omega_1\varepsilon) C_1 + C_2+\rmi C_1\big[\alpha_7C_1\bar{C}_1+\alpha_8\rmi(C_1\bar{C}_2-\bar{C}_1C_2)\big],\\
\dot{C_2}&=\rmi (1+\omega_1\varepsilon) C_2 +\rmi C_2\big[\alpha_7C_1\bar{C}_1+\alpha_8\rmi(C_1\bar{C}_2-\bar{C}_1C_2)\big]
+ C_1\left[-\alpha_0\varepsilon+\rmi \alpha_2(C_1\bar{C_2}-\bar{C_1}C_2)\right],
\end{cases}
\end{equation*}
where $C_1$, $C_2\in\C$, the constants $\omega_1, \alpha_j\in\R$, and the conjugate equations are omitted. It is at this level that the structure of the 
pearling bifurcation is  made clear: the PNF admits a \textit{degenerate $1:1$ resonance}, related to the $1:1$ resonances extensively  investigated in \cite{imd_1989, iooper_1993, harioo}. As in the $1:1$ resonance case, the PNF has two first integrals
\[
K:=\frac{\rmi}{2}(C_1\bar{C}_2-\bar{C}_1C_2), \quad H:=|C_2|^2+\left(-\alpha_0\varepsilon+2\alpha_2K\right)|C_1|^2.
\]
Imposing consistency conditions to the solutions of the PNF slaves $H$ to the scaled parameter $\kappa:=\varepsilon^{-3/2} K$, which remains as a free parameter in the
construction of the pearled solutions. More importantly, the parameter $\alpha_0$ in the PNF, given in \eqref{def-alpha0}, is precisely the critical bifurcation parameter whose sign  characterizes the onset of the pearling bifurcation. For $\alpha_0>0$ we characterize the pearled solutions of the PNF and establish their existence in
the full system through a persistence argument. While the persistence argument is based upon \cite{iooper_1993}, the analysis in this case is more delicate as
the degeneracy corresponds to a distinct singularity requiring different scalings. Moreover the coupling between the pearling modes and the meander modes 
requires the analysis of an eight dimensional problem. In the remainder of this section we make a rigorous statement of these results.

\subsection{Pearling of Extended Flat Bilayers} 
 The existence of a one-dimensional family of flat bilayer solutions, $u_h$, 
parameterized by the Lagrange multiplier, $\gamma$, was established in \cite{doelmanpromislow_2013}. 
Their construction is based upon new coordinates, corresponding to the $\varepsilon$-scaled distance $r$ to $\Gamma_f$ and a tangential variable $\tau$
for which the  Laplacian takes the form
 \begin{equation}
 \label{e:slap}
 \varepsilon^2\Delta=\partial_r^2+\varepsilon^2\partial_\tau^2,
\end{equation}
and the stationary equation (\ref{e:sFCH}) is rewritten as
\begin{equation}
 \label{e:2sFCH}
 \left(\partial_r^2-W^{\prime\prime}(u)+\varepsilon^2\partial_\tau^2+\varepsilon\eta_1\right)
 \left(\partial_r^2u-W^\prime(u)+\varepsilon^2\partial_\tau^2u\right)+\varepsilon\eta_d W^\prime(u)=\varepsilon\gamma.
\end{equation}
For the flat interface, the bilayer profile is independent of the tangential variable, $\tau$, and hence is captured as the first component of a homoclinic solution of the 
$4$-th order extended flat-bilayer ODE system in $r\in\R$,
 \begin{equation}
 \label{e:rODE}
  \begin{cases}
  \partial_r u=p, \\
  \partial_r p=W^\prime(u)+\varepsilon v, \\
  \partial_r v=q, \\
  \partial_r q=W^{\prime\prime}(u)v+\left(\gamma-\eta_dW^\prime(u)\right)-\varepsilon\eta_1v,
  \end{cases}
 \end{equation}
 For sufficiently small $\varepsilon$, this extended flat-bilayer ODE system \eqref{e:rODE} contains $3$ critical points, among which we consider the one with leading order
$(-1,0,-\frac{\gamma}{\mu_-},0)$, which we denote as 
\[
P_-(\varepsilon)=\left(u_-(\varepsilon),0,v_-(\varepsilon),0\right).
\]
Indeed, via \eqref{e:rODE}, it is straightforward to see that the parameter $\gamma$ relates linearly, at leading order, to the far-field density of amphiphilic
material,  $1+u_-(\varepsilon)$, via the expansion
\[
1+ u_-(\varepsilon;\gamma)=\frac{\gamma}{\mu_-^2}\varepsilon+\caO(\varepsilon^2).
\]

In \cite{doelmanpromislow_2013} the existence of the flat homoclinic solution $U_h=(u_h,p_h, v_h,q_h)^T$ is established for $\varepsilon>0$ sufficiently small, 
but independent of $\eta_1, \eta_2,$ and $\gamma$. The construction follows by perturbation off of the $\varepsilon=0$ case, in which case the first component $u_0$ is the solution
of the two-dimensional ODE
\begin{equation}
\label{2d-ODE}
\partial_r^2 u_0 = W'(u_0),
\end{equation}
which is homoclinic to $u_-(0)$. The linearization of (\ref{2d-ODE}) about $u_0$, yields the operator
\begin{equation}
\label{cL_0-def}
 \mathcal{L}_0:=\partial_r^2-W''(u_0),
\end{equation}
which, acting on $L^2(\R)$, has a single positive eigenvalue, $\lambda_0>0$, and a zero eigenvalue, $\lambda_1=0$, with the remainder of the spectrum strictly negative. Denoting
the associated eigenfunctions by $\psi_0$ and $\psi_1$ and introducing, $v_0\in L^\infty(\R)$, the unique, even solution of
\begin{equation}
 v_0 = \gamma \mathcal{L}_0^{-1} 1 -\eta_d \mathcal{L}_0^{-1}W'(u_0),
 \end{equation}
the pearling bifurcation of the bilayer $u_h$  is characterized in terms of the functionalization parameters $\eta_1$ and $\eta_2$ via the sign of the quantity
\begin{equation}
\label{def-alpha0}
\alpha_0=\frac{1}{4\lambda_0^2}\int_{\R}\left(W^\tprime(u_0)v_0-\eta_dW^\pprime(u_0)\right)\psi_0^2\rmd r=  \alpha_{01} \gamma  - \alpha_{02} \eta_d, \end{equation}
where the constants
\begin{equation}
\begin{aligned}
  \alpha_{01} &= \frac{1}{4\lambda_0^2} \int_{\R} W^\tprime(u_0)(\mathcal{L}^{-1} 1)\psi_0^2\rmd r, \\
\alpha_{02} & :=     \int_{\R} \left((\mathcal{L}_0^{-1}W'(u_0)+  W^\pprime(u_0)\right)\psi_0^2\rmd r,
\end{aligned}
\end{equation}
depend only upon the shape of the double well potential, $W$.

Our main result for flat bilayers establishes that a one parameter family of pearled solutions of (\ref{e:2sFCH}) generically bifurcates out
of each stationary flat bilayer for  $\alpha_0>0$.

\begin{Theorem}[existence of extended pearled flat bilayers]
\label{t:main}
 Fix $\eta_1, \eta_2, \gamma \in\R$. Assume that $W$ is a non-degenerate double well potential and that
 $\alpha_0$ defined in (\ref{def-alpha0}) is strictly positive and
\begin{equation}
\label{def-beta0}
\beta_0:=\frac{1}{4\lambda_0^2}\int_\R\left(W^\tprime(u_0)v_0-\eta_dW^\pprime(u_0)\right)\psi_1^2\rmd r\neq 0,\\
\end{equation} 
Then there exist positive constants $\varepsilon_0>0$ and $\kappa_0>0$ such that, for any $\varepsilon\in(0,\varepsilon_0]$, up to translation, 
the extended stationary strong-FCH \eqref{e:2sFCH} 
 admits a smooth one-parameter family of extended pearled solutions, 
 $u_\mathrm{p}(\tau,r;\sqrt[4]{\varepsilon},\sqrt{|\kappa|})$ with period $T_\mathrm{p}(\sqrt[4]{\varepsilon},\sqrt{|\kappa|})$, parameterized by $\kappa\in[-\kappa_0,\kappa_0]$. More specifically, $u_{\mathrm{p}}$ and $T_{\mathrm{p}}$ are smooth with respect to their arguments within the domains expect at $\kappa=0$.  The extended pearled solution  $u_\mathrm{p}$ admits the asymptotic form
 \begin{equation}
 u_\mathrm{p}(\tau, r)=u_h(r)+2\frac{\sqrt{\varepsilon|\kappa|}}{\sqrt[4]{\alpha_0}}\cos\left(\frac{2\pi}{T_\mathrm{p}} \tau\right)\psi_0(r)+
 \caO\left(\varepsilon(\sqrt{\varepsilon}+\sqrt{|\kappa|})\right),
 \end{equation}
where the error is measured in the $L^\infty(\R^2)$-norm and
\begin{equation}
T_\mathrm{p}=\frac{2\pi\varepsilon}{\sqrt{\lambda_0}}\left[1-\sqrt{\alpha_0\varepsilon}+\caO\left(\varepsilon(1+\sqrt{|\kappa|})\right)\right] .
\end{equation}
Moreover,
 the far-field limit of the extended pearled solution is 
 \begin{equation}
  \lim_{r\rightarrow \infty}u_\mathrm{p}(\tau,r)=\lim_{r\rightarrow \infty}u_h(r)=u_-(\varepsilon).
  \end{equation}
 \end{Theorem}

\subsection{Pearling of extended Circular Bilayers}
For a circular co-dimension one interface $\Gamma_{R_0}$ we take the tangential coordinate $s$ to represent the direction with constant curvature $k=-R_0$, and rescale the corresponding independent variable
as $\theta=s/R_0$ which lies in $[0,2\pi]$. 
The Laplacian admits the expression
\begin{equation}
 \label{e:plap}
\varepsilon^2\Delta=\partial_r^2+\frac{\varepsilon}{R_0+\varepsilon r}\partial_r+\frac{\varepsilon^2}{(R_0+\varepsilon r)^2}\partial_\theta^2,
\end{equation}
and the stationary strong-FCH \eqref{e:sFCH} in $(r, \theta)$ takes the form
\begin{equation}
 \label{e:2sFCH2}
 \Big(\partial_r^2-W^{\prime\prime}(u)+\frac{\varepsilon\partial_r}{R_0+\varepsilon r}+\frac{\varepsilon^2\partial_\theta^2}{(R_0+\varepsilon r)^2}+\varepsilon\eta_1\Big)
 \Big(\partial_r^2u-W^\prime(u)+\frac{\varepsilon\partial_ru}{R_0+\varepsilon r}+\frac{\varepsilon^2\partial_\theta^2u}{(R_0+\varepsilon r)^2}\Big)+\varepsilon\eta_d W^\prime(u)=\varepsilon\gamma.
\end{equation}
 Suppressing the tangential variable $\theta$, the stationary strong-FCH \eqref{e:2sFCH2} reduces to the extended circular-bilayer ODE system in $r\in\R$, 
 \begin{equation}
 \label{e:rODE2}
  \begin{cases}
 \partial_r u=p, \\
  \partial_r p=W^\prime(u)+\varepsilon v, \\
  \partial_r v=q, \\
  \partial_r q=W^{\prime\prime}(u)v+[\gamma_1-\eta_dW^\prime(u)]+\varepsilon[\gamma_2-\frac{2}{R_0}q+\frac{1}{R_0^2}W^\prime(u)-\eta_1v-
      \frac{1}{R_0}\eta_1p]+\caO(\varepsilon^2),
  \end{cases}
 \end{equation}
 where $\gamma$ has been expanded as,
 \[\gamma=\gamma_1+\varepsilon\gamma_2+\caO(\varepsilon^2).\]
Like the flat-bilayer system, the extended circular-bilayer ODE system \eqref{e:rODE2} possesses $3$ critical points, of which we single out
the critical point 
\[P_-(\varepsilon)=\left(u_-(\varepsilon),0,v_-(\varepsilon),0\right),\]
which satisfies $P_-(\varepsilon)\to (-1,0,-\frac{\gamma_1}{\mu_-},0)$, as $\varepsilon\to0.$   In \cite{doelmanpromislow_2013}, it was shown
that for fixed $\eta_1, \eta_2$ and $R_0>0$ there exists a unique function 
$\gamma_h=\gamma_1+\caO(\varepsilon)$ for which
 \begin{equation}
 \label{def-gamma_1}
  \gamma_1=(\eta_d-2\eta_1)\frac{\int_\R (u^{\prime }_0 )^2\rmd r}{2\int_\R(u_0+1)\rmd r},
 \end{equation}
such that for the choice $\gamma=\gamma_h(\varepsilon)$ there exists a nontrivial orbit of (\ref{e:rODE2}) which is homoclinic to $P_-(\varepsilon).$
\begin{Remark}
The parameter $\gamma$ is free for flat bilayers while it is prescribed for circular bilayers because
the flat-bilayer ODE system \eqref{e:rODE} is Hamiltonian while the circular-bilayer ODE system \eqref{e:rODE2} is not. 

\end{Remark}

\begin{figure}[H]
 \centering
                \includegraphics[width=0.6\textwidth]{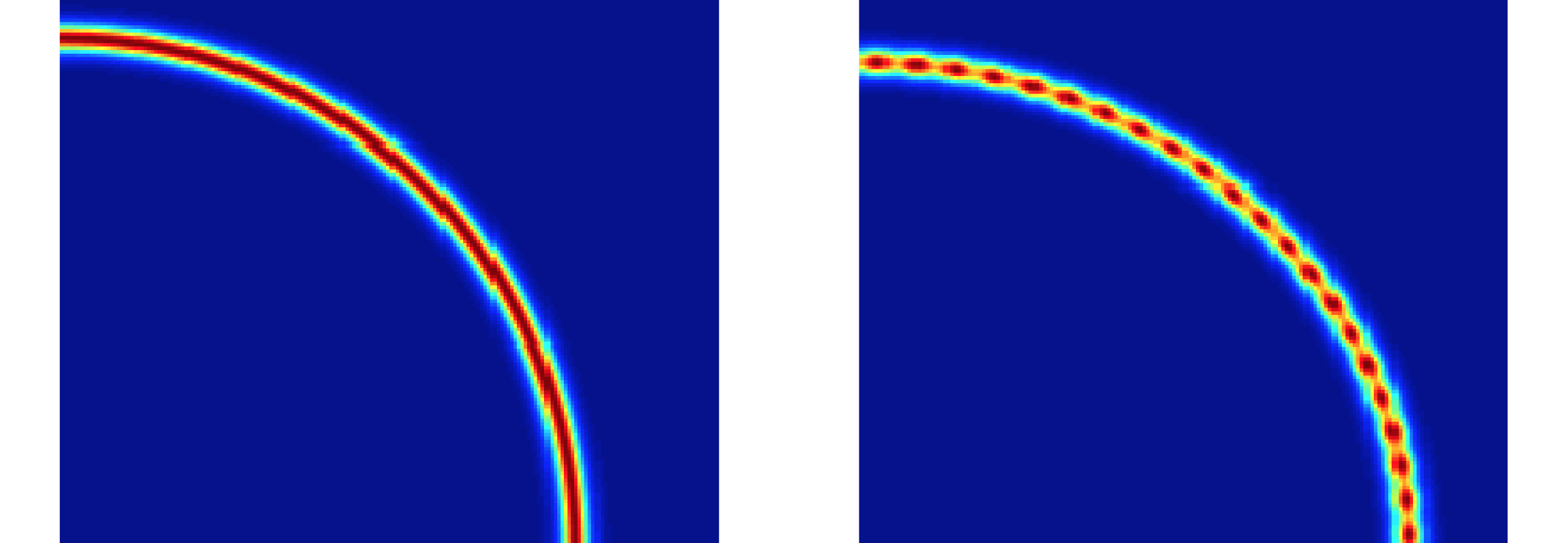}
\caption{ Quarter-plane views of equilibrium of the strong FCH equation \eqref{e:FCH} corresponding to radially symmetric bilayer initialdata with $\varepsilon=0.1$  
and double well potential $W$ as given in Section 5 of \cite{doelmanpromislow_2013}. (left) For $\eta_1= 1$ and $\eta_2= 2$, we have $\alpha_0<0$, and the $t=3000$ 
evolution is a circular bilayer equilibrium.  (right)  For $\eta_1=2$ and $\eta_2=2$, we have $\alpha_0>0$, and the $t=500$ evolution of the initial data yields a circular pearled bilayer.}
    \label{f:cirpearling}
 \end{figure}

Our main result for circular bilayers provides the existence of discrete families of one-parameter, pearled, bilayer solutions of the stationary strong-FCH equation (\ref{e:2sFCH2}); see Figure \ref{f:cirpearling}. Both their radii $R_{0,n}=R_{0,n}(\varepsilon, \kappa)$ and pearling amplitudes are parameterized
by the value of the scaled first-integral $\kappa$ of the Pearling Normal Form equation.
 
\begin{Theorem}[existence of extended pearled circular bilayers]
\label{t:main2}
 Fix $\eta_1, \eta_2\in\R$ and  $R_->0$. Assume that $W$ is a non-degenerate double well potential and that $\alpha_0$ and $\beta_0$, defined in 
 (\ref{def-alpha0}) and (\ref{def-beta0}) respectively, satisfy $\alpha_0>0$, $\beta_0\neq 0$. 
 Then there exist constants $\varepsilon_0, \kappa_0>0$ and $\un>0$ such that, 
 for all $\varepsilon\in(0,\varepsilon_0]$ and each $\mathrm{n}\in \Z^+\cap[\frac{\un}{\varepsilon},+\infty)$, 
the stationary, strong-FCH equation \eqref{e:2sFCH2} in the infinite strip $(\theta,r)\in(\R/2\pi\Z)\times\R$, subject to the choice $\gamma=\gamma_h(\eps)$,  with $\gamma_h$ defined by (\ref{def-gamma_1}),  
admits, up to translation, a finite family of one-parameter pearled solutions
 $u_{\mathrm{p,n}}(\theta,r;\sqrt[4]{\varepsilon},\sqrt{|\kappa|})$ with period $\frac{2\pi}{n}$ and 
 radius $R_{\mathrm{0,n}}(\sqrt[4]{\varepsilon},\sqrt{|\kappa|})\geq R_-$. 
 Each solution is parameterized by $\kappa\in[-\kappa_0,\kappa_0]$, and 
is smooth with respect to its arguments except at $\kappa=0$.  The extended pearled solution  $u_{\mathrm{p,n}}$ admits the asymptotic form
 \begin{equation}
 \label{circ-bilayer}
 u_{\mathrm{p,n}}(\theta, r;\sqrt[4]{\varepsilon},\sqrt{|\kappa|})=u_h(r)+2\frac{\sqrt{\varepsilon|\kappa|}}{\sqrt[4]{\alpha_0}}\cos(\mathrm{n}\theta)\psi_0(r)+\caO\left(\varepsilon(\sqrt{\varepsilon}+\sqrt{|\kappa|})\right),
 \end{equation}
where the radius of the circular bilayer
\begin{equation}
R_{\mathrm{0,n}} = \Frac{\mathrm{n}\varepsilon}{\sqrt{\lambda_0}} \left[1-\sqrt{\alpha_0\varepsilon}+\caO\left(\varepsilon(1+\sqrt{|\kappa|})\right)\right].
\end{equation}
depends only weakly upon $\kappa$. The far-field limit of the extended pearled solution 
 \begin{equation}
  \lim_{r\rightarrow \infty}u_{\mathrm{p,n}}(\theta,r)=\lim_{r\rightarrow \infty}u_h(r)=u_-(\varepsilon),
  \end{equation}
 is independent of $\mathrm{n}$.
 \end{Theorem} 
\begin{Remark}
The number $n$ can be interpreted as the number of ``beads'' within a pearled circular bilayer. The size of each bead--the periodicity in the physical variables-- is
\[
T_{\mathrm{p},n}:=\frac{2\pi R_{0,n}}{n}=\frac{2\pi\varepsilon}{\sqrt{\lambda_0}}\left[1-\sqrt{\alpha_0\varepsilon}+\caO\left(\varepsilon(1+\sqrt{|\kappa|})\right)\right],
\]
depends only weakly upon $\kappa$, at order $\caO(\varepsilon^2\sqrt{|\kappa|})$, while the leading order amplitude of each bead, 
\begin{equation}
\label{Ap-def}
A_\mathrm{p}:=2\frac{\sqrt{\varepsilon|\kappa|}}{\sqrt[4]{\alpha_0}},
\end{equation}
scales with $(\sqrt{\varepsilon|\kappa|})$.
\end{Remark}

For both the flat and circular interfaces, the form of the amplitude of the pearled pattern suggests a divergence as 
$\alpha_0\to0^+$, however this is an anomaly arising from the degeneracy of the $1:1$ resonance in the PNF system, \eqref{e:PNF}. Indeed an analysis
of Lemma \ref{l:29} shows that a necessary condition for the existence of periodic patterns is
\begin{equation}\label{e:super}
\sqrt{\varepsilon_0}\kappa_0< \frac{\alpha_0}{2|\alpha_2|},
\end{equation}
from which we deduce that the pearling bifurcation, while degenerate, retains some supercritical characteristics.

\begin{Proposition}[super-criticality of pearled bilayers]
In addition to the assumptions of either Theorem~\ref{t:main} or \ref{t:main2}, assume that
 $\alpha_2$, defined in \eqref{e:alphas}, satisfies  $\alpha_2\neq 0$. Fix
$\varepsilon\in(0,\varepsilon_0)$ and tune $\eta_1$ and $\eta_2$ so that $\alpha_0$ goes to $0$; then,
under this limit, the pearling amplitude, defined in \eqref{Ap-def}, satisfies 
\[
\lim\limits_{\alpha_0\to0}\sup\limits_{\kappa\in[-\kappa_0,\kappa_0]}\Frac{A_p(\kappa)}{\sqrt[4]{\alpha_0}}\leq C,
\]
for some constant $C>0.$ 
\end{Proposition}

\subsection{Pearling and Degeneracy in Bounded Domains}

The existence results for both bilayers and pearled bilayers  naturally extend to a bounded domain, $\Omega\subset \R^2$ 
so long as the domain possesses the same symmetry as the bilayer interface. Indeed,  for typical homogeneous boundary conditions, 
such as discussed in \cite{PZ-2013}, and for a bilayer interface $\Gamma$ that is an $\caO(1)$ distance from 
$\partial\Omega$ in the unscaled coordinates, then the exponential decay of the extended pearled patterns in $r$ leads to an $\caO(\varepsilon^{-1})$ exponential decay
in the unscaled coordinates, and a standard matching argument; such as in \cite{PB_2011}, permits an extension of the existence result. This
is particularly relevant for the circular bilayers within a concentric circular domain. The adaptation of the extended flat bilayer to a flat bilayer
within a rectangular domain subject to periodic boundary conditions is trivial so long as the flat interface intersects the domain boundary at a right angle; see Figure \ref{f:domain} for an illustration.
The construction of the associated pearled solutions requires a tuning of the periodicity of the pearled pattern, as in the case of
the circular bilayer. 
  
For the gradient flow (\ref{e:FCH}), the total mass $\int_\Omega u(x)\mathrm{d}x$ is conserved under time evolution, and as such it is natural to search
for equilibria with prescribed total mass. For circular bilayers; see Figure \ref{f:domain}, the far-field value of $u$ is prescribed, 
and the mass of a circular bilayer is an increasing function of the radius $R_0$. Moreover the mass is independent of the pearling correction, 
at least to leading order, thus the total mass of the circular bilayer $u_{p,n}$ in \eqref{circ-bilayer} increases monotonically with its
radius $R_{0,n}$; however the admissible  radii 
$$\left \{ R_{\mathrm{0,n}}(\kappa) \, \Bigl|\, \mathrm{n}\in \Z^+\cap[\frac{\un}{\varepsilon},+\infty), \, \kappa\in[-\kappa_0,\kappa_0]\right\}.$$
depend only weakly upon the internal parameter $\kappa$. Indeed the gaps between consecutive radii satisfy
\[R_{\mathrm{0,n}}(\kappa)-R_{\mathrm{0,n+1}}(\kappa)= \frac{\varepsilon}{\sqrt{\lambda_0}}+\caO\left(\varepsilon^{\frac32}\right),\]
while the range of the radii over the values of $\kappa$ is bounded by 
$|R_{\mathrm{0,n}}(\kappa_0)-R_{\mathrm{0,n}}(0)|\leq \caO(\varepsilon^2).$
While we have established the existence of radii $R_0$ which support pearled bilayers, there also may exist radii, and corresponding
total masses,for which no pearled circular bilayer solutions exist local to the associated circular bilayer; see Figure \ref{f:R0Kappa}.

\begin{minipage}{0.58\textwidth}
 \begin{figure}[H]
 \centering
     \includegraphics[width=\textwidth]{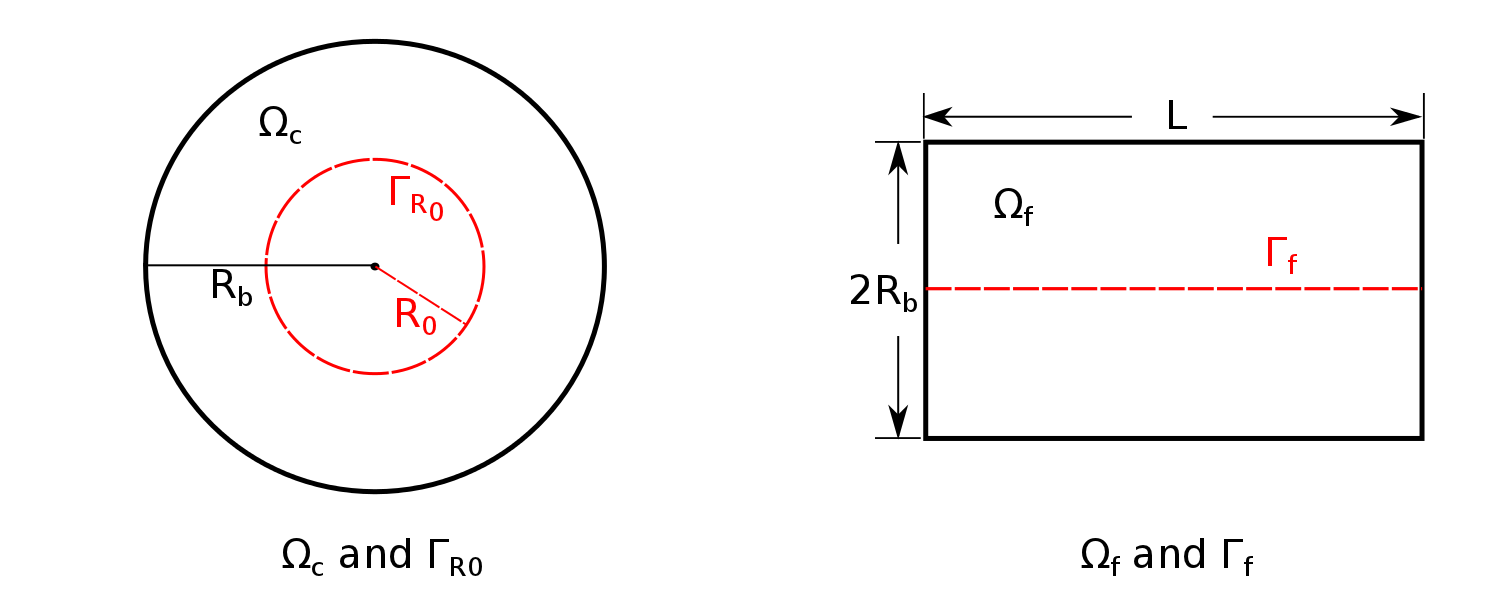}
                \caption{(left) A circular bilayer with interface $\Gamma_{R_0}$ in a concentric domain of radius $R_b$. (right) A flat
                bilayer with interface $\Gamma_f$ which intersects the rectangular domain at a right angle.}
    \label{f:domain}
        \end{figure}
\end{minipage}\qquad
\begin{minipage}{0.36\textwidth}
        \begin{figure}[H]
                \includegraphics[width=\textwidth]{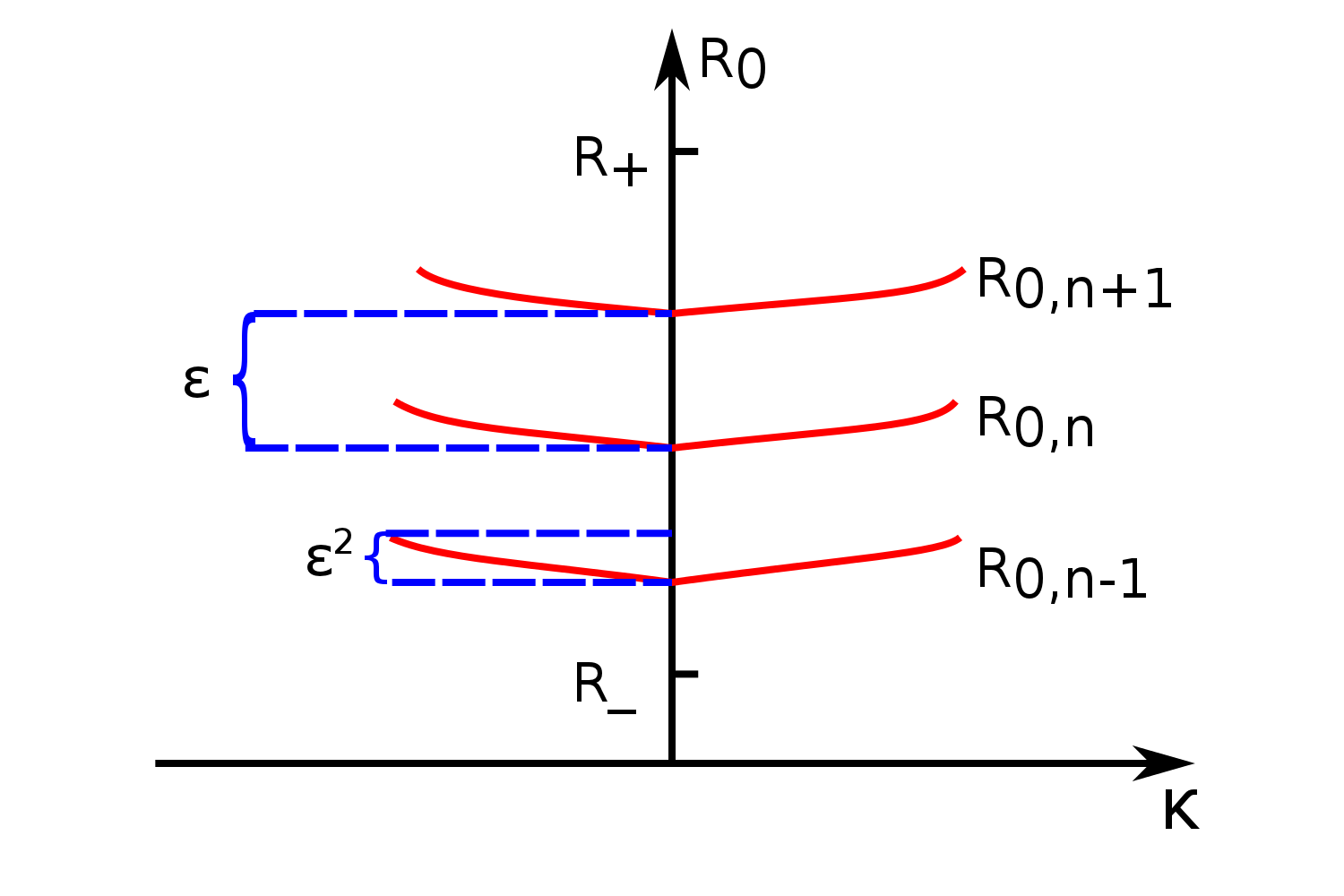}
                \caption{The admissible radii $\{R_{0,\textrm{n}}\}$ graphed verses $\kappa$ for fixed $\varepsilon$. The gaps between
                successive radii are $\caO(\eps)$ while the variation in $R_{0,{\textrm{n}}}$ with $\kappa$ is $\caO(\varepsilon^2).$}
    \label{f:R0Kappa}
        \end{figure}
\end{minipage}
 
As an existence problem, these scalings  imply that an $\caO(\varepsilon^3)$ change in the mass faction, which corresponds
to an $\caO(\varepsilon^2)$ change in the bilayer radius $R_0$, can induce an $\caO(1)$ impact on $\kappa$, and hence an
$\caO(\sqrt{\varepsilon})$ influence on the pearling amplitude of the associated equilibrium. This sensitivity of the pearling amplitude
to the mass fraction exemplifies the degeneracy  of the pearled morphologies. The size of the pearled ``beads" is fixed, but the amplitude
of the pearling pattern couples sensitively to the full system.  In particular for the strong FCH gradient flow, \eqref{e:FCH}, the possibility of 
non-existence of pearled morphologies at particular mass fractions and the delicate interaction between the radius of a circular bilayer and the 
amplitude of the high-frequency pearled morphology suggest a complex problem whose resolution may be quite sensitive to 
numerical truncation error.

\section{Pearling of the Flat Planar Bilayer}\label{s:2}

This section presents the construction of the pearled solutions $u_\mathrm{p}$ to the stationary strong-FCH \eqref{e:2sFCH} about
an infinite, flat, co-dimension one interface, $\Gamma_f$ embedded in $\R^2.$ 
The extended pearled solutions $u_\mathrm{p}$ are small-amplitude modulations of the extended flat bilayers $u_h$, 
periodic in the flat direction $\tau$. The construction is organized as follows: In Section \ref{ss:21}, the application of spatial dynamics techniques, together with a center manifold reduction, reduces the FCH equation to an 8th order ODE system; the derivation of the leading-order terms of the reduced ODE system are summarized in Section \ref{ss:22} with the details relegated to the Appendix. A normal form analysis presented in Section \ref{ss:23} reveals the pearling bifurcation structure; and in section \ref{ss:24}, it is shown that the pearling norm form admits a family of periodic orbits, which persist as solutions of  the full reduced ODE system, yielding the extended pearled solutions $u_\mathrm{p}$ of Theorem~\ref{t:main}.

\subsection{Spatial dynamics and center manifold reduction}\label{ss:21}
The spatial dynamics analysis begins by re-writing equation \eqref{e:EPP} as an infinite-dimension dynamical system in the rescaled 
$\tau$ variable followed by a normal form reduction on the associated center manifold.

To this end,  we rescale $\tau$ by $t=\frac{\sqrt{\lambda_0}}{\varepsilon}\tau$ and search for extended pearled solutions 
$u_\mathrm{rp}$ of
\begin{equation}
 \label{e:EPP}
 \left(\partial_r^2-W^{\prime\prime}(u)+\lambda_0\partial_t^2+\varepsilon\eta_1\right)
 \left(\partial_r^2u-W^\prime(u)+\lambda_0\partial_t^2u\right)+\varepsilon\eta_d W^\prime(u)-\varepsilon\gamma=0,
\end{equation}
which satisfy boundary conditions at infinity, 
\begin{equation}
 \lim_{r\rightarrow \pm\infty}|u_{\mathrm{rp}}(t,r)-u_-(\varepsilon)|=0, \text{ for all }t\in\R,
\end{equation}
and are even and  $T_\mathrm{rp}$-periodic in $t$,
\begin{equation}\label{e:pc}
 u_\mathrm{rp}(-t,r)=u_\mathrm{rp}(t,r), \quad u_\mathrm{rp}(t+T_\mathrm{rp},r)=u_\mathrm{rp}(t,r), \text{ for all }(t,r)\in\R^2,
\end{equation}
where  $T_\mathrm{rp}$ is to be determined.

We replace $u$ with $u_h+\delta u$ in \eqref{e:EPP} and consider the equation of the perturbation $\delta u$. For brevity, we denote the 
the perturbation by ``$u$'', instead of ``$\delta u$''. The perturbation solves the system
\begin{equation}
 \label{e:EPU}
 \mathcal{L}u+\mathcal{F}(u)=0,
\end{equation}
where the linear operator
\begin{equation}
 \mathcal{L}:=\left(\caL_h+\lambda_0\partial_t^2+\varepsilon\eta_1\right)\left(\caL_h+\lambda_0\partial_t^2\right)+\mathcal{M},
\end{equation}
is expressed in terms of the second order operator, $\caL_h:=\partial_r^2-W^{\pprime}(u_h)$ and the potential 
\[\mathcal{M}:=\varepsilon\eta_d W^{\pprime}(u_h)-\left(\partial_r^2u_h-W^\prime(u_h)\right)W^{\tprime}(u_h),\] 
while the nonlinearity given by
\begin{equation}
\label{e:F}
\begin{aligned}
 \mathcal{F}(u,\varepsilon):=&-\lambda_0W^{\tprime}(u_h+u)\left(\partial_tu\right)^2-2\lambda_0\left(W^{\prime\prime}(u_h+u)-W^{\prime\prime}(u_h)\right)\partial_t^2u-\\
 &\left[\caL_h+\varepsilon(\eta_1-\eta_d)-\left(W^{\pprime}(u_h+u)-W^{\pprime}(u_h)\right)\right]\left(W^\prime(u_h+u)-W^\prime(u_h)-W^{\pprime}(u_h)u\right)-\\
 &\left(W^{\pprime}(u_h+u)-W^{\pprime}(u_h)\right)\caL_h u-
 \left(\partial_r^2u_h-W^\prime(u_h)\right)\left(W^{\pprime}(u_h+u)-W^\pprime(u_h)-W^\tprime(u_h)u\right).
\end{aligned}
\end{equation}
We recast the system \eqref{e:EPU} in the vector form 
\begin{equation}
 \label{e:2ddFCH}
 \dot{U}=\bbL(\varepsilon)U+\bbF(U,\varepsilon),
\end{equation}
using the transformation $U_1=u$, $U_2=u_t$, $U_3=\lambda_0u_{tt}+\caL_h u$, $U_4=\partial_t\left(\lambda_0u_{tt}+\caL_h u\right)$ and introducing
\[
U = \begin{pmatrix}U_1\\U_2\\U_3\\U_4\end{pmatrix},\quad
\bbL(\varepsilon) = \begin{pmatrix}
0&1&0&0\\ -\frac{1}{\lambda_0}\caL_h&0&\frac{1}{\lambda_0}&0\\ 0&0&0&1\\
-\frac{1}{\lambda_0}\mathcal{M}& 0&
-\frac{1}{\lambda_0}(\caL_h+\varepsilon\eta_1)&0\end{pmatrix},\quad
\bbF(U,\varepsilon) = \begin{pmatrix}0\\0\\0\\-\frac{1}{\lambda_0}\mathcal{F}\end{pmatrix}.
\]
\begin{Remark}
To avoid technicalities we search for $u_p$ for a fixed value of $\gamma$. It is straightforward to recover the smooth dependence 
of $u_p$ with respect to $\gamma$. 
\end{Remark}
We observe that, for given small $\varepsilon$, 
$\bbL(\varepsilon): \mathcal{D}(\bbL)\rightarrow \mathcal{X}$ is a closed operator defined in the Hilbert space $\mathcal{X}$ with 
its domain $\mathcal{D}(\bbL)=\mathcal{Y}$, where
\[
 \mathcal{X}=H^3(\R)\times H^2(\R) \times H^1(\R) \times L^2(\R), \quad 
 \mathcal{Y}=H^4(\R) \times H^3(\R)\times H^2(\R) \times H^1(\R).
\] 

In the sequel we replace $\partial_t u$ and $\partial_t^2 u$ with $U_2$ and 
$\frac{1}{\lambda_0}(U_3-\caL_h U_1)$, respectively, in equation \eqref{e:F} for $\mathcal{F}$.
The map $\bbF: \mathcal{Y}\times [-\varepsilon_0,\varepsilon_0]\rightarrow\mathcal{Y}$ is smooth, for
$\varepsilon_0>0$ is sufficiently small.

\begin{figure}[H]
\centering
    \includegraphics[width=0.40\textwidth]{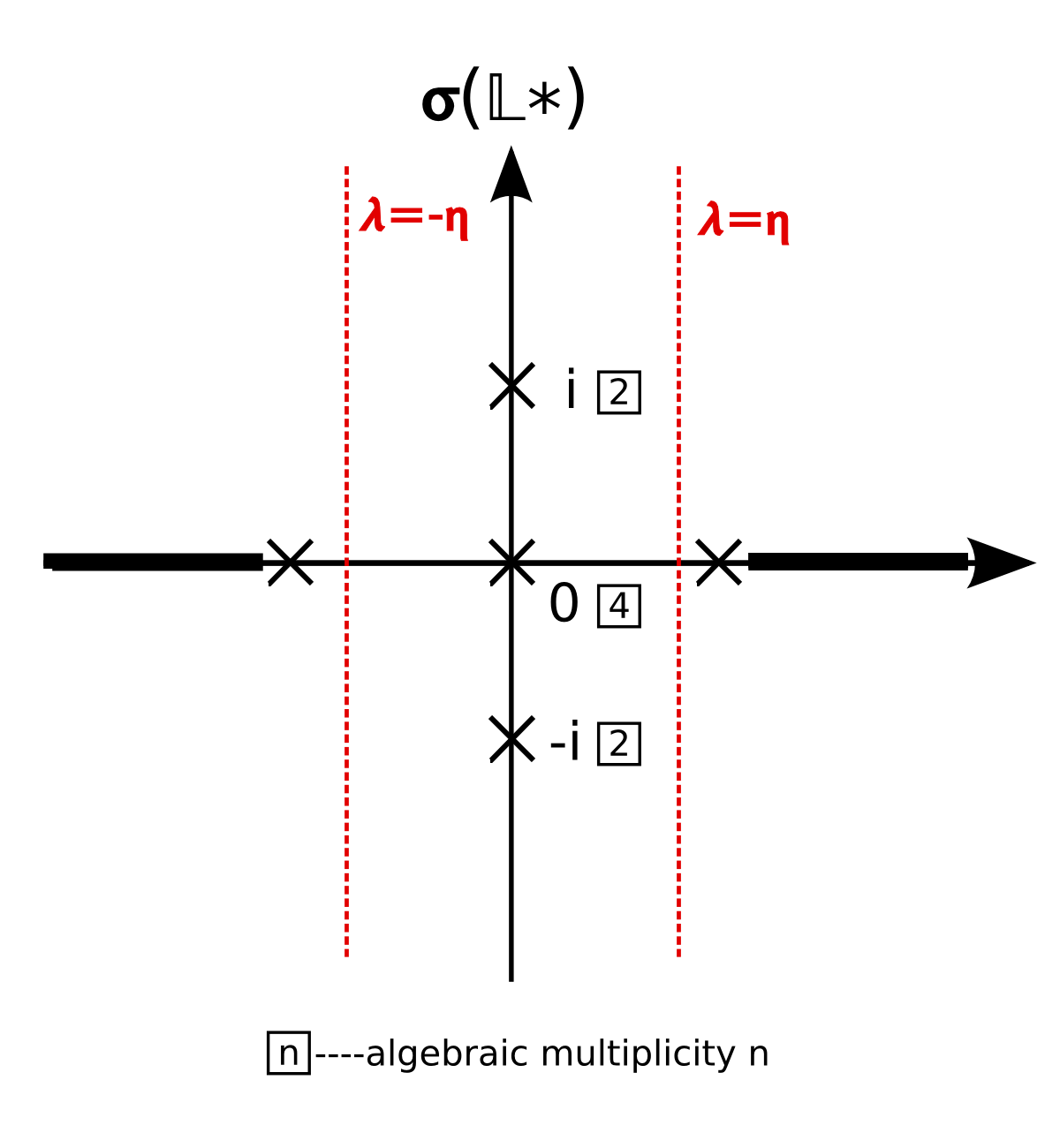}
  \caption{The spectrum of $\bbL_*$ indicating the center eigenvalues and their multiplicity.}
\label{f:specL}
\end{figure}

\begin{Lemma}
\label{l:spec}
The spectrum of $\bbL_*:=\bbL(0,0)$, $\sigma(\bbL_*)$, as shown in Figure \ref{f:specL}, satisfies
\begin{itemize}
 \item[(\rmnum{1})]$\sigma_c(\bbL_*):=\sigma(\bbL_*)\cap i\R=\{0, \pm i\}$, where eigenvalue $0$ has geometric multiplicity 1 and algebraic multiplicity 4, and eignvalues
 $\pm \rmi$ have geometric multiplicity 1 and algebraic multiplicity 2.
 \item[(\rmnum{2})]There exists $\eta>0$ such that $\sigma(\bbL_*)\cap \{|\Re\lambda|\leq \eta\}=\sigma_c(\bbL_*)$.
\end{itemize}
\end{Lemma}

\begin{Proof}
We first introduce the operator 
\begin{equation*}
 \begin{matrix}
  \mathcal{L}^\lambda: & H^4(\R) & \longrightarrow & L^2(\R) \\
  & u & \longmapsto & \left(\caL_0+\lambda_0\lambda^2\right)^2u
 \end{matrix}
\end{equation*}
 which, for any $\lambda\in\C$, has the same Fredholm properties as the operator $\bbL_*-\lambda\id$; see a similar case in \cite{ssmorse_2008} for
 a detailed proof. More specifically, $\bbL_*-\lambda\id$ is Fredholm if and only if $\mathcal{L}^\lambda$ is Fredholm. In addition, if Fredholm, then $\bbL_*-\lambda\id$ and $\mathcal{L}^\lambda$ have the same Fredholm index. We omit the technical details required to establish that 
 $\dim\cok (\bbL_*-\lambda\id)=\dim\cok \mathcal{L}^\lambda$; however it is straightforward to see that
 $$\dim\ker (\bbL_*-\lambda\id)=\dim\ker \mathcal{L}^\lambda,$$
 since
 \[
 (\bbL_*-\lambda\id)\begin{pmatrix}U_1\\ U_2 \\ U_3 \\ U_4 \end{pmatrix}=0 \Longleftrightarrow \mathcal{L}^\lambda U_1=0.
 \]
 To obtain the spectral properties of $\bbL_*$,  the dispersion relation of $\mathcal{L}^\lambda$  implies that
 \[
  \sigma(\bbL_*)=\{\lambda\in\C\mid (\mu+\lambda_0\lambda^2)^2=0, \text{ for some }\mu\in\sigma(\caL_0)\},
 \]
 where $\caL_0$, defined in \eqref{cL_0-def}, is of Sturm-Liouville type with simple, real spectrum thats satisfies 
 \[
 \sigma(\caL_0)\cap \{\Re\lambda\geq 0\}=\{0,\lambda_0\}, \quad
 \sigma(\caL_0)\cap \{\Re\lambda< 0\}\subset (\infty,-c), \text{ for some }c>0.
 \]
These observations conclude the proof.
\end{Proof}

The center space $\mathcal{X}_c$ of $\bbL_*$, that is, the spectral subspace associated to $\sigma_c(\bbL_*)$, is 8-dimensional
and spanned by the eigenfunctions $\{E_1, E_2,\bar{E}_1,\bar{E}_2, F_1, F_2, F_3, F_4\}$, where
\begin{equation}
\begin{aligned}
E_1=\begin{pmatrix}1\\ \rmi \\ 0 \\ 0 \end{pmatrix}\psi_0,\quad
E_2=\begin{pmatrix}\rmi\\ 0    \\ 2\lambda_0\rmi \\ -2\lambda_0 \end{pmatrix}\psi_0,\quad
F_1=\begin{pmatrix}1 \\ 0    \\ 0 \\ 0 \end{pmatrix}\psi_1, \\
F_2=\begin{pmatrix}0 \\ 1    \\ 0 \\ 0 \end{pmatrix}\psi_1,\quad
F_3=\begin{pmatrix}0 \\ 0    \\ \lambda_0 \\ 0 \end{pmatrix}\psi_1,\quad
F_4=\begin{pmatrix}0 \\ 0    \\ 0 \\ \lambda_0 \end{pmatrix}\psi_1.
\end{aligned}
\end{equation}
Moreover, these generalized eigenfunctions of $\bbL_*$ satisfies 
\begin{equation}\label{e:eigen}
\begin{aligned}
 (\bbL_*-\rmi)E_1=0, \quad  (\bbL_*-\rmi)E_2=E_1, \quad \bbL_*F_1=0, \quad \bbL_*F_2=F_1,\text{ }\\
 (\bbL_*+\rmi)\bar{E}_1=0, \quad  (\bbL_*+\rmi)\bar{E}_2=\bar{E}_1,  \quad \bbL_*F_3=F_2, \quad  \bbL_*F_4=F_3,\\
S_1^2=\id,\quad S_1E_1=\bar{E}_1,\quad S_1E_2=-\bar{E}_2,\quad S_1F_j=F_j, \quad S_1F_k=-F_k, j=1,3; k=2,4,\\
S_2^2=\id,\quad S_2E_j=E_j,\quad S_2\bar{E}_j=\bar{E}_j,\quad S_2F_k=-F_k, j=1,2; k=1,2,3,4. 
\end{aligned}
\end{equation}
where $S_1$ and $S_2$ are the symmetries inherited from the $t\rightarrow-t$ and $r\rightarrow -r$ symmetries of the original PDE \eqref{e:EPP}. Here $S_1$ is a reversible symmetry and plays a crucial role in the subsequent bifurcation analysis. From \eqref{e:eigen} we develop an explicit expression of the spectral projection $\mathbb{P}_c:\mathcal{X}\rightarrow \mathcal{X}_c$, 
\begin{equation}\label{e:uc0}
\begin{aligned}
 U_c:=\mathbb{P}_c U=& \langle U, E_1^{\ad}\rangle E_1+ \langle U,E_2^{\ad}\rangle E_2+ \langle U,\bar{E}_1^{\ad}\rangle \bar{E}_1+ \langle U,\bar{E}_2^{\ad}\rangle \bar{E}_2+\\
 & \langle U, F_1^{\ad}\rangle F_1+ \langle U,F_2^{\ad}\rangle F_2+ \langle U,F_3^{\ad}\rangle F_3+ \langle U,F_4^{\ad}\rangle F_4,
\end{aligned}
\end{equation}
where 
\begin{equation}
\begin{aligned}
E_1^{\ad}=\begin{pmatrix} \frac{1}{2} \\ \frac{\rmi}{2}  \\ -\frac{1}{4\lambda_0} \\ 0 \end{pmatrix}\psi_0,\quad
E_2^{\ad}=\begin{pmatrix} 0 \\ 0    \\ \frac{\rmi}{4\lambda_0} \\ -\frac{1}{4\lambda_0} \end{pmatrix}\psi_0,\quad
F_1^{\ad}=\begin{pmatrix}1 \\ 0    \\ 0 \\ 0 \end{pmatrix}\psi_1, \\
F_2^{\ad}=\begin{pmatrix}0 \\ 1    \\ 0 \\ 0 \end{pmatrix}\psi_1,\quad
F_3^{\ad}=\begin{pmatrix}0 \\ 0    \\ \frac{1}{\lambda_0} \\ 0 \end{pmatrix}\psi_1,\quad
F_4^{\ad}=\begin{pmatrix}0 \\ 0    \\ 0 \\ \frac{1}{\lambda_0} \end{pmatrix}\psi_1.
\end{aligned}
\end{equation}
These vector functions with superscript ``$\ad$'' are generalized eigenfunctions of the adjoint operator $\bbL_*^\ad$ associated to $0$ and $\pm\rmi$ 
in $(L^2(\R))^4$ with canonical inner product $\langle\cdot,\cdot\rangle$.
Moreover, a standard calculation \cite{kappro} shows that, for any given $w_0>1$, there exists $C\geq1$ such that 
\begin{equation}
\label{e:nest}
 \|(\rmi w - \bbL_*)^{-1}U\|_{\mathcal{X}}\leq \frac{C}{|w|} \|U\|_{\mathcal{X}}, \text{ for all } |w|\geq w_0, w\in\R, U\in (\id-\mathbb{P}_c)\mathcal{X}.
\end{equation}
Therefore, based on Lemma \ref{l:spec} and the norm estimate \eqref{e:nest} on $\bbL_*|_{(\id-\mathbb{P}_c)\mathcal{X}}$,
we can apply the center manifold reduction theorem to the system \eqref{e:2ddFCH} and obtain the following proposition (see \cite[Theorem 2.9]{harioo}).
\begin{Proposition}
 Given any fixed $\gamma$ and $k\in\Z^+$, there exist open sets containing the origin $\mathcal{U}\subset \mathcal{X}_c$, $\mathcal{V}\subset (\id-\mathbb{P}_c)\mathcal{Y}$, $\mathcal{W}\in\R$,
 and a $C^k$-smooth map $\Psi:\mathcal{U}\times\mathcal{W}\rightarrow\mathcal{V}$, for any fixed nonnegative integer $k$, such that the center manifold $\mathcal{M}_c$,
 that is, the graph of the map $\Psi$, has the following properties.
 \begin{itemize}
  \item[(\rmnum{1})]The center manifold $\mathcal{M}_c$ is tangent to the center eigenspace $\mathcal{X}_c$,
  \begin{equation}
   \|\Psi(U_c,\varepsilon)\|_\mathcal{Y}=\caO(|\varepsilon|\|U_c\|+\|U_c\|^2).
  \end{equation}
  \item[(\rmnum{2})]The center manifold $\mathcal{M}_c$ is locally invariant, that is, 
  if $U$ is a solution to \eqref{e:2ddFCH} with $U(0)\in\mathcal{M}_c$ and
  $U(t)\in \mathcal{U}\times\mathcal{V}$ for $t\in[0, T]$, then $U(t)\in\mathcal{M}_c$ for all $t\in[0,T]$.
  \item[(\rmnum{2})]The center manifold $\mathcal{M}_c$ contains all bounded solutions to \eqref{e:2ddFCH} with $\R$ as the existence interval, that is, 
  if $U$ is a solution to \eqref{e:2ddFCH} satisfying $\{U(t)\mid t\in\R\}\subset \mathcal{U}\times\mathcal{V}$, then $\{U(t)\mid t\in\R\}\subset\mathcal{M}_c$.
 \end{itemize}
\end{Proposition}

\subsection{Reduced center manifold ODE}\label{ss:22}
In this section we calculate the reduced ODE system obtained by restricting \eqref{e:2ddFCH} to the center manifold. From the analysis presented
in Section~\ref{ss:21} and summarized in Figure~\ref{f:specL} it follows that the reduced ODE system is of $8$-th order which can be viewed as a 
coupling of two four-dimensional systems which exhibit the  so-called ``reversible-Hopf bifurcation'' and the ``reversible $0^{4+}$ bifurcation''. 
Moreover, the coupling occurs at the nonlinear level and is weak. 
On the linear level, the  $S_1$-reversibility of the reduction to the $\pm\rmi$-eigenspace gives rise to the ``reversible-Hopf bifurcation'', which is 
well-studied; see \cite{glebsky_1995,iooper_1993};  while the  $S_1$-reversibility of the $0$-eigenspace gives rise to the ``reversible $0^{4+}$ bifurcation'', whose study is quite open; see \cite{harioo}. Fortunately, extended pearled solutions result from the ``reversible-Hopf bifurcation''. 
Moreover, it is known that the analysis of this bifurcation relies on the coefficients of the cubic terms in the norm form \cite{iooper_1993}.
Therefore, all the necessary terms of the reduced ODE system, up to cubic order, are explicitly determined in this section. 

To restrict the system \eqref{e:2ddFCH} to the center manifold we consider $U$ in the form 
\begin{equation}
\label{e:sch}
 U=U_c+\Psi(U_c,\varepsilon).
\end{equation}
Substituting this form \eqref{e:sch} into \eqref{e:2ddFCH} and applying the projection $\mathbb{P}_c$, we obtain the reduced equation,
\begin{equation}
\label{e:2dred}
 \dot{U}_c=\bbL_*U_c+\mathbb{P}_c\Big(\mathbb{M}(\varepsilon)\big(U_c+\Psi(U_c,\varepsilon)\big)+\bbF\big(U_c+\Psi(U_c,\varepsilon),\varepsilon\big)\Big),
\end{equation}
where $\mathbb{M}(\varepsilon):=\bbL(\varepsilon)-\bbL_*$.
Moreover,  from \eqref{e:uc0}, we note that $U_c$ admits the general expression
\begin{equation}
\label{e:ucg}
 U_c(t)=\sum_{j=1}^2\left(A_j(t)E_j+\bar{A}_j(t)\bar{E}_j\right)+\sum_{k=1}^4B_k(t)F_k,
\end{equation}
Using this expression of $U_c$, we rewrite the reduced system \eqref{e:2dred} explicitly in terms of 
\begin{equation}\label{e:bfA}
 \mathbf{A}:=(A_1,A_2,\bar{A}_1,\bar{A}_2,B_1,B_2,B_3,B_4).
\end{equation}
We summarize the essential result into Lemma \ref{l:sim}, relegating the detailed results and concomitant calculations to Appendix \ref{ss:31}.  
The principle technicality in the calculation lies in finding the explicit expression of
$\Psi_{(2,0,0)}(U_c,U_c)$ in terms of $\mathbf{A}$; see Lemma \ref{l:31} for details.
\begin{Lemma}\label{l:sim}
The reduced system \eqref{e:2dred}, in terms of ${\bf A}$, called the reduced ODE system, admits the expression 
\begin{equation}
 \label{e:2dcredA}
 \dot{\mathbf{A}}=\mathbf{L}(\varepsilon)\mathbf{A}+\mathbf{R}_2(\mathbf{A})
 +\mathbf{R}_3(\mathbf{A})+ \caO\left(|\varepsilon|^2\|{\bf A}\|+|\varepsilon|\|{\bf A}\|^2+\|{\bf A}\|^4\right),
\end{equation}
where the linear term ${\bf L}$, the quadratic term ${\bf R}_2$, the cubic term ${\bf R}_3$ are of the following expressions.
\begin{equation}\label{e:bfL}
\mathbf{L}(\varepsilon)=\begin{pmatrix}
\rmi(1+\mu_1\varepsilon) &  1-\mu_1\varepsilon  & \rmi \mu_1\varepsilon  & \mu_1\varepsilon
& 0 & 0 & 0 &0 \\
\mu_2\varepsilon & \rmi \left(1+\mu_3\varepsilon \right) &  \mu_2\varepsilon & -\rmi \mu_3 \varepsilon
& 0 & 0 & 0 &0 \\
-\rmi \mu_1\varepsilon & \mu_1 \varepsilon& -\rmi\left(1+\mu_1\varepsilon\right) &  1+\mu_1\varepsilon
& 0 & 0 & 0 &0 \\
\mu_2\varepsilon & \rmi \mu_3\varepsilon &\mu_2 \varepsilon& -\rmi \left(1+\mu_3 \varepsilon\right) &   0 & 0 & 0 &0 \\
0 & 0 & 0 &0 & 0 & 1 & 0 & 0\\
0 & 0 & 0 &0 & \mu_4 \varepsilon& 0 & 1 & 0\\
0 & 0 & 0 &0 & 0 &0 & 0 & 1\\
0 & 0 & 0 &0 & \mu_5\varepsilon & 0 & \mu_6\varepsilon & 0\\
\end{pmatrix},
\end{equation}
\begin{equation*}
\mathbf{R}_2(\mathbf{A})=\left(0,R_{2,2},0,\bar{R}_{2,2},0,0,0,R_{2,8}\right)^T, \quad
\mathbf{R}_3({\bf A})=\left(0,R_{3,2},0,\bar{R}_{3,2},0,0,0,R_{3,8}\right)^T, 
\end{equation*}
where the expressions of every $\mu_j\in \R$ and $R_{2\backslash 3, 2\backslash 8}$ in terms of ${\bf A}$ can be found in Lemma \ref{l:reducedA}.
\end{Lemma}

\subsection{Norm forms}\label{ss:23}
We obtain a normal form of the leading-order-term reduced system via a composition of a linear versal transformation and a near-identity nonlinear transformation. The versal transformation allows a Jordan-form type decomposition which is smooth in the parameters, see \cite{viarnold} for full details. 

\begin{Lemma}\label{l:versal}
 For sufficiently small $\varepsilon$, there exists a smooth linear map $\mathbf{T}(\varepsilon)$ with $\mathbf{T}(0)=\id$ such that 
 under the transformation
\[
 \mathbf{A}=\mathbf{T}(\varepsilon)\mathbf{C}, \quad \mathbf{C}=(C_1,C_2,\bar{C}_1,\bar{C}_2, D_1,D_2,D_3,D_4)^T,
\]
the linear part of \eqref{e:2dcredA} in $\mathbf{A}$, that is,
 \begin{equation}
 \label{e:l2dcred}
 \dot{\mathbf{A}}=\mathbf{L}(\varepsilon)\mathbf{A},
\end{equation}
takes the versal normal  form 
\begin{equation}
 \dot{\mathbf{C}}=\mathscr{L}(\varepsilon)\mathbf{C}+\caO\left(|\varepsilon|^2\|{\bf C}\|\right),
\end{equation}
where
\begin{equation}
 \label{e:ln2dcred}
 \mathscr{L}(\varepsilon)=\begin{pmatrix}
  \rmi\left(1+\omega_1\varepsilon\right) & 1 & 0 & 0 & 0 & 0 & 0 & 0 \\
  \omega_2\varepsilon& \rmi\left(1+\omega_1\varepsilon\right) & 0 & 0 & 0 & 0 & 0 & 0 \\
  0 & 0 & -\rmi\left(1+\omega_1\varepsilon\right) & 1 & 0 & 0 & 0 & 0 \\
  0 & 0 & \omega_2 \varepsilon & -\rmi\left(1+\omega_1\varepsilon\right) & 0 & 0 & 0 & 0 \\
  0 & 0 & 0 & 0 & 0 & 1 & 0 & 0 \\
  0 & 0 & 0 & 0 & 0 & 0 & 1 & 0 \\
  0 & 0 & 0 & 0 & 0 & 0 & 0 & 1 \\
  0 & 0 & 0 & 0 & \omega_3\varepsilon & 0 & \omega_4\varepsilon & 0 \\
 \end{pmatrix}.
\end{equation}
Here we have introduced
\begin{equation}\label{e:ls}
\omega_1=\frac{1}{2}(\mu_1+\mu_3),\quad \omega_2=\mu_2,\quad \omega_3=\mu_5,\quad \omega_4=\mu_4+\mu_6,
\end{equation}
where the expression of each $\mu_j\in\R$ can be found in Lemma \ref{l:reducedA}. 
\end{Lemma}
\begin{Proof}
We point out that $\mathbf{L}(\varepsilon)$ inherits the symmetries $\tau\rightarrow -\tau$ and $r\rightarrow -r$ of the original PDE \eqref{e:2sFCH}, that is,
\[
 S_1\mathbf{L}(\varepsilon)=-\mathbf{L}(\varepsilon)S_1,\quad S_2\mathbf{L}(\varepsilon)=\mathbf{L}(\varepsilon)S_2,
\]
where 
\begin{equation*}
  \begin{aligned}
  S_1(A_1,A_2,\bar{A}_1,\bar{A}_2, B_1,B_2,B_3,B_4)^T=(\bar{A}_1,-\bar{A}_2,A_1,-A_2, B_1,-B_2,B_3,-B_4)^T,\\
 S_2(A_1,A_2,\bar{A}_1,\bar{A}_2, B_1,B_2,B_3,B_4)^T=(A_1,A_2,\bar{A}_1,\bar{A}_2, -B_1,-B_2,-B_3,-B_4)^T.
  \end{aligned}
 \end{equation*}
Then, according to \cite[Theorem 4.4]{viarnold}, a versal deformation of the Jordan normal form $\mathbf{L}$ keeping the symmetries can be chosen in the form
\begin{equation*}
\begin{pmatrix}
  \rmi\left(1+\widetilde{\omega}_1\right) & 1 & 0 & 0 & 0 & 0 & 0 & 0 \\
  \widetilde{\omega}_2& \rmi\left(1+\widetilde{\omega}_1\right) & 0 & 0 & 0 & 0 & 0 & 0 \\
  0 & 0 & -\rmi\left(1+\widetilde{\omega}_1\right) & 1 & 0 & 0 & 0 & 0 \\
  0 & 0 & \widetilde{\omega}_2 & -\rmi\left(1+\widetilde{\omega}_1\right) & 0 & 0 & 0 & 0 \\
  0 & 0 & 0 & 0 & 0 & 1 & 0 & 0 \\
  0 & 0 & 0 & 0 & 0 & 0 & 1 & 0 \\
  0 & 0 & 0 & 0 & 0 & 0 & 0 & 1 \\
  0 & 0 & 0 & 0 & \widetilde{\omega}_3 & 0 & \widetilde{\omega}_4 & 0 \\
 \end{pmatrix}.
\end{equation*}
where $\widetilde{\omega}_j(\varepsilon)\in\R$ with $\omega_j(0)=0$ for $j=1,2,3,4$.
Comparing the coefficients of the characteristic polynomials of the two $4\times 4$ diagonal blocks associated to $(A_1,A_2, \bar{A}_1, \bar{A}_2)$ in \eqref{e:bfL} and \eqref{e:ln2dcred}, we have
\begin{equation*}
 \begin{cases}
  (1+\widetilde{\omega}_1)^2-\widetilde{\omega}_2=1+(\mu_1-\mu_2+\mu_3)\varepsilon,\\
  \left((1+\widetilde{\omega}_1)^2+\widetilde{\omega}_2\right)^2=1+2(\mu_1+\mu_2+\mu_3)\varepsilon-4(\mu_2-\mu_3)\mu_1\varepsilon^2,
 \end{cases}
\end{equation*}
from which we have
\begin{equation*}
\begin{aligned}
\widetilde{\omega}_1(\varepsilon)&=\frac{1}{2}(\mu_1+\mu_3)\varepsilon+\caO\left(|\varepsilon|^2\right),
&  \widetilde{\omega}_2(\varepsilon)&=\mu_2\varepsilon+\caO\left(|\varepsilon|^2\right).
 \end{aligned}
\end{equation*}
Similarly, we have
\begin{equation*}
\begin{aligned}
\widetilde{\omega}_3(\varepsilon)&=\mu_5\varepsilon+\caO\left(|\varepsilon|^2\right),
& \widetilde{\omega}_4(\varepsilon)&=(\mu_4+\mu_6)\varepsilon+\caO\left(|\varepsilon|^2\right).
\end{aligned}
\end{equation*}
We truncate this versal deformation up to linear terms in $\varepsilon$, denote it as $\mathscr{L}(\varepsilon)$ and conclude our proof.
\end{Proof}

On the other hand, we have the following nonlinear normal form.
\begin{Lemma}
  There exist smooth families of degree-$2$ polynomials  
  \[
   \Phi_2=(\Phi_{2,1},\Phi_{2,2},\Phi_{2,3},\Phi_{2,4},\Phi_{2,5},\Phi_{2,6},\Phi_{2,7},\Phi_{2,8})^T,
  \]
  and degree-$3$ polynomials
  \[
   \Phi_3=(\Phi_{3,1},\Phi_{3,2},\Phi_{3,3},\Phi_{3,4},\Phi_{3,5},\Phi_{3,6},\Phi_{3,7},\Phi_{3,8})^T,
  \]
  in terms of $\mathbf{C}$ such that 
  such that under the near-identity transformation
  \begin{equation}\label{e:nonl}
 \mathbf{A}=\mathbf{C}+\Phi_2(\mathbf{C})+\Phi_3(\mathbf{C}),
\end{equation}
 the nonlinear part of \eqref{e:2dcredA}, that is,
 \begin{equation}
 \label{e:q2dcred}
 \dot{\mathbf{A}}=\mathbf{L}(0)\mathbf{A}+\mathbf{R}_2(\mathbf{A},\mathbf{A})
 +\mathbf{R}_3(\mathbf{A},\mathbf{A},\mathbf{A}),
\end{equation}
takes the normal form
\begin{equation}\label{e:nfnonl}
 \dot{\mathbf{C}}=\mathbf{L}(0)\mathbf{C}+\mathscr{R}_2(\mathbf{C},\mathbf{C})
 +\mathscr{R}_3(\mathbf{C},\mathbf{C},\mathbf{C})+\caO(|\mathbf{C}|^4).
\end{equation}
Here $\mathscr{R}_2=0$ and
$\mathscr{R}_3=(\mathscr{R}_{3,1},\mathscr{R}_{3,2},\mathscr{R}_{3,3},\mathscr{R}_{3,4},\mathscr{R}_{3,5},\mathscr{R}_{3,6},\mathscr{R}_{3,7},\mathscr{R}_{3,8})^T$ is of the form
\begin{equation}
\label{e:nq2dcred}
\begin{aligned}
 \mathscr{R}_{3,1}=&\rmi \Big\{C_1\big[\alpha_7C_1\bar{C}_1+\alpha_8\rmi(C_1\bar{C}_2-\bar{C}_1C_2)\big]+\alpha_9C_1D_1^2+
 \alpha_{10}\rmi D_1(C_2D_1-C_1D_2)+\\
& \alpha_{11}C_1(2D_1D_3-D_2^2)+\alpha_{12}\rmi\big[C_1(D_2D_3-3D_1D_4)+C_2(2D_1D_3-D_2^2)\big]\Big\};\\
 \mathscr{R}_{3,2}=&\Big\{C_1\big[\alpha_1C_1\bar{C}_1+\alpha_2\rmi(C_1\bar{C}_2-\bar{C}_1C_2)\big]+\alpha_3C_1D_1^2+
 \alpha_{4}\rmi D_1(C_2D_1-C_1D_2)+\\
 &\alpha_5C_1(2D_1D_3-D_2^2)+\alpha_6\rmi\big[C_1(D_2D_3-3D_1D_4)+C_2(2D_1D_3-D_2^2)\big]\Big\}+\\
 &\rmi\Big\{C_2\big[\alpha_7C_1\bar{C}_1+\alpha_8\rmi(C_1\bar{C}_2-\bar{C}_1C_2)\big]+\alpha_9C_2D_1^2+
 \alpha_{10}\rmi D_2(C_2D_1-C_1D_2)+\\
& \alpha_{11}C_1(3D_1D_4-D_2 D_3)+\alpha_{12}\rmi\big[2C_1(2D_3^2-3D_2D_4)+C_2(3D_1D_4-D_2D_3)\big]\Big\};\\
 \mathscr{R}_{3,8}=&D_1(\beta_{1}C_1\bar{C}_1+\beta_{2}C_2\bar{C}_2)+\rmi(C_1\bar{C}_2-\bar{C}_1C_2)(\beta_3D_1+\beta_4D_3)+\beta_5C_1\bar{C}_1D_3+\\
 &\beta_6D_2(C_1\bar{C}_2+\bar{C}_1C_2)+\beta_7[3(C_1\bar{C}_2+\bar{C}_1C_2)D_4-2C_2\bar{C}_2D_3]+\beta_8D_1D_2^2+\\
& D_1^2(\beta_9D_1+\beta_{10}D_3)+
\beta_{11}(D_2^2D_3-2D_1D_3^2)+\beta_{12}(D_2^2D_3-3D_1D_2D_4)+\\
&\beta_{13}(9D_2D_3D_4-9D_1D_4^2-4D_3^3);\\
\mathscr{R}_{3,{j+2}}=&\bar{\mathscr{R}}_{3,j}, \quad j=1,2; \quad \quad \mathscr{R}_{3,k}=0,\quad k= 5,6,7; \\
\end{aligned}
\end{equation}
where the explicit expressions of the coefficients $\alpha_j$,$\beta_k\in\R$, are given in Lemma \ref{l:32}. Moreover, the transformation preserves the reversibility $S_1$ and the symmetry $S_2$.
\end{Lemma}
\begin{Proof}
Following \cite[Chapter 3]{harioo}, we cast the normal form problem as a solvability issue on a space of polynomials in ${\bf C}$ which is expressed in terms of the Fredholm alternative of the operator 
\begin{equation}
 \mathcal{D}=\left(\rmi C_1+C_2\right)\frac{\partial}{\partial C_1}+\rmi C_2\frac{\partial}{\partial C_2}+
  \left(-\rmi \bar{C}_1+\bar{C}_2\right)\frac{\partial}{\partial \bar{C}_1}+(-\rmi \bar{C}_2)\frac{\partial}{\partial \bar{C}_2}+
  \sum_{j=1}^3D_{j+1}\frac{\partial}{\partial D_i}.
\end{equation}
For convenience, we introduce the polynomial space $\mathbf{P}_j$, $j=2,3$, which is the set of all degree-$j$ homogeneous polynomials in $\mathbf{C}$, 
with the inner product
\[
 \langle P\mid Q\rangle=P(\partial_\mathbf{C})\bar{Q}(\mathbf{C})|_{\mathbf{C}=0}.
\]
We point out here that the conjugacy $\bar{Q}$ only acts on the coefficients, in the sense that, for example, for $Q({\bf C})=\rmi C_1^2$, $\bar{Q}({\bf C})=-\rmi C_1^2$.

More specifically, plugging \eqref{e:nonl} and \eqref{e:nfnonl} into \eqref{e:q2dcred}, we obtain the following two equalities.
 \begin{eqnarray}
   \big(\mathcal{D}-\mathbf{L}(0)\big)\Phi_2&=&\mathbf{R}_2-\mathscr{R}_2,\qquad\qquad\qquad\qquad\qquad\quad\quad \label{e:nf2}\\ 
   \big(\mathcal{D}-\mathbf{L}(0)\big)\Phi_3&=&\mathbf{R}_3+
   2\mathbf{R}_2(\mathbf{C},\Phi_2)-\mathscr{R}_3-\big(D_{\mathbf{C}}\Phi_2\big)\mathscr{R}_2, \label{e:nf3}
 \end{eqnarray}

From the Fredholm alternative, we may solve for $\Phi_{2}$ and  $\mathscr{R}_{2}$ uniquely in \eqref{e:nf2} subject to
\[
 \mathscr{R}_2\in\ker\left(\big(\mathcal{D}^{\ad}-\mathbf{L}^{\ad}(0)\big)|_{\mathbf{P}_2^8}\right),\quad \Phi_2\in \left(\ker\left((\mathcal{D}-\mathbf{L}(0))|_{\mathbf{P}_2^8}\right)\right)^\perp,
\]
where
\[
 \mathcal{D}^{\ad}=(-\rmi C_1)\frac{\partial}{\partial C_1}+\left( C_1 - \rmi C_2\right)\frac{\partial}{\partial C_2}+
  \rmi \bar{C}_1\frac{\partial}{\partial \bar{C}_1}+\left( \bar{C}_1+ \rmi \bar{C}_2\right)\frac{\partial}{\partial \bar{C}_2}+
  \sum_{j=2}^4D_{j-1}\frac{\partial}{\partial D_i}.
\]
In fact,  we claim that $\mathscr{R}_2=0$. To show this, we only need to verify that 
\[
 \mathbf{R}_2\in \Rg\left((\mathcal{D}-\mathbf{L}(0))|_{\mathbf{P}_2^8}\right)=\left(\ker\left((\mathcal{D}^{\ad}-\mathbf{L}^{ad}(0))|_{\mathbf{P}_2^8}\right)\right)^\perp,
\]
which follows from the expression of ${\bf R}_2$ in \eqref{e:Rjs} and the fact that 
\begin{equation}\
\label{e:ker2}
 \begin{aligned}
 \ker\big(\mathcal{D}^{\ad}|_{\mathbf{P}_2}\big)=\span\{&C_1\bar{C}_1,C_1\bar{C}_2-\bar{C}_1C_2,D_1^2,2D_1D_3-D_2^2\},\\
 \ker\big((\mathcal{D}^{\ad}+\rmi)|_{\mathbf{P}_2}\big)=\span\{&C_1D_1,C_1D_2-C_2D_1\}.\\
   \end{aligned}
\end{equation}
As a result, we  obtain $\Phi_2\in \left(\ker\left((\mathcal{D}-\mathbf{L}(0))|_{\mathbf{P}_2^8}\right)\right)^\perp$ with coefficients,
 \begin{equation}
 \label{e:phi2}
  \begin{aligned}
  \Phi_{2,1}(\mathbf{C})=& \frac{2\nu_1}{9}\Big[\big(3 c_{1,+}^2-12 c_{1,-}^2 +14 c_{1,+}c_{2,-}-40 c_{1,-}c_{2,+}-9 c_{2,-}^2\big)\\
             & +\rmi\big(12c_{1,+}c_{1,-} + 32 c_{1,+}c_{2,+} +4 c_{1,-}c_{2,-} \big)\Big]\\
             & +\nu_2\Big[\big(-\frac{1}{2}D_1^2+D_1D_3+2D_2^2-2D_2D_4-3D_3^2\big) +2\rmi\big(D_1D_2-2D_2D_3+2D_3D_4\big)\Big],\\
  \Phi_{2,2}(\mathbf{C})=& \frac{2\nu_1}{9}\Big[\big(-18c_{1,+}c_{1,-} + 12 c_{1,+}c_{2,+} + 6c_{1,-}c_{2,-}- 44 c_{2,+}c_{2,-}\big)\\
             & +\rmi\big(9 c_{1,+}^2 -30 c_{1,+}c_{2,-}+24 c_{1,-}c_{2,+}+32 c_{2,+}^2+13 c_{2,-}^2\big)\Big]\\
             & +\nu_2\Big[\big(D_1D_2+D_1D_4+D_2D_3-4D_3D_4\big)+\rmi\big(\frac{1}{2}D_1^2+D_1D_3-2D_2D_4-D_3^2+4D_4^2\big)\Big],\\
  \Phi_{2,3}(\mathbf{C})=& \Phi_{2,1}(S_1\mathbf{C})=\bar{\Phi}_{2,1}(\mathbf{C}), \quad \quad
  \Phi_{2,4}(\mathbf{C})= -\Phi_{2,2}(S_1\mathbf{C})=\bar{\Phi}_{2,2}(\mathbf{C}),\\
  \Phi_{2,5}(\mathbf{C})=& 8\nu_2\left[(c_{1,+}-c_{2,-})D_1-2c_{1,-}D_2-4c_{1,+}D_3+8(c_{1,-}+c_{2,+})D_4 \right],\\
  \Phi_{2,6}(\mathbf{C})=& 8\nu_2\left[-c_{1,-}D_1-(c_{1,+}+3c_{2,-})D_2+2(c_{1,-}-2c_{2,+})D_3+4c_{1,+}D_4 \right],\\
  \Phi_{2,7}(\mathbf{C})=& 8\nu_2\left[-(c_{1,+}+c_{2,-})D_1-4c_{2,+}D_2+(c_{1,+}+3c_{2,-})D_3-2c_{1,-}D_4 \right],\\
  \Phi_{2,8}(\mathbf{C})=& 8\nu_2\left[(c_{1,-}-2c_{2,+})D_1-(c_{1,+}-3c_{2,-})D_2-c_{1,-}D_3-(c_{1,+}-c_{2,-})D_4 \right].\\
  \end{aligned}
 \end{equation}

Conversely, it is less straightforward to obtain the explicit expression of $\mathscr{R}_3$.  
We start by determining a representative form for $\mathscr{R}_3$.
Similar to the quadratic case, from the Fredholm alternative, we solve \eqref{e:nf3} uniquely for  $\Phi_{3}$ and $\mathscr{R}_{3}$  
subject to
\begin{equation*}
\begin{aligned}
&\begin{pmatrix} \mathscr{R}_{3,1} \\ \mathscr{R}_{3,2}\end{pmatrix}\in\ker\left(\begin{pmatrix}\mathcal{D}^{\ad}+\rmi&0\\-1&\mathcal{D}^{\ad}+\rmi\end{pmatrix}\mid_{\mathbf{P}_3^2}\right), \quad \mathscr{R}_{3,3}=\bar{\mathscr{R}}_{3,1}, \quad \mathscr{R}_{3,4}=\bar{\mathscr{R}}_{3,2},\\
&
\mathscr{R}_{3,5}=\mathscr{R}_{3,6}=\mathscr{R}_{3,7}=0, \quad \mathscr{R}_{3,8}\in \ker\big((\mathcal{D}^{\ad})^4|_{\mathbf{P}_3}\big);\\
&\begin{pmatrix} \Phi_{3,1} \\ \Phi_{3,2}\end{pmatrix}\in\left(\ker\left(\begin{pmatrix}\mathcal{D}-\rmi&-1\\0&\mathcal{D}-\rmi\end{pmatrix}\mid_{\mathbf{P}_3^2}\right)\right)^\perp, \quad \Phi_{3,3}=\bar{\Phi}_{3,1}, \quad \Phi_{3,4}=\bar{\Phi}_{3,2},\\
&
\Phi_{3,6}=\mathcal{D}\Phi_{3,5}, \quad \Phi_{3,7}=\mathcal{D}^2\Phi_{3,5},\quad \Phi_{3,8}=\mathcal{D}^3\Phi_{3,5}, \quad \Phi_{3,5}\in \left( \ker\big((\mathcal{D})^4|_{\mathbf{P}_3}\big)\right)^\perp.
\end{aligned}
\end{equation*}
Similarly, we point out that 
\begin{equation}\
\label{e:ker3}
 \begin{aligned}
   \ker\big(\mathcal{D}^{\ad}|_{\mathbf{P}_3}\big)=\span\{&C_1\bar{C}_1D_1,(C_1\bar{C}_2-\bar{C}_1C_2)D_1,D_1^3,\\
  &D_1(2D_1D_3-D_2^2),3D_1(D_2D_3-D_1D_4)-D_2^3\},\\
 \ker\big((\mathcal{D}^{\ad}+\rmi)|_{\mathbf{P}_3}\big)=\span\{&C_1^2\bar{C}_1,C_1(C_1\bar{C}_2-\bar{C}_1C_2),C_1D_1^2,C_1(2D_1D_3-D_2^2),\\
 &(C_1D_2-C_2D_1)D_1,C_1(D_2D_3-3D_1D_4)+C_2(2D_1D_3-D_2^2)\}.
 \end{aligned}
\end{equation}
Based on \eqref{e:ker3} and the condition that $\mathscr{R}_{3,1\backslash 2}$ satisfies
\begin{equation*}
\begin{pmatrix} \mathscr{R}_{3,1} \\ \mathscr{R}_{3,2}\end{pmatrix}\in\ker\big(\begin{pmatrix}\mathcal{D}^{\ad}+\rmi&0\\-1&\mathcal{D}^{\ad}+\rmi\end{pmatrix}\mid_{\mathbf{P}_3^2}\big),
\end{equation*}
we obtain $\mathscr{R}_{3,1\backslash 2}$ in the following form,
\begin{equation}
\begin{aligned}
 \mathscr{R}_{3,1}=&C_1\big[\widetilde{\alpha}_7C_1\bar{C}_1+\widetilde{\alpha}_8(C_1\bar{C}_2-\bar{C}_1C_2)\big]+\widetilde{\alpha}_9C_1D_1^2+
 \widetilde{\alpha}_{10}D_1(C_2D_1-C_1D_2)+\\
& \widetilde{\alpha}_{11}C_1(2D_1D_3-D_2^2)+\widetilde{\alpha}_{12}\big[C_1(D_2D_3-3D_1D_4)+C_2(2D_1D_3-D_2^2)\big];\\
 \mathscr{R}_{3,2}=&C_1\big[\widetilde{\alpha}_1C_1\bar{C}_1+\widetilde{\alpha}_2(C_1\bar{C}_2-\bar{C}_1C_2)\big]+\widetilde{\alpha}_3C_1D_1^2+
 \widetilde{\alpha}_{4}D_1(C_2D_1-C_1D_2)+\\
 &\widetilde{\alpha}_5C_1(2D_1D_3-D_2^2)+\widetilde{\alpha}_6\big[C_1(D_2D_3-3D_1D_4)+C_2(2D_1D_3-D_2^2)\big]+\\
 &C_2\big[\widetilde{\alpha}_7C_1\bar{C}_1+\widetilde{\alpha}_8(C_1\bar{C}_2-\bar{C}_1C_2)\big]+\widetilde{\alpha}_9C_2D_1^2+
 \widetilde{\alpha}_{10}D_2(C_2D_1-C_1D_2)+\\
& \widetilde{\alpha}_{11}C_1(3D_1D_4-D_2 D_3)+\widetilde{\alpha}_{12}\big[2C_1(2D_3^2-3D_2D_4)+C_2(3D_1D_4-D_2D_3)\big].\\
 \end{aligned}
\end{equation}
Moreover, based on \eqref{e:ker3} and the condition that $\mathscr{R}_{3,5}=\mathscr{R}_{3,6}=\mathscr{R}_{3,7}=0$, and $ \mathscr{R}_{3,8}\in \ker\big((\mathcal{D}^{\ad})^4|_{\mathbf{P}_3}\big)$, $ \mathscr{R}_{3,8}$ takes the form
\begin{equation*} 
\begin{aligned}
 \mathscr{R}_{3,8}=&\sum_{j, k=1}^2\sum_{\ell=1}^4\widetilde{\beta}_{jk\ell}C_j\bar{C}_kD_\ell+D_1^2\sum_{j=1}^4\widetilde{\beta}_jD_j+
 \widetilde{\beta}_5D_1D_2^2+ \widetilde{\beta}_6D_2^3+ \widetilde{\beta}_7D_1D_2D_3+\\
 &\widetilde{\beta}_8(D_2^2D_3-2D_1D_3^2)+\widetilde{\beta}_9(D_2^2D_3-3D_1D_2D_4)+\widetilde{\beta}_{10}(2D_1D_3D_4-D_2^2D_4)+\\
 &\widetilde{\beta}_{11}(3D_1D_3D_4-D_2D_3^2)+\widetilde{\beta}_{12}(9D_2D_3D_4-9D_1D_4^2-4D_3^3),\\
\end{aligned}
\end{equation*}
where $\widetilde{\beta}_{224}=0$ and $\widetilde{\beta}_{124}+\widetilde{\beta}_{214}+3\widetilde{\beta}_{223}=0$.

We point out that this normal form inherits the symmetries of the original reduced ODE system. More specifically, $\Psi_{2\backslash 3}$ and $\mathscr{R}_{2\backslash 3}$ commute
with $S_1$ and $S_2$; see \cite[3.3]{harioo} for details.
As a result, the preservation of the reversibility $S_1$ further simplifies the cubic term $\mathscr{R}_{3}$. In fact, we have  
 \begin{equation}
  \widetilde{\alpha}_1,\widetilde{\alpha}_3,\widetilde{\alpha}_5, \widetilde{\alpha}_8,\widetilde{\alpha}_{10},\widetilde{\alpha}_{12}\in\R,\quad \widetilde{\alpha}_2,\widetilde{\alpha}_4,\widetilde{\alpha}_6, \widetilde{\alpha}_7,\widetilde{\alpha}_{9},\widetilde{\alpha}_{11}\in\rmi\R,
 \end{equation}
 and
 \begin{equation}
 \begin{aligned}
 \mathscr{R}_{3,8}=&D_1(\beta_{1}C_1\bar{C}_1+\beta_{2}C_2\bar{C}_2)+\rmi(C_1\bar{C}_2-\bar{C}_1C_2)(\beta_3D_1+\beta_4D_3)+\beta_5C_1\bar{C}_1D_3+\\
 &\beta_6D_2(C_1\bar{C}_2+\bar{C}_1C_2)+\beta_7[3(C_1\bar{C}_2+\bar{C}_1C_2)D_4-2C_2\bar{C}_2D_3]+\beta_8D_1D_2^2+\\
& D_1^2(\beta_9D_1+\beta_{10}D_3)+
\beta_{11}(D_2^2D_3-2D_1D_3^2)+\beta_{12}(D_2^2D_3-3D_1D_2D_4)+\\
&\beta_{13}(9D_2D_3D_4-9D_1D_4^2-4D_3^3),\\
\end{aligned}
 \end{equation}
 where all $\beta_j\in\R$, $j=1, 2, \dots, 13$.
 The expression for $\Phi_3$ is not required in the sequel and omitted.  \end{Proof}

Applying the composition of the linear and nonlinear normal form transformation, that is,
\[
 \mathbf{A}=\mathbf{T}(\varepsilon)\left(\mathbf{C}+\Psi_2(\mathbf{C})+\Psi_3(\mathbf{C})\right),
\]
the system \eqref{e:2dcredA} admits the normal form 
\begin{equation}
 \label{e:2dnfeq}
 \begin{aligned}
 \dot{\mathbf{C}}=\mathscr{L}(\varepsilon)\mathbf{C}+\mathscr{R}_3(\mathbf{C})+
  \caO\Big(|\varepsilon|^2\|\mathbf{C}\|+|\varepsilon|\|\mathbf{C}\|^2+\|\mathbf{C}\|^4\Big).
 \end{aligned}
\end{equation}
Truncating the normal form system at cubic terms in $\mathbf{C}$ and leading order in $\varepsilon$, yields 
\begin{equation}\label{e:2dnfc}
\dot{\mathbf{C}}=\mathscr{L}(\varepsilon)\mathbf{C}+\mathscr{R}_3(\mathbf{C}).
\end{equation}
The truncated normal form gains an extra rotational symmetry $R_\theta$, given by,
 \begin{equation}\label{e:symm}
  R_\theta(\mathbf{C})=(\rme^{\rmi \theta}C_1,\rme^{\rmi \theta}C_2,\rme^{-\rmi \theta}\bar{C}_1,\rme^{-\rmi \theta}\bar{C}_2,D_1,D_2,D_3,D_4).
 \end{equation}
\begin{Remark}
This additional symmetry $R_\theta$ results from the form of the linear term in the original 8th-order ODE \eqref{e:2dcredA} and our particular choice of the normal form transformation; see \cite[Chapter 3]{harioo} for details. Moreover, this additional symmetry $R_\theta$ fails to hold for the full normal form system \eqref{e:2dnfeq} while the reversibility $S_1$ and the symmetry $S_2$ hold.
\end{Remark}

\subsection{Construction of extended pearled solutions}\label{ss:24}
We adapt the techniques of  \cite[Section 3.1, 4.1]{iooper_1993}, employing rescalings and the implicit function theorem to construct periodic solutions to the normal form system \eqref{e:2dnfeq}, which correspond to extended pearled solutions of the flat-bilayer system \eqref{e:2sFCH}.

Restricting the truncated normal form system \eqref{e:2dnfc}  to the subspace 
\[\widetilde{\R^4}:=\{(C_1, C_2, \bar{C}_1, \bar{C}_2, 0,0,0,0)\mid C_1, C_2\in\C\},\]
yields 
the 1:1 resonant normal form; see \cite{imd_1989, iooper_1993, harioo},
\begin{equation}
\label{e:rhopf}
\begin{cases}
\dot{C_1}=\rmi (1+\omega_1\varepsilon) C_1 + C_2+\rmi C_1\big[\alpha_7C_1\bar{C}_1+\alpha_8\rmi(C_1\bar{C}_2-\bar{C}_1C_2)\big],\\
\dot{C_2}=\rmi (1+\omega_1\varepsilon) C_2 +\rmi C_2\big[\alpha_7C_1\bar{C}_1+\alpha_8\rmi(C_1\bar{C}_2-\bar{C}_1C_2)\big]
+ C_1\left[\omega_2\varepsilon+\alpha_1C_1\bar{C_1}+\rmi \alpha_2(C_1\bar{C_2}-\bar{C_1}C_2)\right].
\end{cases}
\end{equation}

\begin{Remark}
The 4th-order system in \cite{iooper_1993} admits a very general normal form in which the even-order terms automatically vanishes. We can not make this generalization here since the invariance of the pearling modes is not guaranteed when we push the normal form to high orders.
\end{Remark}

The construction of the extended pearled solutions relies crucially on two  properties of the 1:1 resonant normal form. First,
the 1:1 resonant normal form \eqref{e:rhopf} admits  two \textit{first integrals},
\[
K=\frac{\rmi}{2}(C_1\bar{C}_2-\bar{C}_1C_2), \quad H=|C_2|^2-\left[\frac{\alpha_1}{2}|C_1|^2-\left(\omega_2\varepsilon+2\alpha_2K\right)\right]|C_1|^2,
\]
and as a consequence may be reduced to a 2nd-order ODE in the variable $u_1:=|C_1|^2$.  The pearled morphologies we seek correspond to periodic
solutions of \eqref{e:rhopf},  which are temporal equilibrium of the 2nd-order ODE for $u_1$. As a second point, the 1:1 resonant normal form is autonomous in the pearling 
modes $(C_1, C_2, \bar{C}_1,\bar{C}_2)$ and thus the subspace $\widetilde{\R^4}$ is invariant under the truncated normal form flow \eqref{e:2dnfc}. 
In this sense, the pearling modes 
$(C_1, C_2, \bar{C}_1, \bar{C}_2)$  and the meandering modes $(D_1, D_2, D_3, D_4)$ exhibit a weak coupling. Accordingly, we 
anticipate that structures in the 1:1 resonant normal form
will persist in the full normal form system. 

A complication in the persistence argument arises through the degeneracy of the particular 1:1 resonant normal form studied here. 
The two parameters, $\alpha_1$ and $\omega_2$, characterize the 1:1 resonant normal form, where $\alpha_1$ is the coefficient of $C_1^2\bar{C}_1$ in 
the second entry of the cubic normal form. As shown in Lemma \ref{l:32}, for the pearling problem we have
\[
\alpha_1=0,
\]
which leads to a \textit{degenerate 1:1 resonance}.  For uniformity of notation and the sign consistency with the linear stability condition in 
\cite{doelmanpromislow_2013}, we also introduce 
\begin{equation}\label{e:alpha_0}
\alpha_0:=-\omega_2=-\mu_2=\frac{1}{4\lambda_0^2}\int_\R\left(W^\tprime(u_0)v_0-\eta_dW^\pprime(u_0)\right)\psi_0^2\rmd r.
\end{equation}
With these modifications, we rename the degenerate system  \emph{the pearling normal form} (PNF) system,
\begin{equation}\label{e:PNF}
\hskip -0.2in
\begin{cases}
\dot{C_1}&=\rmi (1+\omega_1\varepsilon) C_1 + C_2+\rmi C_1\big[\alpha_7C_1\bar{C}_1+\alpha_8\rmi(C_1\bar{C}_2-\bar{C}_1C_2)\big],\\
\dot{C_2}&=\rmi (1+\omega_1\varepsilon) C_2 +\rmi C_2\big[\alpha_7C_1\bar{C}_1+\alpha_8\rmi(C_1\bar{C}_2-\bar{C}_1C_2)\big]
+ C_1\left[-\alpha_0\varepsilon+\rmi \alpha_2(C_1\bar{C_2}-\bar{C_1}C_2)\right],
\end{cases}
\end{equation}
For this degenerate case the persistence issue is a singular perturbation problem; removing the singularity
requires two novel proper rescalings.  After the first scaling, we construct a Poincare map, which is well-defined for sufficiently 
small system parameters, including the zeroes.  However the base of the transverse hyper-plane in the Poincare map consists of eigenvectors. As the 
system parameters approach zero, the degeneracy of eigenvalues results in the coalescence of the eigenvectors, which we overcome via 
a second rescaling. The persistence follows from an implicit-function-theorem argument.


The existence results for periodic solutions of the PNF system are summarized in the following lemma, where
for convenience, we assume $\varepsilon>0$ and introduce the rescaled first integral $\kappa:=\varepsilon^{-3/2}K$.
\begin{Lemma}[degenerate 1:1 resonance]\label{l:29}
For fixed $\eta_1$, $\eta_2$, $\gamma\in\R$ and a non-degenerate double-well potential $W$, there exist $\varepsilon_0$, $\kappa_0>0$ such that, for every $\varepsilon\in(0,\varepsilon_0]$, the PNF system \eqref{e:PNF}
admits a \textit{degenerate 1:1 resonance}, characterized by $\alpha_0$, defined in \eqref{e:alpha_0}. More specifically, we have
\begin{enumerate}
\item[(\rmnum{1})] For $\alpha_0<0$ , the PNF system \eqref{e:PNF} has no periodic solutions except for the trivial equilibrium.
\item[(\rmnum{2})] For $\alpha_0>0$, the PNF system \eqref{e:PNF} possesses a family of periodic orbits $(C_1^p, C_2^p,\bar{C}_1^p,\bar{C}_2^p)$,
parameterized by $\kappa\in[-\kappa_0,\kappa_0]$. In fact, the family of periodic orbits is smooth in terms of small $\sqrt{\varepsilon}$ and $\sqrt{|\kappa|}$ except for $\kappa=0$, admitting the form
 \begin{equation}\label{e:peri}
 \begin{aligned}
   C_1^p(t,\theta;\sqrt{\varepsilon}, \sqrt{|\kappa|})=&\sqrt{\varepsilon|\kappa|} r_1\rme^{\rmi(\omega t+\theta)},\\
   C_2^p(t,\theta; \sqrt{\varepsilon},\sqrt{|\kappa|})=&\sgn(\kappa)\rmi \varepsilon\sqrt{|\kappa|} r_2\rme^{\rmi(\omega t+\theta)},
 \end{aligned}
 \end{equation}
 where
 \begin{equation}\label{e:r1}
 r_1(\sqrt{\varepsilon},\sqrt{|\kappa|})=
 (\alpha_0-2\alpha_2\sqrt{\varepsilon}\kappa)^{-1/4},
 \end{equation}
 and
 \begin{equation}
  r_2=\frac{1}{r_1},\quad \omega=1+\omega_1\varepsilon+\sgn(\kappa)\sqrt{\varepsilon} r_2^2+\alpha_7\varepsilon|\kappa| r_1^2+2\alpha_8\varepsilon^{3/2}\kappa,\quad \theta\in\R/[0,2\pi].
 \end{equation}
\end{enumerate}
\end{Lemma}
\begin{Proof}
 Under the polar coordinate change
\[
C_1=\widetilde{r}_1\rme^{\rmi(1+\omega_1\varepsilon+\theta_1)}, \quad C_2=\widetilde{r}_2\rme^{\rmi(1+\omega_1\varepsilon+\theta_2)}, \quad u_1=\widetilde{r}_1^2, \quad u_2=\widetilde{r}_2^2,
\]
the PNF system \eqref{e:PNF} becomes
\begin{equation}\label{e:rPNF}
\begin{cases}
&\left(\frac{\rmd u_1}{\rmd t}\right)^2=4f(u_1),\\
&\frac{\rmd (\theta_2-\theta_1)}{\rmd t}=-K(u_1u_2)^{-1}f^\prime(u_1),\\
&\frac{\rmd \theta_1}{\rmd t}=Ku_1^{-1}+\alpha_7u_1+2\alpha_8K,\\
&\frac{\rmd u_2}{\rmd t}=\left(\alpha_7u_1+2\alpha_8K\right)\frac{\rmd u_1}{\rmd t},
\end{cases}
\end{equation}
where $f(u_1)=(-\alpha_0\varepsilon+2\alpha_2K)u_1^2+Hu_1-K^2$. We observe that a double root of
\[
f(u_1)=0,
\]
corresponds to an equilibrium of the ODE 
\begin{equation}\label{e:u1}
\left(\frac{\rmd u_1}{\rmd t}\right)^2=4f(u_1),
\end{equation}
which corresponds to a periodic solution in the PNF system \eqref{e:PNF}. 
We apply the rescaling
\[
u_1=\varepsilon v_1, \quad K=\varepsilon^{3/2}\kappa, \quad H=\varepsilon^2 h,
\]
to $f(u_1)=0$ and have
\[
(-\alpha_0+2\alpha_2\kappa\sqrt{\varepsilon})v_1^2+hv_1-\kappa^2=0
\]
which admits a double root if and only if  
\begin{equation}\label{e:G}
\left(2(-\alpha_0+2\alpha_2\kappa\sqrt{\varepsilon})v_1+h,(-\alpha_0+2\alpha_2\kappa\sqrt{\varepsilon})v_1^2+\kappa^2\right)^T=0.
\end{equation}
If $\alpha_0<0$, \eqref{e:G} admits only the trivial solution for small $\varepsilon$ and $k$. If $\alpha_0>0$, then we can solve $v_1$ and $h$ in terms of $\varepsilon$
and $\kappa$. In fact, we have, for sufficiently small $\varepsilon$ and $k$, 
\begin{equation}
\begin{cases}
v_1(\varepsilon, \kappa)=\frac{|\kappa|}{\sqrt{\alpha_0-2\alpha_2\sqrt{\varepsilon}\kappa}},\\
h(\varepsilon, \kappa)=2(\alpha_0-2\alpha_2\sqrt{\varepsilon}\kappa)v_1.
\end{cases}
\end{equation}
We conclude our proof by letting $r_1=\sqrt{v_1/|\kappa|}$.
\end{Proof}

In the sequel we assume 
$\alpha_0>0$ and $\kappa\geq 0$. The analysis of the case $\kappa<0$ differs only by a sign change. To demonstrate the persistence of the periodic solutions of the PNF system in the full normal form system \eqref{e:2dnfeq}, it is necessary to remove the singular nature of the bifurcation. 
To this end we apply the rescaling
\begin{equation}
 \label{e:res}
 \mathbf{C}=\sqrt{\varepsilon\kappa}\widetilde{\mathbf{C}},
\end{equation}
to the normal form system \eqref{e:2dnfeq}, obtaining a new ODE system
\begin{equation}
 \label{e:ODE}
 \begin{aligned}
 \dot{\mathbf{C}}=\mathscr{L}(\varepsilon)\mathbf{C}+\varepsilon\kappa\mathscr{R}_3(\mathbf{C})
 +\varepsilon\caO\big(\varepsilon\|\mathbf{C}\|+\sqrt{\varepsilon\kappa}\|\mathbf{C}\|^2+\sqrt{\varepsilon\kappa^3}\|\mathbf{C}\|^4\big),
 \end{aligned}
\end{equation}
where we have dropped the ``tilde'' notation on ${\bf C}$.
To simplify the proof of the persistence we introduce the new small parameter $\zeta$ for which $\zeta=0$ corresponds to the cubic truncation, while $\zeta=\varepsilon$ corresponds to the full normal form. Specifically, we study the system
\begin{equation}
 \label{e:gODE}
 \begin{aligned}
 \dot{\mathbf{C}}=\mathbf{F}(\mathbf{C})
 +\zeta\caO\big(\varepsilon\|\mathbf{C}\|+\sqrt{\varepsilon\kappa}\|\mathbf{C}\|^2+\sqrt{\varepsilon\kappa^3}\|\mathbf{C}\|^4\big),
 \end{aligned}
\end{equation}
where $\mathbf{F}(\mathbf{C}):=
\mathscr{L}(\varepsilon)\mathbf{C}+\varepsilon\kappa\mathscr{R}_3(\mathbf{C})$.
The following Proposition, taken from  \cite{iooper_1993}, greatly simplifies the construction.
\begin{Proposition}\label{p:tpt}
An orbit of an autonomous reversible system is periodic and reversible 
if and only if there exist two different fixed points on this orbit
with respect to the reversibility. 
\end{Proposition}

We lift the scalars $r_1$ and $r_2$, introduced in Lemma \ref{l:29}, which serve as the base point for the periodic solutions of the PNF system, to a vector
in the $8$ dimensional space, defining \textit{the base point} $\mathbf{r}$ as
\[
 \mathbf{r}(\sqrt{\varepsilon},\sqrt{\kappa})=(r_1,\rmi\sqrt{\varepsilon} r_2,r_1,-\rmi\sqrt{\varepsilon}r_2,0,0,0,0)^T,
\]
and obeserve that the periodic solution $R_{\omega t}\mathbf{r}$ to the system \eqref{e:gODE} when $\zeta=0$
has two fixed points under reversibility, that is, 
\[
 S_1\mathbf{r}=\mathbf{r},\quad S_1R_{\pi}\mathbf{r}=R_{\pi}\mathbf{r}.
\]
Here we recall that 
\begin{equation*}
  \begin{aligned}
  S_1(C_1,C_2,\bar{C}_1,\bar{C}_2, D_1,D_2,D_3, D_4)^T&=(\bar{C}_1,-\bar{C}_2, C_1, -C_2, D_1, -D_2, D_3, -D_4)^T,\\
 R_{\theta}(C_1,C_2,\bar{C}_1,\bar{C}_2, D_1,D_2,D_3, D_4)^T&=(\rme^{\rmi\theta}C_1, \rme^{\rmi\theta}C_2, \rme^{-\rmi\theta}\bar{C}_1, \rme^{-\rmi\theta}\bar{C}_2, D_1,D_2,D_3, D_4)^T.
  \end{aligned}
 \end{equation*}
We assign two transversal hyper-planes, $H_1$ and $H_2$, respectively to $\mathbf{r}$ and $R_{\pi}\mathbf{r}$, given as follows.
\begin{equation*}
H_1=\{\mathbf{C}\in\widetilde{\R^4}\times\R^4\mid S_1\mathbf{C}={\bf C}\},\quad
H_2=\{\mathbf{C}\in\widetilde{\R^4}\times\R^4\mid (\mathbf{C}-R_\pi\mathbf{r})\cdot R_\pi\mathbf{G}\mathbf{r}=0\},
\end{equation*}
where $\mathbf{G}$ is the infinitesimal generator of the group $R_\theta$ and ``$\cdot$" represents the Euclidean inner product.
It is then not hard to see that, for the rescaled system \eqref{e:gODE}, there exists a smooth Poincar\'{e} map, denoted as $\Pi$,  from an open neighborhood of the base point $\mathbf{r}$ in $H_1$,
$N(\mathbf{r},H_1)$, into one of $R_\pi\mathbf{r}$ in $H_2$, $N(R_\pi\mathbf{r}, H_2)$. More specifically, we have
\begin{equation}
\Pi(\mathbf{C},\sqrt{\varepsilon}, \sqrt{\kappa},\zeta): N(\mathbf{r},H_1)\times [0,\sqrt{\varepsilon_0}]\times[0,\sqrt{\kappa_0}]\times [-\zeta_0,\zeta_0]\rightarrow N(R_\pi\mathbf{r}, H_2).
\end{equation}
Meanwhile, there is also a smooth ``arrival time" map
\begin{equation}
T(\mathbf{C},\sqrt{\varepsilon}, \sqrt{\kappa},\zeta): N(\mathbf{r},H_1)\times [0,\sqrt{\varepsilon_0}]\times[0,\sqrt{\kappa_0}]\times [-\zeta_0,\zeta_0]\rightarrow \R.
\end{equation}

According to Proposition \ref{p:tpt}, any point in $H_1\cap\Rg(\Pi)$ corresponds to a periodic orbit of the system \eqref{e:gODE}, vice versa.
To further analyze the Poincar\'{e} map, we first linearize the system \eqref{e:gODE} around the periodic orbit $R_{\omega t}\mathbf{r}$.  
We introduce the change of variables local to $R_{\omega t}\mathbf{r}$,
\[
 \mathbf{C}=R_{\omega t}(\mathbf{r}+\mathbf{q}),
\]
and study the flow of $\mathbf{q}$ instead, that is,
\begin{equation}\label{e:gODE1}
 \frac{\mathrm{d}\mathbf{q}}{\mathrm{d}t}=\mathbf{F}(\mathbf{r}+\mathbf{q})-\mathbf{F}(\mathbf{r})-\omega \mathbf{G}\mathbf{q}
 +\caO\big(\sqrt{\varepsilon}(\sqrt{\varepsilon}+\sqrt{\kappa})|\zeta|\big).
\end{equation}
Linearizing the system \eqref{e:gODE1} at $\mathbf{q}=0$ yields the following system
\begin{equation}\label{e:gODEn}
  \dot{\mathbf{q}}=\mathbf{H}\mathbf{q}
 +\caO\big(\varepsilon\kappa\|\mathbf{q}\|^2
  +\sqrt{\varepsilon}(\sqrt{\varepsilon}+\sqrt{\kappa})|\zeta|\big),
\end{equation}
where $\mathbf{H}:=\nabla_{\mathbf{C}}\mathbf{F}(\mathbf{r})-\omega \mathbf{G}$. 
\begin{Remark}
The reversibility holds within the truncated system
$ \dot{\mathbf{q}}=\mathbf{F}(\mathbf{r}+\mathbf{q})-\mathbf{F}(\mathbf{r})-\omega \mathbf{G}\mathbf{q}$,
but not within the full ODE system about $\mathbf{q}$,  since the rotational symmetry $R_{\omega t}$ and the reversibility $S_1$ do not commute. 
As a result, we have
\[
S_1{\bf H}=-{\bf H}S_1.
\]
\end{Remark}

The next step is to obtain the eigenvalues and corresponding eigenmodes of $\mathbf{H}$.
We note that $\mathbf{H}$ is block diagonal. 
The upper diagonal block ${\bf H}_1$ of $\mathbf{H}$ is of the form
\begin{equation*}
\begin{aligned}
{\bf H}_1 =&  \begin{pmatrix}
  0 & 1 & 0 & 0 \\
  0 & 0 & 0 & 0 \\
  0 & 0 & 0 & 1 \\
  0 & 0 & 0 & 0
 \end{pmatrix}+\sqrt{\varepsilon}r_2^2
 \begin{pmatrix}
  -\rmi & 0 & 0 & 0 \\
  0 & -\rmi & 0 & 0 \\
  0 & 0 & \rmi & 0 \\
  0 & 0 & 0  & \rmi
 \end{pmatrix}-\varepsilon
 \begin{pmatrix}
  0 & 0 & 0 & 0 \\
  \alpha_0 & 0 & 0 & 0 \\
  0 & 0 & 0 & 0 \\
  0 & 0 & \alpha_0 & 0
 \end{pmatrix}-\varepsilon^2\kappa r_2^2
  \begin{pmatrix}
  0 & 0 & 0 & 0 \\
 \alpha_8 & 0 & \alpha_8 & 0 \\
  0 & 0 & 0 & 0 \\
  \alpha_8 & 0 & \alpha_8 & 0
 \end{pmatrix}+
\\
 &\varepsilon \kappa r_1^2\begin{pmatrix}
  \rmi\alpha_7 & \alpha_8 & \rmi\alpha_7 & -\alpha_8 \\
  0 & -\rmi\alpha_2 & 0 & \rmi\alpha_2 \\
  -\rmi\alpha_7 & -\alpha_8  & -\rmi\alpha_7 & \alpha_8  \\
0 & -\rmi\alpha_2 & 0 &  \rmi\alpha_2
 \end{pmatrix}+\varepsilon^{3/2}\kappa
\begin{pmatrix}
\rmi\alpha_8 & 0 & \rmi\alpha_8 &0 \\
3\alpha_2-\alpha_7 &\rmi\alpha_8 & \alpha_2-\alpha_7 & -\rmi \alpha_8\\
  -\rmi\alpha_8 & 0 & -\rmi\alpha_8 & 0 \\
 \alpha_2-\alpha_7 & \rmi \alpha_8 & \alpha_2-\alpha_7 &-\rmi\alpha_8
 \end{pmatrix}\\
 =& \begin{pmatrix}
  0 & 1 & 0 & 0 \\
  0 & 0 & 0 & 0 \\
  0 & 0 & 0 & 1 \\
  0 & 0 & 0 & 0
 \end{pmatrix}+\sqrt{\alpha_0\varepsilon}
 \begin{pmatrix}
  -\rmi & 0 & 0 & 0 \\
  0 & -\rmi & 0 & 0 \\
  0 & 0 & \rmi & 0 \\
  0 & 0 & 0  & \rmi
 \end{pmatrix}+\caO(\varepsilon).
  \end{aligned}
\end{equation*}
It is straightforward to see that
\[
 \mathbf{H}\mathbf{G}\mathbf{r}=0,\quad \mathbf{H}\frac{\partial \mathbf{r}}{\partial \sqrt{\kappa}}=\frac{\partial \omega}{\partial \sqrt{\kappa}}\mathbf{G}\mathbf{r},
\]
where $\mathbf{r}:=\mathbf{r}(\sqrt{\varepsilon},\sqrt{\kappa})$ and $\omega:=\omega(\sqrt{\varepsilon},\sqrt{\kappa})$. As a result, $0$ is an eigenvalue to the upper diagonal block ${\bf H}_1$ with algebraic multiplicity $2$. A direct calculation then shows that 
the determinant of ${\bf H}_1$ is 
\[
\det(\lambda-{\bf H}_1)=\lambda^4+4\varepsilon r_2^4\lambda^2,
\]
which indicates that the other two eigenvalues of ${\bf H}_1$ are
\[
 \pm\lambda_{1}=\pm2\rmi\sqrt{\varepsilon}r_2^2
 =\pm2\rmi\sqrt{\alpha_0\varepsilon-2\alpha_2\varepsilon^{3/2}\kappa}=\pm2\rmi\sqrt{\alpha_0\varepsilon}+\caO(\varepsilon\kappa),
\]
with associated eigenvectors ${\bf r}^{\pm}_1$ satisfying
\[
 \mathbf{H}{\bf r}_1^+=\lambda_1{\bf r}_1^+,\quad  S_1\mathbf{r}_1^+=\mathbf{r}_1^-.
\]
More specifically, a nonzero vector of cofactors of any row of ${\bf H}_1-\lambda_1$ is an eigenvector with respect to $\lambda_1$ since the algebraic multiplicity of $\lambda_1$ is 1. We then let $\mathbf{r}_1^+=({\bf r}_{1,1}^+,0,0,0,0)^T$, where ${\bf r}_{1,1}^+$ is the vector of cofactors of the second row of ${\bf H}_1-\lambda_1$ after an $\varepsilon^{3/2}\kappa$-rescaling, that is,
\[
\mathbf{r}_{1,1}^+=\begin{pmatrix} \alpha_7\sqrt{\varepsilon}\kappa r_1^4\\ \rmi\alpha_7\varepsilon\kappa r_1^2 \\
  2-\alpha_7\sqrt{\varepsilon}\kappa r_1^4 \\2\rmi\sqrt{\varepsilon} r_2^{2}+\rmi\alpha_7\varepsilon\kappa r_1^2\end{pmatrix}=\begin{pmatrix} 0\\ 0 \\  2 \\0\end{pmatrix}+\caO(\sqrt{\varepsilon}).
\]

The lower block ${\bf H}_2$ of ${\bf H}$ is of the form
\begin{equation*}
 {\bf H}_2=\begin{pmatrix}
  0 & 1 & 0 & 0 \\
  0 & 0 & 1 & 0 \\
  0 & 0 & 0 & 1 \\
  c & 0 & b & 0
 \end{pmatrix},
\end{equation*}
where
\[
b=\varepsilon \left[\omega_4+\beta_5 \kappa r_1^2-2\beta_7\varepsilon\kappa r_2^2+2\beta_4\sqrt{\varepsilon}\kappa\right] , \quad
c=\varepsilon\left(\beta_0+\beta_1\kappa r_1^2+\beta_2\varepsilon\kappa r_2^2+ 2\beta_3\sqrt{\varepsilon}\kappa\right).
\]
Here we use the fact that $\beta_0=\omega_3$.
Noting that the characteristic polynomial of ${\bf H}_2$ is 
\[
\lambda^4-b\lambda^2-c=0,
\]
we conclude that ${\bf H}_2$ has nonzero eigenvalues if and only if $c\neq0$, which can be guaranteed
by further assuming that $\beta_0\neq 0$ for $\varepsilon$ and $\kappa$ sufficiently small.

We summarize our assumptions on system parameters in the following hypothesis.
\begin{Hypothesis}[generic and non-degeneracy condition]\label{h:gen}
We assume that
\begin{equation}\label{e:gen}
 \alpha_0 >0,\quad \beta_0\neq 0, \quad \varepsilon>0, \quad \kappa>0.
\end{equation}
\end{Hypothesis}

Under this non-degeneracy assumption, we have, for sufficiently small $\kappa$ and $\varepsilon$,
\[
b^2+4c\neq 0,
\]
which implies that ${\bf H}_2$ admits four distinct nonzero eigenvalues $\pm\lambda_2$ and $\pm\lambda_3$ of order $\sqrt[4]{\varepsilon}$ and with associated eigenvectors $\mathbf{r}_2^\pm$ and $\mathbf{r}_3^\pm$ satisfying
 \[
  \mathbf{H}\mathbf{r}_2^+=\lambda_2\mathbf{r}_2^+,\quad \mathbf{H}\mathbf{r}_3^+=\lambda_3\mathbf{r}_3^+,\quad
  S_1\mathbf{r}_2^+=\mathbf{r}_2^-,\quad S_1\mathbf{r}_3^+=\mathbf{r}_3^-.
 \]
 More specifically, we choose
 \begin{equation*}
 \begin{aligned}
& \lambda_2=\left(\frac{b+\sqrt{b^2+4c}}{2}\right)^{1/2}=\sqrt[4]{\beta_0\varepsilon}+\caO(\varepsilon^{3/4}+\sqrt[4]{\varepsilon}\kappa), \qquad \mathbf{r}_2^+=(0,0,0,0,1, \lambda_2,\lambda_2^2,\lambda_2^3)^T,\\
 & \lambda_3=\left(\frac{b-\sqrt{b^2+4c}}{2}\right)^{1/2}=\rmi\sqrt[4]{\beta_0\varepsilon}+\caO(\varepsilon^{3/4}+\sqrt[4]{\varepsilon}\kappa), \qquad \mathbf{r}_3^+=(0,0,0,0,1, \lambda_3,\lambda_3^2,\lambda_3^3)^T.
  \end{aligned}
 \end{equation*}

Based on the spectral information about $\mathbf{H}$ we collected, we denote 
\begin{equation*}
\begin{aligned}
\mathbf{r}_0(\sqrt{\varepsilon},\sqrt{\kappa})&=\sqrt[4]{\alpha_0}\frac{\partial \mathbf{r}}{\partial \sqrt{\kappa}}=(1,0,1,0,0,0,0,0)^T+\caO(\sqrt{\varepsilon}),\\
\mathbf{r}_j(\sqrt{\varepsilon},\sqrt{\kappa})&=\mathbf{r}_j^++\mathbf{r}_j^-, \quad j=1,2,3,\\
\widetilde{\mathbf{r}}_1(\sqrt{\varepsilon},\sqrt{\kappa})&=\mathbf{r}_1^+-\mathbf{r}_1^-,\\
\widetilde{\mathbf{r}}_j(\sqrt[4]{\varepsilon},\sqrt{\kappa})&=\mathbf{r}_j^+-\mathbf{r}_j^-, \quad j=2,3.\\
\end{aligned}
\end{equation*}
We note that every ${\bf r}_j$, $j=0,1,2,3$ is a smooth with respect to its arguments. In particular, even though $\lambda_2$ and $\lambda_3$ are of order $\sqrt[4]{\varepsilon}$,
\begin{equation*}
\begin{aligned}
{\bf r}_j={\bf r}_j^++{\bf r}_j^-=2(0,0,0,0,1,0,\lambda_j^2,0)^T=2(0,0,0,0,1,0,0,0)^T+\caO(\sqrt{\varepsilon}), \quad j=2,3,\\
\end{aligned}
\end{equation*}
is smooth in terms of $\sqrt{\varepsilon}$.
We characterize the two transversal hyperplanes, $H_1$ and $H_2$, by the eigenvectors, that is,
$H_1=\mathbf{r}+\widetilde{H}_1$  and $H_2=R_\pi \mathbf{r}+\widetilde{H}_2$, where
\begin{equation*}
  \widetilde{H}_1=\span\{\mathbf{r}_0,\mathbf{r}_j\mid j=1,2,3\}, \quad
  \widetilde{H}_2=\span\{R_\pi \mathbf{r}_0,R_\pi\mathbf{r}_j^\pm\mid j=1,2,3\}.
\end{equation*}
We also parameterize $\mathbf{q}_1\in \widetilde{H}_1$ and $\mathbf{q}_2\in \widetilde{H}_2$ by
\begin{equation}\label{e:loor}
 \begin{aligned}
  \mathbf{q}_1&=\sum_{j=0}^3q_{1,j}\mathbf{r}_j,\\
 \mathbf{q}_2&=\sum_{j=0}^3q_{2,j}R_\pi\mathbf{r}_j+\sum_{j=1}^3\widetilde{q}_{2,j}R_\pi\widetilde{\mathbf{r}}_j
 \end{aligned}
\end{equation}
where we denote $q_1=(q_{1,0},q_{1,1},q_{1,2},q_{1,3})$ and $q_2=(q_{2,0},q_{2,1},q_{2,2},q_{2,3},\widetilde{q}_{2,1},\widetilde{q}_{2,2},\widetilde{q}_{2,3})$.
\begin{Remark}\label{r:coal}
The parameterization \eqref{e:loor} is singular at $\varepsilon=0$, since the eigenvalues coalesce. 
More specifically, when $\varepsilon=0$,  multiple eigenvectors collapse into one, that is, 
\[2{\bf r}_0(0,\sqrt{\kappa})=\mathbf{r}_1(0,\sqrt{\kappa}), \quad \mathbf{r}_2^+(0,\sqrt{\kappa})=\mathbf{r}_2^-(0,\sqrt{\kappa})=\mathbf{r}_3^+(0,\sqrt{\kappa})=\mathbf{r}_3^-(0,\sqrt{\kappa}).\]
\end{Remark}

Therefore, with this singular parameterization \eqref{e:loor}, we rewrite the Poincar\'{e} map and the arrival time map as follows.
\begin{equation*}
\begin{aligned}
&\begin{matrix}
\widetilde{\Pi}: &[-q_0, q_0]^4\times [0,\sqrt{\varepsilon}_0]\times[0,\sqrt{\kappa}_0]\times [-\zeta_0,\zeta_0]&\longrightarrow &N(R_\pi\mathbf{r}, H_2),\\
&(q_1,\sqrt{\varepsilon},\sqrt{\kappa},\zeta)
&\longmapsto&\Pi(\mathbf{r}+\sum_{j=0}^3q_{1,j}\mathbf{r}_j,\sqrt{\varepsilon}, \sqrt{\kappa},\zeta),\\
\end{matrix}\\
&\begin{matrix}
\widetilde{T}: &[-q_0, q_0]^4\times [0,\sqrt{\varepsilon}_0]\times[0,\sqrt{\kappa}_0]\times [-\zeta_0,\zeta_0]&\longrightarrow &N(R_\pi\mathbf{r}, H_2),\\
&(q_1,\sqrt{\varepsilon},\sqrt{\kappa},\zeta)
&\longmapsto&T(\mathbf{r}+\sum_{j=0}^3q_{1,j}\mathbf{r}_j,\sqrt{\varepsilon}, \sqrt{\kappa},\zeta).\\
\end{matrix}\\
\end{aligned}
\end{equation*}
Note that $\widetilde{\Pi}$ and $\widetilde{T}$ are smooth in terms of their arguments in the domain due to the fact that every ${\bf r}_j$ is smooth in terms of $\sqrt{\varepsilon}$ and $\sqrt{\kappa}$.
Moreover, according to the coalescence of eigenvectors when $\varepsilon=0$ in Remark \ref{r:coal}, it is not hard to verify that 
\begin{equation}\label{e:Tep}
\widetilde{T}(q_1,\sqrt{\varepsilon},\sqrt{\kappa},\zeta)=\pi+\caO\left(\sqrt{\varepsilon}(1+\sqrt{\kappa}+|\zeta|+\|q_1\|)\right)=\frac{\pi}{\omega}+\caO\left(\sqrt{\varepsilon}(|\zeta|+\|q_1\|)\right).
\end{equation}
Applying the variation of constant formula to \eqref{e:gODEn} and the parameterization of ${\bf q}_1$ in \eqref{e:loor}, together with the equality \eqref{e:Tep}, we have that
\begin{equation}\label{e:poex}
 \begin{aligned}
 \widetilde{ \Pi}(q_1,\sqrt{\varepsilon},\sqrt{\kappa},\zeta)=&R_{\omega \widetilde{T}}\left(\mathbf{r}+
  \mathrm{e}^{\mathbf{H} \widetilde{T}}\mathbf{q}_1\right)+\mathcal{O}(\varepsilon\kappa\|q_1\|^2+\sqrt{\varepsilon}(\sqrt{\varepsilon}+\sqrt{\kappa})|\zeta|)\\
  =&R_{\pi}\mathbf{r}+R_{\pi}\exp(\mathbf{H} \frac{\pi}{\omega})\mathbf{q}_1+
  \omega\left(\sum_{j=0}^3\widetilde{T}_jq_{1,j}\right)R_{\pi}\mathbf{G}\mathbf{r}+
  \caO\left(\sqrt{\varepsilon}(|\zeta|+\|q_1\|^2)\right)\\
 =&R_{\pi}\mathbf{r}+q_{1,0}R_{\pi}\mathbf{r}_0+\sum_{j=1}^3q_{1,j}\left[
  \cosh(\lambda_j\frac{\pi}{\omega})R_\pi\mathbf{r}_j+\sinh(\lambda_j\frac{\pi}{\omega})
  R_\pi\widetilde{\mathbf{r}}_j\right]+\\
 &\omega\left[\left(\widetilde{T}_0 +
  \alpha_0^{1/4}\frac{\pi}{\omega^2}\frac{\partial \omega}{\partial \sqrt{\kappa}}\right)q_{1,0}+
 \sum_{j=1}^3\widetilde{T}_jq_{1,j}\right]R_{\pi}\mathbf{G}\mathbf{r}+\caO\left(\sqrt{\varepsilon}(|\zeta|+\|q_1\|^2)\right),
 \end{aligned}
\end{equation}
where $\cosh$ and $\sinh$ take their natural analytic extension onto $\C$. Moreover, we have
\[
\widetilde{T}_0:=\frac{\partial \widetilde{T}}{\partial q_{1,0}}(0,\sqrt{\varepsilon},\sqrt{\kappa},0), \quad 
 \widetilde{T}_j:=\frac{\partial \widetilde{T}}{\partial q_{1,j}}(0,\sqrt{\varepsilon},\sqrt{\kappa},0), j=1,2,3.
\]
Noting that $R_{\pi}\mathbf{G}\mathbf{r}$ is transverse to $H_2$ and $\Pi(\mathbf{r}+\mathbf{q}_1,\zeta)\in H_2$, 
we conclude that the coefficient of $R_{\pi}\mathbf{G}\mathbf{r}$ is zero, that is, in leading order,
\begin{equation}\label{e:tildeT}
\widetilde{T}_0=-\sqrt[4]{\alpha_0}\frac{\pi}{\omega^2}\frac{\partial \omega}{\partial \sqrt{\kappa}},\quad
 \widetilde{T}_j=0, \quad j=1,2,3.
\end{equation}

Expressing the expansion of $\widetilde{\Pi}$ in \eqref{e:poex} in terms of ${\bf q}_2$ as in \eqref{e:loor}, we have
\begin{equation}
 \begin{aligned}
  q_{2,0}=&q_{1,0}+\caO\left(\sqrt{\varepsilon}(|\zeta|+\|q_1\|^2)\right),\\
  q_{2,j}=&\cosh(\lambda_j\frac{\pi}{\omega})q_{1,j}+\caO\left(\sqrt{\varepsilon}(|\zeta|+\|q_1\|^2)\right),\quad j=1,2,3,\\
  \widetilde{q}_{2,j}=&\sinh(\lambda_j\frac{\pi}{\omega})q_{1,j}+\caO\left(\sqrt{\varepsilon}(|\zeta|+\|q_1\|^2)\right),\quad j=1,2,3.
 \end{aligned}
\end{equation}
Therefore, $\mathbf{q}_2\in H_1\cap \Rg(\widetilde{\Pi})$ if and only if 
\begin{equation}\label{e:imp}
\widetilde{q}_{2,j}(q_1,\sqrt[4]{\varepsilon}, \sqrt{\kappa}, \zeta)=0, \quad j=1,2,3.
\end{equation}
Moreover, noting that, under the assumption \eqref{e:gen} and the assumption that $\varepsilon$ and $\kappa$ are sufficiently small, 
\begin{equation}\label{e:jacobian}
\begin{aligned}
&\varepsilon^{-1/2}\sinh(\lambda_1\frac{\pi}{\omega})
=\rmi\varepsilon^{-1/2}\sin\left(\frac{2\pi\sqrt{\varepsilon}r_2^2}{\omega}\right)
=2\sqrt{\alpha_0}\pi\rmi+\caO(\sqrt{\varepsilon}),\\
&\varepsilon^{-1/4}\sinh(\lambda_2\frac{\pi}{\omega})
=\varepsilon^{-1/4}\lambda_2\frac{\pi}{\omega}+\caO(\sqrt{\varepsilon})
=\sqrt[4]{\beta_0}\pi+\caO(\sqrt{\varepsilon}+\kappa);\\
&\varepsilon^{-1/4}\sinh(\lambda_3\frac{\pi}{\omega})
=\varepsilon^{-1/4}\lambda_3\frac{\pi}{\omega}+\caO(\sqrt{\varepsilon})
=\sqrt[4]{\beta_0}\pi\rmi+\caO(\sqrt{\varepsilon}+\kappa),
\end{aligned}
\end{equation}
we apply the rescalings
\[
\widetilde{q}_{2,1}=\sqrt{\varepsilon}p_{2,1}, \qquad \widetilde{q}_{2,j}=\sqrt[4]{\varepsilon}p_{2,j},  \quad j=2,3,
\]
to the system \eqref{e:imp} and have
\begin{equation}\label{e:reimp}
\begin{cases}
&p_{2,1}(q_1,\sqrt[4]{\varepsilon}, \sqrt{\kappa},\zeta)
=2\sqrt{\alpha_0}\pi\rmi q_{1,1}+\caO(\sqrt{\varepsilon}\|q_1\|+\|q_1\|^2+|\zeta|)=0,\\
&p_{2,2}(q_1,\sqrt[4]{\varepsilon}, \sqrt{\kappa},\zeta)
=\sqrt[4]{\beta_0}q_{1,2}+\caO(\sqrt{\varepsilon}\|q_1\|+\kappa\|q_1\|+\sqrt[4]{\varepsilon}\|q_1\|^2+\sqrt[4]{\varepsilon}|\zeta|)=0,\\
& p_{2,3}(q_1,\sqrt[4]{\varepsilon}, \sqrt{\kappa},\zeta)
=\sqrt[4]{\beta_0}\rmi q_{1,3}+\caO(\sqrt{\varepsilon}\|q_1\|+\kappa\|q_1\|+\sqrt[4]{\varepsilon}\|q_1\|^2+\sqrt[4]{\varepsilon}|\zeta|)=0.
\end{cases}
\end{equation}
Since the Jacobian of the rescaled system \eqref{e:reimp} with respect to $(q_{1,1},q_{1,2}, q_{1,3})$ at $(q_1,\sqrt[4]{\varepsilon}, \sqrt{\kappa},\zeta)=(0,0,0,0)$ is nonzero,
we may apply the implicit function theorem to the rescaled system \eqref{e:reimp}, determining that,
\begin{itemize}
 \item[(\rmnum{1})]for fixed small $\varepsilon\in[0, \varepsilon_0]$, $\kappa\in[0, \kappa_0]$ and $\zeta\in[\zeta_0, \zeta_0]$, 
there exists a one-parameter family of persistent reversible periodic orbits in \eqref{e:gODE}, parametrized by $q_{1,0}$. The periodic orbit is smooth with respect to $(q_{1,0}, \sqrt[4]{\varepsilon},\sqrt{\kappa}, \zeta)$. If we ignore both the cases $\varepsilon=0$ and $\kappa=0$, then the periodic orbit is smooth with respect to $(q_{1,0}, \varepsilon,\kappa,\zeta)$. In addition, we have
\begin{equation}\label{e:00}
q_{1,j}(0,\sqrt[4]{\varepsilon}, \sqrt{\kappa},0)=0, \quad j=1,2,3,
\end{equation}
due to the fact that, $p_{2,j}(0,\sqrt[4]{\varepsilon}, \sqrt{\kappa},0)=0$, for $j=1,2,3$.
 \item[(\rmnum{2})]for fixed small $\varepsilon\in[0, \min\{\varepsilon_0, \zeta_0\}]$ and $\kappa\in[0, \kappa_0]$, 
there exists a one-parameter family of persistent reversible periodic orbits in \eqref{e:ODE}, parametrized by $q_{1,0}$. The periodic orbit is smooth with respect to $(q_{1,0}, \sqrt[4]{\varepsilon},\sqrt{\kappa})$. If we ignore both of cases $\varepsilon=0$ and $\kappa=0$, then the periodic orbit is smooth with respect to $(q_{1,0}, \varepsilon,\kappa)$.
\end{itemize}

The fact that $\kappa$ is a free-parameter seems to contradict the uniqueness of the $q_{0,1}$-family; however, by its definition, $q_{0,1}$ is
effectively a shift of $\kappa$ and thus there is no contradiction. More specifically, for fixed $\zeta$ and $\varepsilon$, the uniqueness of the $q_{1,0}$-family in \eqref{e:gODE} implies that, for sufficiently small $\kappa>0,q_{1,0}\in\R$,
\begin{equation}\label{e:shift}
 \mathbf{r}(\sqrt{\varepsilon},\sqrt{\kappa})+\sum_{j=1}^3q_{1,j}(0,\sqrt[4]{\varepsilon},\sqrt{\kappa},\zeta)\mathbf{r}_j(\sqrt{\varepsilon},\sqrt{\kappa})=\mathbf{r}(\sqrt{\varepsilon},0)+q_{1,0}\mathbf{r}_0(\sqrt{\varepsilon},0) +\sum_{j=1}^3q_{1,j}(q_{1,0},\sqrt[4]{\varepsilon},0,\zeta)\mathbf{r}_j(\sqrt{\varepsilon},0).
\end{equation}
Setting $\zeta=\varepsilon$ and using the left hand side of \eqref{e:shift} as the initial condition to the system \eqref{e:ODE}, the initial value problem 
\begin{equation}\label{e:IVP}
\begin{cases}
& \dot{\mathbf{C}}=\mathscr{L}(\varepsilon)\mathbf{C}+\varepsilon\kappa\mathscr{R}_3(\mathbf{C})
+\varepsilon\caO\big(\varepsilon\|\mathbf{C}\|+\sqrt{\varepsilon\kappa}\|\mathbf{C}\|^2+\sqrt{\varepsilon\kappa^3}\|\mathbf{C}\|^4\big),\\
& \mathbf{C}(0)={\bf r}(\sqrt{\varepsilon},\sqrt{\kappa})+\sum_{j=1}^3q_{1,j}(0, \sqrt[4]{\varepsilon},\sqrt{\kappa},\varepsilon){\bf r}_j(\sqrt{\varepsilon},\sqrt{\kappa}),
\end{cases}
\end{equation}
admits a periodic solution, denoted as ${\bf C}^\mathrm{rp}$,  with the period 
\[
T_\mathrm{rp}(\sqrt[4]{\varepsilon},\sqrt{\kappa})=2\widetilde{T}(0,q_{1,1}(0, \sqrt[4]{\varepsilon},\sqrt{\kappa},\varepsilon),q_{1,2}(0, \sqrt[4]{\varepsilon},\sqrt{\kappa},\varepsilon),q_{1,3}(0, \sqrt[4]{\varepsilon},\sqrt{\kappa},\varepsilon),\sqrt[4]{\varepsilon},\sqrt{\kappa},\varepsilon).
\]
According to \eqref{e:Tep}, \eqref{e:tildeT} and \eqref{e:00}, we have the estimate
\[
T_\mathrm{rp}(\sqrt[4]{\varepsilon},\sqrt{\kappa})=\frac{2\pi}{\omega}+\caO\left(\sqrt{\varepsilon}(\varepsilon+\|q_1\|^2)\right)=
\frac{2\pi}{\omega}+\caO\left(\sqrt{\varepsilon^3}\right).
\]
Using the transformation 
\[
 \mathbf{C}=R_{\omega t}(\mathbf{r}+\mathbf{q}),
\]
the initial value problem \eqref{e:IVP} becomes
\begin{equation}\label{e:IVP1}
\begin{cases}
&  \dot{\mathbf{q}}=\mathbf{H}\mathbf{q}
 +\caO\big(\varepsilon\|\mathbf{q}\|^2
  +\varepsilon^{1/2}|\zeta|\big),\\
& \mathbf{q}(0)=\sum_{j=1}^3q_{1,j}(0, \sqrt[4]{\varepsilon},\sqrt{\kappa},\varepsilon){\bf r}_j(\sqrt{\varepsilon},\sqrt{\kappa}),
\end{cases}
\end{equation}
which admits a bounded solution $\|{\bf q}(t)\|_{\infty}=\caO(\varepsilon)$.

We summarize the results above in the following lemma.
\begin{Lemma}\label{l:pODE}
 For fixed $\varepsilon\in[0,\varepsilon_0]$, up to translation,
 the rescaled normal form ODE system \eqref{e:ODE}, 
 \begin{equation*}
 \begin{aligned}
 \dot{\mathbf{C}}=\mathscr{L}(\varepsilon)\mathbf{C}+\varepsilon\kappa\mathscr{R}_3(\mathbf{C})
 +\varepsilon\caO\big(\varepsilon\|\mathbf{C}\|+\sqrt{\varepsilon\kappa}\|\mathbf{C}\|^2+\sqrt{\varepsilon\kappa^3}\|\mathbf{C}\|^4\big),
 \end{aligned}
\end{equation*}
 admits a one-parameter family of persistent reversible periodic orbits, $\mathbf{C}^\mathrm{rp}(t;\sqrt[4]{\varepsilon}, \sqrt{\kappa})$, parametrized by $\kappa\in[0,\kappa_0]$. 
The periodic orbit $\mathbf{C}^\mathrm{rp}$ is smooth with respect to all parameters $(t;\sqrt[4]{\varepsilon}, \sqrt{\kappa})$. When neither $\varepsilon=0$ nor $\kappa=0$, then $\mathbf{C}^\mathrm{rp}$ is smooth with respect to $(\varepsilon,\kappa)$ and admits the form
\begin{equation}
{\bf C}^\mathrm{rp}(t;\sqrt[4]{\varepsilon},\sqrt{\kappa})=R_{\omega t}{\bf r}(\sqrt{\varepsilon},\sqrt{\kappa})+\caO\left(\varepsilon\right),
\end{equation}
where the error is measured in the $L^\infty$ norm.
The period of $\mathbf{C}^\mathrm{rp}$, denoted by $T_\mathrm{rp}$, admits the expansion
\begin{equation}
 T_\mathrm{rp}(\sqrt[4]{\varepsilon}, \sqrt{\kappa})=
 \frac{2\pi}{\omega}+\caO\left(\sqrt{\varepsilon^3}\right).
\end{equation}
\end{Lemma}

\begin{Remark}
We can also prove this lemma by using the right hand side of \eqref{e:shift} as the initial condition to \eqref{e:ODE}. But we then have to take a detour to find out the expressions of each $q_{1,j}$, $j=0,1,2,3$, in terms of $(\sqrt[4]{\varepsilon}, \sqrt{\kappa}, \zeta)$. In fact, a direct calculation using \eqref{e:shift} shows that 
\begin{equation*}
\begin{aligned}
 q_{1,0}(\sqrt[4]{\varepsilon},\sqrt{\kappa},\zeta) 
 &=\frac{\sqrt{\kappa}}{2}(r_1+\frac{r_2}{\sqrt{\alpha_0}})+q_{1,1}(0,\sqrt[4]{\varepsilon},\sqrt{\kappa},\zeta)(1-\frac{r_2^2}{\sqrt{\alpha_0}})=\frac{\sqrt{\kappa}}{\sqrt[4]{\alpha_0}}+\caO\left(\varepsilon\kappa^{5/2}+\sqrt{\varepsilon}\kappa|\zeta|\right),\\
q_{1,1}(q_{1,0},\sqrt[4]{\varepsilon},0,\zeta)&=\frac{\sqrt{\kappa}}{4}(r_1-\frac{r_2}{\sqrt{\alpha_0}})+\frac{1}{2}q_{1,1}(0,\sqrt[4]{\varepsilon},\sqrt{\kappa},\zeta)(1+\frac{r_2^2}{\sqrt{\alpha_0}})=\frac{\alpha_2\sqrt{\varepsilon\kappa^3}}{2\sqrt[4]{\alpha_0^5}}+\caO\left(\varepsilon\kappa^{5/2}+|\zeta|\right),\\
q_{1,2}(q_{1,0},\sqrt[4]{\varepsilon},0,\zeta)&=\frac{\sum_{j=2}^3q_{1,j}(0,\sqrt[4]{\varepsilon},\sqrt{\kappa},\zeta)(\lambda_j^2(\sqrt{\varepsilon},\sqrt{\kappa})-\lambda_3^2(\sqrt{\varepsilon},0))}{\lambda_2^2(\sqrt{\varepsilon},0)-\lambda_3^2(\sqrt{\varepsilon},0)}=\caO\left(|\zeta|(1+\sqrt[4]{\varepsilon}+\sqrt{\kappa})\right),\\
q_{1,3}(q_{1,0},\sqrt[4]{\varepsilon},0,\zeta)&=\frac{\sum_{j=2}^3q_{1,j}(0,\sqrt[4]{\varepsilon},\sqrt{\kappa},\zeta)(\lambda_j^2(\sqrt{\varepsilon},\sqrt{\kappa})-\lambda_2^2(\sqrt{\varepsilon},0))}{\lambda_3^2(\sqrt{\varepsilon},0)-\lambda_2^2(\sqrt{\varepsilon},0)}=\caO\left(|\zeta|(1+\sqrt[4]{\varepsilon}+\sqrt{\kappa})\right).\\
\end{aligned}
\end{equation*}
Therefore, in the system \eqref{e:ODE}, by setting $\zeta=\varepsilon$, we obtain that
\begin{equation*}
\begin{aligned}
 q_{1,0}(\sqrt[4]{\varepsilon}, \sqrt{\kappa})&:=q_{1,0}(\sqrt[4]{\varepsilon}, \sqrt{\kappa},\varepsilon)=\frac{\sqrt{\kappa}}{\sqrt[4]{\alpha_0}}+\caO\left(\varepsilon\kappa(\sqrt{\kappa^3}+\sqrt{\varepsilon})\right),\\
 q_{1,1}(\sqrt[4]{\varepsilon}, \sqrt{\kappa})&:=q_{1,1}(q_{1,0}(\sqrt{\varepsilon},\sqrt{\kappa},\varepsilon),\sqrt[4]{\varepsilon}, 0,\varepsilon)=\caO\left(\sqrt{\varepsilon}(\sqrt{\kappa^3}+\sqrt{\varepsilon})\right),\\
  q_{1,j}(\sqrt[4]{\varepsilon},\sqrt{\kappa})&:=q_{1,j}(q_{1,0}(\sqrt{\varepsilon},\sqrt{\kappa},\varepsilon),\sqrt[4]{\varepsilon}, 0,\varepsilon)=\caO(\varepsilon+\varepsilon\sqrt{\kappa}), \quad j=2,3.
\end{aligned}
\end{equation*}
\end{Remark}

Summarizing the results, we can now prove the main theorem--Theorem \ref{t:main}.
\begin{Proof}[Proof of Theorem \ref{t:main}]
The periodic solution ${\bf C}^{\mathrm{rp}}(t;\sqrt[4]{\varepsilon},\sqrt{\kappa})$ of the system \eqref{e:ODE} corresponds to a periodic solution $u_{\mathrm{rp}}(t,r; \sqrt[4]{\varepsilon},\sqrt{\kappa})$ of the PDE \eqref{e:EPP},
\begin{equation*}
 \left(\partial_r^2-W^{\prime\prime}(u)+\lambda_0\partial_t^2+\varepsilon\eta_1\right)
 \left(\partial_r^2u-W^\prime(u)+\lambda_0\partial_t^2u\right)+\varepsilon\eta_d W^\prime(u)-\varepsilon\gamma=0.
\end{equation*}
In fact, based on the center manifold reduction, the normal form transformation and the rescalings, especially Lemma \ref{l:pODE}, we have
\begin{equation}
\begin{aligned}
u_\mathrm{rp}(t, r)=&u_h(r)+\left[(A_1(t)+\bar{A}_1(t))+\rmi(A_2(t)-\bar{A}_2(t))\right]\psi_0(r)+\\
&B_1(t)\psi_1(r)+\caO(\sqrt{\varepsilon}\|{\bf A}\|+\|{\bf A}\|^2)\\
=&u_h(r)+\sqrt{\varepsilon\kappa}\left[(C_1^\mathrm{rp}(t)+\bar{C}_1^\mathrm{rp}(t))+\rmi(C_2^\mathrm{rp}(t)-\bar{C}_2^\mathrm{rp}(t))\right]\psi_0(r)+\\
&\sqrt{\varepsilon\kappa}D_1^\mathrm{rp}(t)\psi_1(r)+\caO(\varepsilon\sqrt{\kappa}\|{\bf C}^\mathrm{rp}\|)\\
=&u_h(r)+2\sqrt{\varepsilon\kappa}r_1\cos(\omega t)\psi_0(r)+\caO\left(\varepsilon(\sqrt{\varepsilon}+\sqrt{\kappa})\right)\\
=&u_h(r)+2\frac{\sqrt{\varepsilon\kappa}}{\sqrt[4]{\alpha_0}}\cos(\omega t)\psi_0(r)+\caO\left(\varepsilon(\sqrt{\varepsilon}+\sqrt{\kappa})\right),\\
\end{aligned}
\end{equation}
where we have the expression of $\omega$ from Lemma \ref{l:29}, that is,
\[
\omega=1+\omega_1\varepsilon+\sqrt{\varepsilon}r_2^2+\alpha_7\varepsilon\kappa  r_1^2+2\alpha_8\varepsilon^{3/2}\kappa.
\]
Moreover, the period of $u_{\mathrm{rp}}$, denoted by $T_{\mathrm{rp}}$, admits the expansion
\[
T_\mathrm{rp}(\sqrt[4]{\varepsilon}, \sqrt{\kappa})=
 \frac{2\pi}{\omega}+\caO\left(\sqrt{\varepsilon^3}\right).
\]
Furthermore, since the PDE \eqref{e:EPP} is a rescaled version of the stationary FCH \eqref{e:2sFCH} with the rescaling
$t=\frac{\sqrt{\lambda_0}}{\varepsilon}\tau$, the periodic solution $u_{\mathrm{rp}}$ of the PDE \eqref{e:EPP} corresponds to 
a periodic solution of the PDE \eqref{e:2sFCH}, denoted as $u_{\mathrm{p}}$ with a period $T_{\mathrm{p}}$. In fact,
\begin{equation}
\begin{aligned}
T_{\mathrm{p}}(\sqrt[4]{\varepsilon},\sqrt{\kappa})&=\frac{\varepsilon}{\sqrt{\lambda_0}}T_{\mathrm{rp}}(\sqrt[4]{\varepsilon},\sqrt{\kappa})=\frac{2\pi\varepsilon}{\sqrt{\lambda_0}}\left[1-\sqrt{\alpha_0\varepsilon}+\caO\left(\varepsilon(1+\sqrt{\kappa})\right)\right],\\
u_\mathrm{p}(\tau, r;\sqrt[4]{\varepsilon},\sqrt{\kappa})&=u_\mathrm{rp}(\frac{\sqrt{\lambda_0}}{\varepsilon}\tau, r;\sqrt[4]{\varepsilon},\sqrt{\kappa})
=u_h(r)+2\frac{\sqrt{\varepsilon\kappa}}{\sqrt[4]{\alpha_0}}\cos(\frac{2\pi}{T_{\mathrm{p}}} \tau)\psi_0(r)+\caO\left(\varepsilon(\sqrt{\varepsilon}+\sqrt{\kappa})\right),
\end{aligned}
\end{equation}
which concludes the proof of Theorem \ref{t:main}.
\end{Proof}

\section{Pearling of the Circular Planar Bilayer}\label{s:3}
In this section we consider the  case in which the bilayer sinterface $\Gamma_{R_0}$ is a circle in $\R^2$,
and construct the extended pearled solutions to the extended stationary strong FCH equation \eqref{e:2sFCH2} in $(r, \theta)\in\R\times\R/2\pi\Z$,   
\begin{equation*}
\Big(\partial_r^2-W^{\prime\prime}(u)+\frac{\varepsilon\partial_r}{R_0+\varepsilon r}+\frac{\varepsilon^2\partial_\theta^2}{(R_0+\varepsilon r)^2}+\varepsilon\eta_1\Big)
 \Big(\partial_r^2u-W^\prime(u)+\frac{\varepsilon\partial_ru}{R_0+\varepsilon r}+\frac{\varepsilon^2\partial_\theta^2u}{(R_0+\varepsilon r)^2}\Big)+\varepsilon\eta_d W^\prime(u)=\varepsilon\gamma.
\end{equation*}
 To exploit the analysis in the Section~\ref{s:2}, we rescale $\theta$ by $\vartheta=\frac{R_0\sqrt{\lambda_0}}{\varepsilon}\theta$ and search for extended pearled solutions 
$u_\mathrm{rp}$ of
\begin{equation}
 \label{e:EPP2}
 \big(\partial_r^2-W^{\prime\prime}(u)+\frac{\varepsilon\partial_r}{R_0+\varepsilon r}+\frac{R_0^2\lambda_0\partial_\vartheta^2}{(R_0+\varepsilon r)^2}+\varepsilon\eta_1\big)
 \big(\partial_r^2u-W^\prime(u)+\frac{\varepsilon\partial_ru}{R_0+\varepsilon r}+\frac{R_0^2\lambda_0\partial_\vartheta^2u}{(R_0+\varepsilon r)^2}\big)+\varepsilon\eta_d W^\prime(u)=\varepsilon\gamma.
\end{equation}
which satisfy the boundary conditions at infinity, 
\begin{equation}
 \lim_{r\rightarrow \pm\infty}|u_{\mathrm{rp}}(\vartheta,r)-u_\infty|=0, \text{ for all }\vartheta\in\R,
\end{equation}
and an even and periodic in $\vartheta$,
\begin{equation}\label{e:pc2}
 u_\mathrm{rp}(-\vartheta,r)=u_\mathrm{rp}(\vartheta,r), \quad u_\mathrm{rp}(\vartheta+T_{\mathrm{rp}},r)=u_\mathrm{rp}(\vartheta,r), \text{ for all }(\vartheta,r)\in\R^2,
\end{equation}
where $T_{\mathrm{rp}}$ and $u_\infty$ are constants to be determined.

We first prove the following proposition, which is similar to the Theorem \ref{t:main}.
\begin{Proposition}\label{p:pcb}
 Fix $\eta_1, \eta_2\in\R$ and $R_0>0$. Assume that $W$ is a non-degenerate double well potential and $\alpha_0>0$, $\beta_0\neq 0$. Then there exist positive constants $\varepsilon_0>0$ and $\kappa_0>0$ such that, for any $\varepsilon\in(0,\varepsilon_0]$, up to translation, the extended stationary FCH \eqref{e:2sFCH2} in the plane $(\theta,r)\in\R^2$,
\begin{equation*}
 \Big(\partial_r^2-W^{\prime\prime}(u)+\frac{\varepsilon\partial_r}{R_0+\varepsilon r}+\frac{\varepsilon^2\partial_\theta^2}{(R_0+\varepsilon r)^2}+\varepsilon\eta_1\Big)
 \Big(\partial_r^2u-W^\prime(u)+\frac{\varepsilon\partial_ru}{R_0+\varepsilon r}+\frac{\varepsilon^2\partial_\theta^2u}{(R_0+\varepsilon r)^2}\Big)+\varepsilon\eta_d W^\prime(u)=\varepsilon\gamma.
\end{equation*}
 admits a smooth one-parameter family of extended pearled solutions, 
 $u_\mathrm{p}(\theta,r;\sqrt[4]{\varepsilon},\sqrt{|\kappa|})$ with period $T_\mathrm{p}(\sqrt[4]{\varepsilon},\sqrt{|\kappa|})$, parameterized by $\kappa\in[-\kappa_0,\kappa_0]$. In fact, $u_{\mathrm{p}}$ and $T_{\mathrm{p}}$ are smooth with respect to their arguments within the domains expect at $\kappa=0$.  The extended pearled solution  $u_\mathrm{p}$ admits the asymptotic form
 \begin{equation}
 u_\mathrm{p}(\theta, r;\sqrt[4]{\varepsilon},\sqrt{|\kappa|})=u_h(r)+2\frac{\sqrt{\varepsilon\kappa}}{\sqrt[4]{\alpha_0}}\cos(\frac{2\pi}{T_{\mathrm{p}}} \theta)\psi_0(r)+\caO\left(\varepsilon(\sqrt{\varepsilon}+\sqrt{\kappa})\right),
 \end{equation}
where
\begin{equation}
T_{\mathrm{p}}(\sqrt[4]{\varepsilon},\sqrt{\kappa})=\frac{2\pi\varepsilon}{R_0\sqrt{\lambda_0}}\left[1-\sqrt{\alpha_0\varepsilon}+\caO\left(\varepsilon(1+\sqrt{\kappa})\right)\right].
\end{equation}
 The far-field limit of the extended pearled solution is 
 \begin{equation}
  \lim_{r\rightarrow \infty}u_\mathrm{p}(\theta,r)=\lim_{r\rightarrow \infty}u_h(r)=u_-(\varepsilon).
  \end{equation}
Moreover, for any $\varepsilon\in(0,\varepsilon_0]$, the extended stationary FCH \eqref{e:2sFCH2} in the infinite periodic strip $(\theta,r)\in\left(\R/2\pi\Z\right)\times \R$, admits a discrete family of extended pearled solutions,
 $u_\mathrm{p}(\theta,r;\sqrt[4]{\varepsilon},\sqrt{|\kappa_j|})$with period $T_\mathrm{p}(\sqrt[4]{\varepsilon},\sqrt{|\kappa_j|})$, where 
\[
\kappa_j\in\{\kappa\in[-\kappa_0,\kappa_0]\backslash\{0\}\mid 
\frac{2\pi}{T_{\mathrm{p}}(\sqrt[4]{\varepsilon},\sqrt{\kappa})}\in\Z^+\}.
\] 
 \end{Proposition} 
\begin{Proof}
The analysis of the circular interface system \eqref{e:EPP2} differs from that of the interface flat system \eqref{e:EPP} 
in two major points:
\begin{itemize}
\item[(\rmnum{1})] The circular system \eqref{e:EPP2} has different linear terms in $\varepsilon$ than the flat system \eqref{e:EPP}.
\item[(\rmnum{2})] The $S_2$ symmetry does not hold for the extended circular bilayers as it does for the flat case.
\end{itemize}

These differences only require that we  recompute  the versal normal form. More specifically, we replace $u$ with $u_h+\delta u$ in \eqref{e:EPP2} and consider the equation of the perturbation $\delta u$ (again repurposing ``$u$'' to denote the perturbation).
\begin{equation}
 \label{e:EPU2}
 \widetilde{\mathcal{L}}u+\widetilde{\mathcal{F}}(u)=0,
\end{equation}
where
\begin{equation}
 \widetilde{\mathcal{L}}:=\left(\widetilde{\caL}_h+\frac{R_0^2\lambda_0\partial_\vartheta^2}{(R_0+\varepsilon r)^2}+\varepsilon\eta_1\right)\left(\widetilde{\caL}_h+\frac{R_0^2\lambda_0\partial_\vartheta^2}{(R_0+\varepsilon r)^2}\right)+\widetilde{\mathcal{M}},
\end{equation}
with $\widetilde{\caL}_h:=\partial_r^2-W^{\pprime}(u_h)+\frac{\varepsilon\partial_r}{R_0+\varepsilon r}$, 
$\widetilde{\mathcal{M}}:=\varepsilon\eta_d W^{\pprime}(u_h)-\left(\partial_r^2u_h-W^\prime(u_h)+\frac{\varepsilon\partial_ru_h}{R_0+\varepsilon r}\right)W^{\tprime}(u_h)$, and
\begin{equation}
\label{e:F2}
\begin{aligned}
\widetilde{ \mathcal{F}}(u,\varepsilon):=&-\frac{R_0^2\lambda_0}{(R_0+\varepsilon r)^2}W^{\tprime}(u_h+u)\left(\partial_\vartheta u\right)^2-2\frac{R_0^2\lambda_0}{(R_0+\varepsilon r)^2}\left(W^{\prime\prime}(u_h+u)-W^{\prime\prime}(u_h)\right)\partial_\vartheta^2u-\\
 &\left[\caL_h+\varepsilon(\eta_1-\eta_d)-\left(W^{\pprime}(u_h+u)-W^{\pprime}(u_h)\right)\right]\left(W^\prime(u_h+u)-W^\prime(u_h)-W^{\pprime}(u_h)u\right)-\\
 & \left(\partial_r^2u_h-W^\prime(u_h)+\frac{\varepsilon\partial_ru_h}{R_0+\varepsilon r}\right)\left(W^{\pprime}(u_h+u)-W^\pprime(u_h)-W^\tprime(u_h)u\right)-\\
 &\left(W^{\pprime}(u_h+u)-W^{\pprime}(u_h)\right)\caL_h u.
\end{aligned}
\end{equation}
we recast the system as 
\begin{equation}
 \label{e:2ddFCH2}
 \dot{U}=\widetilde{\bbL}(\varepsilon)U+\widetilde{\bbF}(U,\varepsilon),
\end{equation}
where
\[
\widetilde{\bbL}(\varepsilon) = \begin{pmatrix}
0&1&0&0\\ -\frac{(R_0+\varepsilon r)^2}{R_0^2\lambda_0}\widetilde{\caL}_h&0&\frac{(R_0+\varepsilon r)^2}{R_0^2\lambda_0}&0\\ 0&0&0&1\\
-\frac{(R_0+\varepsilon r)^2}{R_0^2\lambda_0}\widetilde{\mathcal{M}}& 0&
-\frac{(R_0+\varepsilon r)^2}{R_0^2\lambda_0}(\widetilde{\caL}_h+\varepsilon\eta_1)&0\end{pmatrix},\quad
\widetilde{\bbF}(U,\varepsilon) = \begin{pmatrix}0\\0\\0\\-\frac{(R_0+\varepsilon r)^2}{R_0^2\lambda_0}\widetilde{\mathcal{F}}\end{pmatrix}.
\]
We then have
\[
\frac{\partial \widetilde{\mathbb{L}}}{\partial \varepsilon}(0)=\frac{1}{\lambda_0}
\begin{pmatrix}
 0 & 0 & 0 & 0 \\
-\frac{2r}{R_0}\caL_0 -\frac{1}{R_0}\partial_r+W^\tprime(u_0)u_1 & 0 & \frac{2r}{R_0} & 0\\
  0 & 0 & 0 & 0 \\
 -\eta_dW^\pprime(u_0)+W^\tprime(u_0)(\caL_0u_1+\frac{\partial_r u_0}{R_0}) & 0 & -\frac{2r}{R_0}\caL_0 -\frac{1}{R_0}\partial_r+W^\tprime(u_0)u_1-\eta_1 & 0
\end{pmatrix},
\]
which, after direct computation, leads to that the linear part in the reduced system in terms of ${\bf A}$, denoted as $ \widetilde{{\bf L}}(\varepsilon)$ just like its counterpart in \eqref{e:2dcredA}, is of a more complicated form
\begin{equation*}
\widetilde{\mathbf{L}}(\varepsilon)=\begin{pmatrix}
\rmi(1+\mu_1\varepsilon) &  1-\mu_1\varepsilon  & \rmi \mu_1\varepsilon  & \mu_1\varepsilon
& \rmi \widetilde{\mu}_1\varepsilon & 0 & \rmi \widetilde{\mu}_2\varepsilon &0 \\
\mu_2\varepsilon& \rmi \left(1+\mu_3\varepsilon \right) &  \mu_2\varepsilon & -\rmi \mu_3\varepsilon 
& \widetilde{\mu}_3\varepsilon & 0 & \widetilde{\mu}_4\varepsilon &0 \\
-\rmi \mu_1\varepsilon & \mu_1\varepsilon & -\rmi\left(1+\mu_1\varepsilon\right) &  1-\mu_1\varepsilon  
& -\rmi \widetilde{\mu}_1\varepsilon & 0 & -\rmi \widetilde{\mu}_2\varepsilon &0 \\
\mu_2\varepsilon & \rmi \mu_3\varepsilon &\mu_2\varepsilon & -\rmi \left(1+\mu_3\varepsilon \right) 
& \widetilde{\mu}_3\varepsilon & 0 & \widetilde{\mu}_4\varepsilon &0 \\
0 & 0 & 0 &0 & 0 & 1 & 0 & 0\\
 \widetilde{\mu}_5\varepsilon & \rmi \widetilde{\mu}_6\varepsilon & \widetilde{\mu}_5\varepsilon & -\rmi\widetilde{\mu}_6\varepsilon & \mu_4\varepsilon & 0 & 1 & 0\\
0 & 0 & 0 &0 & 0 &0 & 0 & 1\\
 \widetilde{\mu}_7\varepsilon & \rmi \widetilde{\mu}_8\varepsilon & \widetilde{\mu}_7\varepsilon & -\rmi\widetilde{\mu}_8\varepsilon & \mu_5\varepsilon & 0 & \mu_6\varepsilon & 0\\
\end{pmatrix},
\end{equation*}
Nevertheless, up to linear terms in $\varepsilon$, there exists a versal normal form of $\widetilde{{\bf L}}(\varepsilon)$ preserving the reversibility $S_1$, which takes the exact expression as its counterpart \eqref{e:ln2dcred} in the flat case, 
that is,
\begin{equation}
 \label{e:ln2dcred2}
 \mathscr{L}(\varepsilon)=\begin{pmatrix}
  \rmi\left(1+\omega_1\varepsilon\right) & 1 & 0 & 0 & 0 & 0 & 0 & 0 \\
  \omega_2\varepsilon& \rmi\left(1+\omega_1\varepsilon\right) & 0 & 0 & 0 & 0 & 0 & 0 \\
  0 & 0 & -\rmi\left(1+\omega_1\varepsilon\right) & 1 & 0 & 0 & 0 & 0 \\
  0 & 0 & \omega_2\varepsilon & -\rmi\left(1+\omega_1\varepsilon\right) & 0 & 0 & 0 & 0 \\
  0 & 0 & 0 & 0 & 0 & 1 & 0 & 0 \\
  0 & 0 & 0 & 0 & 0 & 0 & 1 & 0 \\
  0 & 0 & 0 & 0 & 0 & 0 & 0 & 1 \\
  0 & 0 & 0 & 0 & \omega_3\varepsilon & 0 & \omega_4\varepsilon & 0 \\
 \end{pmatrix},
\end{equation}
where $\omega_1=\frac{1}{2}(\mu_1+\mu_3)$, $\omega_2=\mu_2$, $\omega_3=\mu_5$, $\omega_4=\mu_4+\mu_6$. The rest of the proof is the same as the flat case.
\end{Proof}

Theorem \ref{t:main2} is drived from Proposition \ref{p:pcb} by rescaling and inverting the relation between the radius and $\kappa.$
\begin{Proof}[Proof of the Theorem \ref{t:main2}]
The stationary FCH \eqref{e:2sFCH2} on the extended plane $(\theta,r)\in\R^2$ admits a pearled solution for any $R_0\in[R_-, \infty]$, $\varepsilon\in(0,\varepsilon_0]$ and $\kappa\in[-\kappa_0,\kappa_0]$. On the other hand, the stationary FCH \eqref{e:2sFCH2} on the infinite strip $(\theta,r)\in(\R/2\pi\Z)\times\R$ requires that 
\[
T_\mathrm{p}=\frac{2\pi}{\textrm{n}}, \text{ for some }\textrm{n}\in\Z^+.
\]
Therefore, we have
\[
\textrm{n}=\frac{R_0\sqrt{\lambda_0}}{\varepsilon}\left[1+\sqrt{\alpha_0\varepsilon}+\caO\left(\varepsilon(1+\sqrt{\kappa})\right)\right],
\]
which indicates that there exists $\un>0$ so that $\textrm{n}\in[\frac{\un}{\varepsilon},\infty)$.
\end{Proof}

\section{Appendix}
We perform the calculations omitted in Section \ref{ss:22} and Section \ref{ss:23}. We begin by computing the leading order terms of the reduced 
8th-order ODE system in Appendix \ref{ss:31}. The calculation of explicit expressions for $\alpha_1$ and $\alpha_2$ follows in Appendix \ref{ss:32}.

\subsection{The reduced-ODE system in terms of ${\bf A}$}\label{ss:31}
\begin{Lemma}\label{l:reducedA}
The reduced system \eqref{e:2dred}, 
\begin{equation*}
 \dot{U}_c=\bbL_*U_c+\mathbb{P}_c\Big(\mathbb{M}(\varepsilon)\big(U_c+\Psi(U_c,\varepsilon)\big)+\bbF\big(U_c+\Psi(U_c,\varepsilon),\varepsilon,\big)\Big),
\end{equation*}
in terms of ${\bf A}:=(A_1, A_2,\bar{A}_1,\bar{A}_2, B_1, B_2, B_3, B_4)$, admits the expression 
\begin{equation*}
 \dot{\mathbf{A}}=\mathbf{L}(\varepsilon)\mathbf{A}+\mathbf{R}_2(\mathbf{A},\mathbf{A})
 +\mathbf{R}_3(\mathbf{A},\mathbf{A},\mathbf{A})+ \caO\left(|\varepsilon|^2\|{\bf A}\|+|\varepsilon|\|{\bf A}\|^2+\|{\bf A}\|^4\right),
\end{equation*}
where the linear term ${\bf L}$, the quadratic term ${\bf R}_2$, the cubic term ${\bf R}_3$ are of the following expressions.
\begin{equation*}
\mathbf{L}(\varepsilon)=\begin{pmatrix}
\rmi(1+\mu_1\varepsilon) &  1-\mu_1\varepsilon  & \rmi \mu_1\varepsilon  & \mu_1\varepsilon
& 0 & 0 & 0 &0 \\
\mu_2\varepsilon & \rmi \left(1+\mu_3\varepsilon \right) &  \mu_2\varepsilon & -\rmi \mu_3 \varepsilon
& 0 & 0 & 0 &0 \\
-\rmi \mu_1\varepsilon & \mu_1 \varepsilon& -\rmi\left(1+\mu_1\varepsilon\right) &  1+\mu_1\varepsilon
& 0 & 0 & 0 &0 \\
\mu_2\varepsilon & \rmi \mu_3\varepsilon &\mu_2 \varepsilon& -\rmi \left(1+\mu_3 \varepsilon\right) &   0 & 0 & 0 &0 \\
0 & 0 & 0 &0 & 0 & 1 & 0 & 0\\
0 & 0 & 0 &0 & \mu_4 \varepsilon& 0 & 1 & 0\\
0 & 0 & 0 &0 & 0 &0 & 0 & 1\\
0 & 0 & 0 &0 & \mu_5\varepsilon & 0 & \mu_6\varepsilon & 0\\
\end{pmatrix},
\end{equation*}
\begin{equation*}
\mathbf{R}_2(\mathbf{A},\mathbf{A})=\left(0,R_{2,2},0,\bar{R}_{2,2},0,0,0,R_{2,8}\right)^T, \quad
\mathbf{R}_3(\mathbf{A},\mathbf{A},{\bf A})=\left(0,R_{3,2},0,\bar{R}_{3,2},0,0,0,R_{3,8}\right)^T.
\end{equation*}
Here we have
\begin{equation}\label{e:linearco}
 \begin{aligned}
  &\mu_1=-\frac{1}{2\lambda_0}\int_\R W^\tprime(u_0)u_1\psi_0^2\rmd r,\quad
  \quad \mu_2=-\frac{1}{4\lambda_0^2}\int_\R\left(W^\tprime(u_0)\caL_0u_1-\eta_dW^\pprime(u_0)\right)\psi_0^2\rmd r,\\
&\mu_3=\frac{\eta_1}{2\lambda_0}-\frac{1}{4\lambda_0^2}\int_\R\left(W^\tprime(u_0)(\caL_0+2\lambda_0)u_1-\eta_dW^\pprime(u_0)\right)\psi_0^2\rmd r,\\
&\mu_4=\frac{1}{\lambda_0} \int_\R W^\tprime(u_0)u_1\psi_1^2\rmd r,
\quad \mu_5=\frac{1}{\lambda_0^2}\int_\R\left(W^\tprime(u_0)\caL_0u_1-\eta_dW^\pprime(u_0)\right)\psi_1^2\rmd r,\\
& \mu_6=-\frac{\eta_1}{\lambda_0}+\frac{1}{\lambda_0} \int_\R W^\tprime(u_0)u_1\psi_1^2\rmd r,\\
 \end{aligned}
\end{equation}
and 
\begin{equation}\label{e:Rjs}
 \begin{aligned}
 R_{2,2}=&2\nu_1\left(-a_{1,+}^2-6a_{1,+}a_{2,-}+2a_{1,-}^2+7a_{2,-}^2\right)+\nu_2(\frac{1}{2}B_1^2+B_2^2+2B_1B_3), \\
 R_{2,8}=&8\nu_2\left[\left(a_{1,+} + 3 a_{2,-}\right)B_1+2 a_{1,-}B_2 -2(a_{1,+}-a_{2,-})B_3\right],\\
 R_{3,2}=&\left(-\frac{2\nu_3}{3}+\nu_6\right)(a_{1,+}-a_{2,-})^3+2\nu_3\left[a_{1,-}^2(a_{1,+}-a_{2,-})- 2 a_{2,-}(a_{1,+}-a_{2,-})^2\right]+\\
           &\left[\frac{3}{4}\nu_7(a_{1,+}-a_{2,-})-\nu_4 a_{2, -} \right]B_1^2+\nu_4\left[-a_{1,-}B_1B_2+2(a_{1,+}-a_{2,-})(2B_1B_3+ B_2^2)\right]+\rho({\bf A}),\\
R_{3,8}=&(a_{1,+}-a_{2,-})^2\left[4\nu_4(B_1-B_3)-6\nu_7B_1\right]+8\nu_4(a_{1,+}-a_{2,-})(2a_{2,-}B_1+a_{1,-}B_2)-\\
           &4\nu_4a_{1,-}^2B_1-\left[\frac{1}{2}\nu_8B_1^3+\nu_5(B_1^2B_3+B_1B_2^2)\right]+\widetilde{\rho}({\bf A}),\\
             \end{aligned}
\end{equation} 
where
\begin{equation}\label{e:coco}
 \begin{aligned}
  &a_{1,+}=\frac{A_1+\bar{A}_1}{2},\quad a_{1,-}=\frac{A_1-\bar{A}_1}{2\rmi},\quad a_{2,+}=\frac{A_2+\bar{A}_2}{2},\quad  a_{2,-}=\frac{A_2-\bar{A}_2}{2\rmi},\\
  &\nu_1=-\frac{1}{4\lambda_0}\int_\R W^\tprime(u_0)\psi_0^3\rmd r,\quad
 \nu_2=-\frac{1}{4\lambda_0}\int_\R W^\tprime(u_0)\psi_0\psi_1^2\rmd r,\\
  &\nu_3=-\frac{1}{\lambda_0}\int_\R W^\qprime(u_0)\psi_0^4\rmd r, \quad 
  \nu_4=- \frac{1}{\lambda_0}\int_\R W^\qprime(u_0)\psi_0^2\psi_1^2\rmd r,\quad
  \nu_5=- \frac{1}{\lambda_0}\int_\R W^\qprime(u_0)\psi_1^4\rmd r,\\
  &\nu_6=\frac{1}{\lambda_0^2}\int_\R \big(W^\tprime(u_0)\big)^2\psi_0^4\rmd r,\quad
  \nu_7=\frac{1}{\lambda_0^2}\int_\R \big(W^\tprime(u_0)\big)^2\psi_0^2\psi_1^2\rmd r,\quad
  \nu_8=\frac{1}{\lambda_0^2}\int_\R \big(W^\tprime(u_0)\big)^2\psi_1^4\rmd r,\\
  &\rho({\bf A})=\int_\R Z({\bf A})\cdot ({\bf A}^T{\bf X}{\bf A})\rmd r,\quad
  \widetilde{\rho}({\bf A})=\int_\R \widetilde{Z}({\bf A})\cdot \left( {\bf A}^T{\bf X}{\bf A}\right)\rmd r.\\
 \end{aligned}
\end{equation}
In the last two expressions of $\rho$ and $\widetilde{\rho}$, the notation ``$\cdot$" denotes the Euclidean inner product in $\R^4$ and the expression of $X$ is as shown in \eqref{e:Xen} . Moreover, $Z({\bf A})$ and $\widetilde{Z}({\bf A})$ admits the forms of
\begin{equation}\label{e:Z}
\begin{aligned}
Z({\bf A})=&\frac{1}{2\lambda_0^2}W^\tprime(u_0)\psi_0^2\left[\begin{pmatrix}\caL_0\\0\\-2\\0\end{pmatrix}a_{1,+}+\begin{pmatrix}0\\2\lambda_0\\0\\0\end{pmatrix}a_{1,-}+\begin{pmatrix}4\lambda_0-\caL_0\\0\\2\\0\end{pmatrix}a_{2,-}\right]+\\
&\frac{1}{4\lambda_0^2}W^\tprime(u_0)\psi_0\psi_1\left[\begin{pmatrix}\caL_0-\lambda_0\\0\\-2\\0\end{pmatrix}B_1+\begin{pmatrix}0\\-2\lambda_0\\0\\0\end{pmatrix}B_2+\begin{pmatrix}-2\lambda_0\\0\\0\\0\end{pmatrix}B_3\right],\\
\widetilde{Z}({\bf A})=&\frac{2}{\lambda_0^2}W^\tprime(u_0)\psi_0\psi_1\left[\begin{pmatrix}\caL_0-\lambda_0\\0\\2\\0\end{pmatrix}a_{1,+}+\begin{pmatrix}0\\-2\lambda_0\\0\\0\end{pmatrix}a_{1,-}+\begin{pmatrix}-3\lambda_0-\caL_0\\0\\-2\\0\end{pmatrix}a_{2,-}\right]+\\
&\frac{1}{\lambda_0^2}W^\tprime(u_0)\psi_1^2\left[\begin{pmatrix}-\caL_0\\0\\2\\0\end{pmatrix}B_1+\begin{pmatrix}0\\2\lambda_0\\0\\0\end{pmatrix}B_2+\begin{pmatrix}2\lambda_0\\0\\0\\0\end{pmatrix}B_3\right],\\
\end{aligned}
\end{equation}
\end{Lemma}

\begin{Proof}
To simplify the calculation of the leading order terms of \eqref{e:2dred} in terms of ${\bf A}$ we introduce the following notation.
For any given integer $k\in\Z^+$, Banach spaces
$\{\mathcal{X}_j\}_{j=0}^k$ and a smooth map $F:\Pi_{j=1}^k\mathcal{X}_j\rightarrow \mathcal{X}_0$, we define
\[
 F_p:=\big(\Pi_{j=1}^k(p_j)!\big)^{-1}\partial_{x_1}^{p_1}\cdot\ldots\cdot\partial_{x_k}^{p_k}F(0,\dots,0),
\]
where 
\[
 p=(p_1,\dots,p_k)\in \Z^k,\quad p_j\geq 0, \quad j=1,\dots,k.
\]
We note
\[
 \mathbb{M}(0)=0,\quad \bbF(0,\varepsilon)\equiv 0, \quad \Psi(0,\varepsilon)\equiv 0, \quad \bbF_{(1,0)}=0,\quad \Psi_{(1,0)}=0,
\]
and conclude that the reduced system, up to cubic terms of $U_c$, is of the form
\begin{equation}
\label{e:2dred1}
\begin{aligned}
 \dot{U}_c=&\bbL_*U_c+\varepsilon\mathbb{P}_c\mathbb{M}_{1} U_c+\mathbb{P}_c\mathbb{F}_{(2,0)}(U_c,U_c)+
 \mathbb{P}_c\left(2\mathbb{F}_{(2,0)}(U_c,\Psi_{(2,0)}(U_c,U_c))+\mathbb{F}_{(3,0)}(U_c,U_c,U_c)\right),
\end{aligned}
\end{equation}
with the higher order terms in the form of 
$ \caO\left(|\varepsilon|^2\|U_c\|+|\varepsilon|\|U_c\|^2+\|U_c\|^4\right)$.
A direct calculation shows that
\begin{equation}
 \label{e:coeffs}
 \begin{aligned}
 &\mathbb{M}_{1}=\frac{1}{\lambda_0}
\begin{pmatrix}
 0 & 0 & 0 & 0 \\
 W^\tprime(u_0)u_1 & 0 & 0 & 0\\
  0 & 0 & 0 & 0 \\
 -\eta_dW^\pprime(u_0)+W^\tprime(u_0)\caL_0u_1 & 0 & W^\tprime(u_0)u_1-\eta_1 & 0
\end{pmatrix},\\
& \bbF_{(2,0)}(U_c,U_c)=\frac{1}{\lambda_0}\left(
 W^\tprime(u_0)\left(2u_cv_c-u_c\caL_0u_c+\lambda_0p_c^2\right)+\frac{1}{2}\caL_0\left(W^\tprime(u_0)u_c^2\right)\right)\caE_4,\\
 & \bbF_{(3,0)}(U_c,U_c,U_c)=\frac{1}{\lambda_0}\left(
 W^\qprime(u_0)u_c\left(u_cv_c-\frac{1}{2}u_c\caL_0u_c+\lambda_0p_c^2\right)+\frac{1}{6}\caL_0\left(W^\qprime(u_0)u_c^3\right)
 -\frac{1}{2}(W^\tprime(u_0))^2u_c^3\right)\caE_4,\\
 &\mathbb{F}_{(2,0)}(U_c,\Psi_{(2,0)}(U_c,U_c))=
 \frac{1}{2\lambda_0}\left( V_c\cdot \Psi_{(2,0,0)}(U_c,U_c) \right)\caE_4,\\
 \end{aligned}
\end{equation}
where $\caE_4=(0,0,0,1)^T$,  $u_1$ comes from the Taylor expansion,
\[
 u_h(\varepsilon)=u_0+\varepsilon u_1+\caO(\varepsilon^2),
\]
and
\begin{equation}
\label{e:vc}
 V_c=\left(2W^\tprime(u_0)v_c-W^\tprime(u_0)\caL_0u_c+[\caL_0, W^\tprime(u_0)u_c],2\lambda_0W^\tprime(u_0)p_c, 2W^\tprime(u_0)u_c, 0\right)^T.
\end{equation}
We also use the notation that $U_c=(u_c,p_c,v_c,q_c)^T$, where
\begin{equation}
\label{e:uc}
\begin{array}{ll}
u_c=2(a_{1,+}-a_{2,-})\psi_0+B_1\psi_1, & p_c=-2a_{1,-}\psi_0+B_2\psi_1, \\
v_c=-4\lambda_0a_{2,-}\psi_0+\lambda_0B_3\psi_1, & q_c=-4\lambda_0a_{2,+}\psi_0+\lambda_0B_4\psi_1.
\end{array}
\end{equation}
A direct calculation, using \eqref{e:coeffs}-\eqref{e:uc} and the expression of ${\bf X}$ \eqref{e:Xen},  leads to explicit expressions of the linear part ${\bf L}$, the quadratic part ${\bf R}_2$ and the cubic term ${\bf R}_3$. Relegating the calculation of $\Psi_{(2,0)}(U_c,U_c)$ into Lemma \ref{l:31}, we conclude our proof.
\end{Proof}

\begin{Lemma}\label{l:31}
The quadratic term of the center manifold, $\Psi_{(2,0)}(U_c,U_c)$, is a quadratic form of ${\bf A}$ and thus takes the form
\begin{equation}\label{e:qcm}
\Psi_{(2,0)}(U_c,U_c)=\mathbf{A}^\mathbf{T}{\bf X}{\bf A}
\end{equation}
where ${\bf X}=\{X_{jk}\}_{j,k=1}^8$ is symmetric
and every entry $X_{jk}\in \mathbb{P}_h\mathcal{Y}$. 
More specifically, we have
 we have
\begin{equation}\label{e:Xen}
\begin{aligned}
&\begin{cases}
&X_{11}=\bar{X}_{33}=\left(2\rmi-\bbL_*\right)^{-1}Y_{1}, \\
&X_{12}=\bar{X}_{34}=\left(2\rmi-\bbL_*\right)^{-1}\left(Y_{2}-X_{11}\right), \\
&X_{22}=\bar{X}_{44}=\left(2\rmi-\bbL_*\right)^{-1}\left(Y_{4}-2X_{12}\right) ,
\end{cases}
\begin{cases}
&X_{13}=-\bbL_*^{-1}Y_{3}, \\
&X_{14}=-\bbL_*^{-1}\left(-Y_{2}-X_{13}\right), \\
&X_{23}=-\bbL_*^{-1}\left(Y_{2}-X_{13}\right) ,\\
&X_{24}=-\bbL_*^{-1}\left(-Y_{4}-X_{14}-X_{23}\right),
\end{cases}
\end{aligned}
\end{equation}
\begin{equation}
\begin{aligned}
&\begin{cases}
&X_{15}=\bar{X}_{35}=\left(\rmi-\bbL_*\right)^{-1}Y_{5}, \\
&X_{16}=\bar{X}_{36}=\left(\rmi-\bbL_*\right)^{-1}\left(\rmi Y_{7}-X_{15}\right), \\
&X_{17}=\bar{X}_{37}=\left(\rmi-\bbL_*\right)^{-1}\left(Y_{7}-X_{16}\right) ,\\
&X_{18}=\bar{X}_{38}=\left(\rmi-\bbL_*\right)^{-1}\left(-X_{17}\right), \\
&X_{25}=\bar{X}_{45}=\left(\rmi-\bbL_*\right)^{-1}\left(Y_{6}-X_{15}\right), \\
&X_{26}=\bar{X}_{46}=\left(\rmi-\bbL_*\right)^{-1}\left(-X_{16}-X_{25}\right), \\
&X_{27}=\bar{X}_{47}=\left(\rmi-\bbL_*\right)^{-1}\left(\rmi Y_{7}-X_{17}-X_{26}\right) ,\\
&X_{28}=\bar{X}_{48}=\left(\rmi-\bbL_*\right)^{-1}\left(-X_{18}-X_{27}\right),
\end{cases}\quad
\begin{cases}
&X_{55}=-\bbL_*^{-1}Y_{8}, \\
&X_{56}=-\bbL_*^{-1}\left(-X_{55}\right), \\
&X_{57}=-\bbL_*^{-1}\left(Y_{9}-X_{56}\right) ,\\
&X_{58}=-\bbL_*^{-1}\left(-X_{57}\right),\\
&X_{66}=-\bbL_*^{-1}\left(Y_{9}-2X_{56}\right),\\
&X_{67}=-\bbL_*^{-1}\left(-X_{57}-X_{66}\right),\\
&X_{68}=-\bbL_*^{-1}\left(-X_{58}-X_{67}\right),\\
&X_{77}=-\bbL_*^{-1}\left(-2X_{67}\right),\\
&X_{78}=-\bbL_*^{-1}\left(-X_{68}-X_{77}\right),\\
&X_{88}=-\bbL_*^{-1}\left(-2X_{78}\right),
\end{cases}
\end{aligned}
\end{equation}
where, introducing $\caE_4=(0,0,0,1)^T$ we have
\begin{equation}\label{e:Yen}
\begin{aligned}
&Y_1=\left(\frac{1}{2\lambda_0}\mathcal{L}_0-2\right) W^{\tprime}(u_0)\psi_0^2 \caE_4-
6\lambda_0\nu_1\psi_0 \caE_4, \\
&Y_2=\rmi\left(\frac{1}{2\lambda_0}\mathcal{L}_0+1\right) W^{\tprime}(u_0)\psi_0^2\caE_4+
6\rmi\lambda_0\nu_1 \psi_0\caE_4, \\
&Y_3=\frac{1}{2\lambda_0}\mathcal{L}_0 W^{\tprime}(u_0)\psi_0^2\caE_4+
2\lambda_0\nu_1 \psi_0\caE_4,\\
&Y_4=-\left(\frac{1}{2\lambda_0}\mathcal{L}_0+3\right)W^{\tprime}(u_0)\psi_0^2\caE_4-
14\lambda_0\nu_1 \psi_0\caE_4, \\
&Y_5=\left(\frac{1}{2\lambda_0}\mathcal{L}_0-\frac{1}{2}\right) W^{\tprime}(u_0)\psi_0\psi_1\caE_4-
2\lambda_0\nu_2 \psi_1,\caE_4 \\
&Y_6=\rmi\left(\frac{1}{2\lambda_0}\mathcal{L}_0+\frac{3}{2}\right)W^{\tprime}(u_0)\psi_0\psi_1\caE_4+
6\rmi\lambda_0\nu_2 \psi_1\caE_4, \\
&Y_7=W^{\tprime}(u_0)\psi_0\psi_1\caE_4+
4\lambda_0\nu_2 \psi_1\caE_4, \\
&Y_8=\frac{1}{2\lambda_0}\mathcal{L}_0W^{\tprime}(u_0)\psi_1^2\caE_4+
2\lambda_0\nu_2 \psi_0\caE_4, \\
&Y_9=W^{\tprime}(u_0)\psi_1^2\caE_4+
4\lambda_0\nu_2 \psi_0\caE_4.
\end{aligned}
\end{equation}
\end{Lemma}
\begin{Proof}
To find the explicit expression of $\Psi_{(2,0)}(U_c,U_c)$ in terms of ${\bf A}$, we first recall \eqref{e:2ddFCH} and \eqref{e:sch} as follows.
\begin{equation*}
\begin{aligned}
\dot{U}=\bbL(\varepsilon)U+\bbF(U,\varepsilon), \quad U=U_c+\Psi(U_c,\varepsilon).
\end{aligned}
\end{equation*}
Plugging \eqref{e:sch} into \eqref{e:2ddFCH}, applying the projection $\mathbb{P}_h:=\id-\mathbb{P}_c$ and setting $\varepsilon=0$, we obtain
\begin{equation}
 \label{e:teh}
 \dot{\Psi}(U_c,0)=\mathbb{L}_*\Psi(U_c,0)+\mathbb{P}_h\bbF\big(U_c+\Psi(U_c,0),0\big).
\end{equation}
For simplicity, we note that $\mathbb{P}_h\bbF_{(2,0)}(U_c,U_c)$ is a quadratic form of ${\bf A}$ and thus takes the form
\[
 \mathbb{P}_h\bbF_{(2,0)}(U_c,U_c)=\mathbf{A}^\mathbf{T}{\bf Y}{\bf A},
\]
where ${\bf Y}=\{Y_{jk}\}_{j,k=1}^8$ is symmetric
and every entry $Y_{jk}\in \mathbb{P}_h\mathcal{Y}$. 
Restricting \eqref{e:teh} to the quadratic terms of $U_c$ and plugging in \eqref{e:qcm}, we have, for all ${\bf A}$,
\begin{equation*}
 {\bf A}^T\left({\bf L}(0)^T{\bf X}+{\bf X}{\bf L}(0)\right){\bf A}= {\bf A}^T(\mathbb{L}_*{\bf X}){\bf A}+\mathbf{A}^\mathbf{T}{\bf Y}{\bf A},
\end{equation*}
that is,
\begin{equation}\label{e:XY}
{\bf L}(0)^T{\bf X}+{\bf X}{\bf L}(0)-\mathbb{L}_*{\bf X}={\bf Y},
\end{equation}
from which it is not hard to compute all entries of ${\bf X}$ recursively.
More explicitly, ${\bf Y}$ admits the form
\begin{equation}\label{e:Y}
\begin{pmatrix}
Y_1 & Y_2 & Y_3 & -Y_2 & Y_5 & \rmi Y_7 &Y_7 &0 \\
Y_2 & Y_4 & Y_2 & -Y_4 & Y_6 & 0 & \rmi Y_7 &0 \\
Y_3 & Y_2 & Y_1 &-Y_2 & Y_5 & -\rmi Y_7 & Y_7 &0 \\
-Y_2 & -Y_4 & -Y_2 & Y_4 & -Y_6 & 0 & -\rmi Y_7 &0 \\
Y_5 & Y_6 & Y_5 & -Y_6 & Y_{8} & 0 & Y_{9} &0 \\
\rmi Y_7 & 0 & -\rmi Y_7 &0 & 0 & Y_{9} & 0 & 0\\
Y_7 & \rmi Y_7 & Y_7 & -\rmi Y_7 &  Y_{9} & 0 & 0 & 0\\
0 & 0 & 0 &0 & 0 & 0 & 0 &0\\
\end{pmatrix},
\end{equation}
where $Y_j$'s admit the expressions as in \eqref{e:Yen}.
Plugging \eqref{e:Y} into \eqref{e:XY}, we obtain the expression of $X_{ij}$'s as in \eqref{e:Xen}. 
\end{Proof}

\subsection{Explicit expressions of $\alpha_1$ and $\alpha_2$}\label{ss:32}
\begin{Lemma}\label{l:32}
Among the coefficients of cubic terms of the normal form system \eqref{e:2dnfeq}, we have
\begin{equation}\label{e:alphas}
\begin{cases}
&\alpha_1=0,\\ 
&\alpha_2=-\frac{\nu_3}{3}+\frac{80}{9}\nu_1^2+\int_\R\left(W^{\tprime}(u_0)\psi_0^2+4\lambda_0\nu_1\psi_0\right)\widetilde{\caL}\left(W^{\tprime}(u_0)\psi_0^2+4\lambda_0\nu_1\psi_0\right)\rmd r,
\end{cases}
\end{equation}
where $\widetilde{\caL}:=\frac{1}{3\lambda_0^2}\left(\frac{1}{2}+2\lambda_0\caL_0^{-1}+2\lambda_0(\caL_0-4\lambda_0)^{-1}-\lambda_0(\caL_0-4\lambda_0)^{-2}\right)$ is a self-adjoint operator.
\end{Lemma}
\vskip -0.4in
\begin{Remark}
The techniques used in the proof of Lemma \ref{l:32} permit  the calculation of explicit expressions for each $\alpha_j$ and $\beta_k$. Nevertheless, 
we only present the calculations of $\alpha_1$ and $\alpha_2$, as the other coefficients are not needed in the sequel.
\end{Remark}
\vskip -0.4in
\begin{Proof}
To calculate all these coefficients, we first recall the equality \eqref{e:nf3} with $\mathscr{R}_2=0$,
\begin{equation*}
   \big(\mathcal{D}-\mathbf{L}(0,0)\big)\Phi_3=\mathbf{R}_3+2\mathbf{R}_2(\mathbf{C},\Phi_2)-\mathscr{R}_3,
\end{equation*}
and the restrictions
\[
\begin{pmatrix}\mathbf{R}_{3,1}+2\mathbf{R}_{2,1}(\mathbf{C},\Phi_2)-\mathscr{R}_{3,1}\\
\mathbf{R}_{3,2}+2\mathbf{R}_{2,2}(\mathbf{C},\Phi_2)-\mathscr{R}_{3,2}\end{pmatrix}\in
\left(\ker\left(\begin{pmatrix}\mathcal{D}^{\ad}+\rmi&0\\-1&\mathcal{D}^{\ad}+\rmi\end{pmatrix}\mid_{\mathbf{P}_3^2}\right)\right)^\perp,
\]
 we have that $\alpha_1$ is exactly the coefficient of $C_1^2\bar{C}_1$ in 
\[\mathbf{R}_{3,2}+2\mathbf{R}_{2,2}(\mathbf{C},\Phi_2),\]
that is,
\begin{equation}\label{e:alpha0}
\alpha_1=\frac{1}{2}\langle \mathbf{R}_{3,2}+2\mathbf{R}_{2,2}(\mathbf{C},\Phi_2)\mid C_1^2\bar{C}_1\rangle,
\end{equation}
where we recall that this inner product is the one of polynomials, defined as $\langle P \mid Q \rangle=P(\partial_{{\bf C}})\bar{Q}(C)$.
According to \eqref{e:Rjs} and \eqref{e:phi2}, we have
\begin{equation*}
\begin{aligned}
\alpha_1=\frac{1}{2}\langle \mathbf{R}_{3,2}+2\mathbf{R}_{2,2}(\mathbf{C},\Phi_2)\mid C_1^2\bar{C}_1\rangle=-6\nu_1^2+\frac{3}{8}\nu_6+\frac{1}{2}\langle\rho({\bf C})\mid C_1^2\bar{C}_1\rangle.
\end{aligned}
\end{equation*}
From the expression of $\rho({\bf A})$ in \eqref{e:coco}, it is straight forward to see that 
\begin{equation}
\label{e:rhoalpha1}
\begin{aligned}
\frac{1}{2}\langle\rho({\bf C})\mid  C_1^2\bar{C}_1\rangle&=\frac{1}{2}\int_\R \langle Z({\bf C})\cdot(X_{11}C_1^2+2X_{13}C_1\bar{C}_1)\mid C_1^2\bar{C}_1\rangle\rmd r,
\end{aligned}
\end{equation}
Based on \eqref{e:Xen} and \eqref{e:Yen}, a direct calculation shows that 
\begin{equation*}
\begin{aligned}
X_{13}=-\bbL_*^{-1}Y_3&=\frac{1}{2}\begin{pmatrix}\mathcal{L}_0^{-1}\\0\\1\\0\end{pmatrix}\left(W^{\tprime}(u_0)\psi_0^2+4\lambda_0\nu_1\psi_0\right),\\
X_{11}=-(\bbL_*-2\rmi)^{-1}Y_1&=\frac{1}{2}\begin{pmatrix}\left(\mathcal{L}_0-4\lambda_0\right)^{-1}\\2\rmi\left(\mathcal{L}_0-4\lambda_0\right)^{-1}\\1\\2\rmi\end{pmatrix}\left(W^{\tprime}(u_0)\psi_0^2+4\lambda_0\nu_1\psi_0\right).
\end{aligned}
\end{equation*}
Plugging \eqref{e:Z} and the above expressions into \eqref{e:rhoalpha1}, we have
\begin{equation*}
\begin{aligned}
\frac{1}{2}\langle\rho({\bf C})\mid C_1^2\bar{C}_1\rangle&=-\frac{3}{8}\nu_6+6\nu_1^2.
\end{aligned}
\end{equation*}
Therefore, we deduce that $\alpha_1=0.$
A similar calculation shows that
\[
\alpha_2=-\frac{\nu_3}{3}+\frac{80}{9}\nu_1^2+\int_\R\left(W^{\tprime}(u_0)\psi_0^2+4\lambda_0\nu_1\psi_0\right)\widetilde{\caL}\left(W^{\tprime}(u_0)\psi_0^2+4\lambda_0\nu_1\psi_0\right)\rmd r,
\] 
where $\widetilde{\caL}:=\frac{1}{3\lambda_0^2}\left(\frac{1}{2}+2\lambda_0\caL_0^{-1}+2\lambda_0(\caL_0-4\lambda_0)^{-1}-\lambda_0(\caL_0-4\lambda_0)^{-2}\right)$ is a self-adjoint operator.
\end{Proof}

\subsection*{Acknowledgment}
The first author acknowledges support from  NSF DMS grants 1109127 and 1409940.

\bibliographystyle{siam}
\bibliography{myref}

\begin{thebibliography}{10}

\bibitem{amphiphilic2000}
{\sc P.~Alexandridis and B.~Lindman}, {\em Amphiphilic block copolymers:
  self-assembly and applications}, Elsevier, 2000.

\bibitem{viarnold}
{\sc V.~Arnold}, {\em Matrices depending on parameters}, Uspehi Mat. Nauk, 26
  (1971), pp.~101--114.

\bibitem{szostak}
{\sc I.~Budin and J.~Szostak}, {\em Physical effects underlying the transition
  from primitive to modern cell membranes}, Proceedings of the National Academy
  of Sciences, 108 (2011), pp.~5249--5254.

\bibitem{CH}
{\sc J.~Cahn and J.~Hilliard}, {\em Free energy of a nonuniform system. i.
  interfacial free energy}, The Journal of Chemical Physics, 28 (1958),
  pp.~258--267.

\bibitem{keithdai_2013}
{\sc S.~Dai and K.~Promislow}, {\em Geometric evolution of bilayers under the
  functionalized cahn--hilliard equation}, Proceedings of the Royal Society A:
  Mathematical, Physical and Engineering Science, 469 (2013).

\bibitem{doelmanpromislow_2013}
{\sc A.~Doelman, G.~Hayrapetyan, K.~Promislow, and B.~Wetton}, {\em Meander and
  pearling of single-curvature bilayer interfaces in the functionalized
  {C}ahn-{H}illiard equation}, SIAM Journal Math. Analysis, to appear (2014).

\bibitem{KeithNiGr-11}
{\sc N.~Gavish, G.~Hayrapetyan, K.~Promislow, and L.~Yang}, {\em Curvature
  driven flow of bi-layer interfaces}, Physica D: Nonlinear Phenomena, 240
  (2011), pp.~675 -- 693.

\bibitem{glebsky_1995}
{\sc L.~Glebsky and L.~Lerman}, {\em On small stationary localized solutions
  for the generalized {$1$}-{D} {S}wift-{H}ohenberg equation}, Chaos, 5 (1995),
  pp.~424--431.

\bibitem{GS-90}
{\sc G.~Gompper and M.~Schick}, {\em Correlation between structural and
  interfacial properties of amphiphilic systems}, Phys. Rev. Lett., 65 (1990),
  pp.~1116--1119.

\bibitem{harioo}
{\sc M.~Haragus and G.~Iooss}, {\em Local bifurcations, center manifolds, and
  normal forms in infinite-dimensional dynamical systems}, Universitext,
  Springer-Verlag London Ltd., London, 2011.

\bibitem{keithgreg_2014}
{\sc G.~Hayrapetyan and K.~Promislow}, {\em Spectra of functionalized operators
  arising from hypersurfaces}, Z. Angew. Math. Phys., to appear (2014).

\bibitem{imd_1989}
{\sc G.~Iooss, A.~Mielke, and Y.~Demay}, {\em Theory of steady
  {G}inzburg-{L}andau equation, in hydrodynamic stability problems}, European
  J. Mech. B Fluids, 8 (1989), pp.~229--268.

\bibitem{iooper_1993}
{\sc G.~Iooss and M.-C. P{\'e}rou{\`e}me}, {\em Perturbed homoclinic solutions
  in reversible {$1:1$} resonance vector fields}, J. Differential Equations,
  102 (1993), pp.~62--88.

\bibitem{kappro}
{\sc T.~Kapitula and K.~Promislow}, {\em Spectral and dynamical stability of
  nonlinear waves}, vol.~185 of Applied Mathematical Sciences, Springer, New
  York, 2013.
\newblock With a foreword by Christopher K. R. T. Jones.

\bibitem{KeithBrian-09}
{\sc K.~Promislow and B.~Wetton}, {\em Pem fuel cells: A mathematical
  overview}, SIAM Journal on Applied Mathematics, 70 (2009), pp.~369--409.

\bibitem{PZ-2013}
{\sc K.~Promislow and H.~Zhang}, {\em Critical points of functionalized
  lagrangians}, Discrete and Continuous Dynamical Systems, 33 (2013),
  pp.~1--16.

\bibitem{DeGiorgi}
{\sc M.~R\"{o}ger and R.~Sch\"{a}tzle}, {\em On a modified conjecture of de
  giorgi}, Mathematische Zeitschrift, 254 (2006), pp.~675--714.

\bibitem{ssmorse_2008}
{\sc B.~Sandstede and A.~Scheel}, {\em Relative {M}orse indices, {F}redholm
  indices, and group velocities}, Discrete Contin. Dyn. Syst., 20 (2008),
  pp.~139--158.

\bibitem{TS-87}
{\sc M.~Teubner and R.~Strey}, {\em Origin of the scattering peak in
  microemulsions}, The Journal of Chemical Physics, 87 (1987), pp.~3195--3200.

\bibitem{PB_2011}
{\sc P.~van Heijster and B.~Sandstede}, {\em Planar radial spots in a
  three-component {F}itz{H}ugh-{N}agumo system}, J. Nonlinear Sci., 21 (2011),
  pp.~705--745.

\bibitem{hayward}
{\sc J.~Zhu and R.~Hayward}, {\em Wormlike micelles with microphase-separated
  cores from blends of amphiphilic ab and hydrophobic bc diblock copolymers},
  Macromolecules, 41 (2008), pp.~7794--7797.

\end{thebibliography}

\end{document}